\documentclass{ws-cm}
\usepackage{graphicx}
\usepackage{psfrag}
\usepackage{amsfonts,amssymb
}
\usepackage{amsmath}
\usepackage{float}
\usepackage{color}
\usepackage{enumitem}

\renewcommand{\vec}[1]{\mbox{\boldmath $#1$}}
\newtheorem{prop}{Proposition}[section]
\newcommand{\ds}{\displaystyle}
\newcommand{\NN}{\mathbb{N}}
\newcommand{\RR}{\mathbb{R}}
\newcommand{\CC}{\mathbb{C}}

\newcommand{\VV}{\mathbb{V}}

\newcommand{\GG}{\mathbb{G}}
\newcommand{\QQ}{\mathbb{Q}}

\newcommand{\pd} { \partial_}

\newcommand{\im}{\mathrm{Range} \,}
\renewcommand{\Im}{\mathrm{Im} \,}
\renewcommand{\ker}{\mathrm{Ker} \,}
\newcommand{\rk}{\mathrm{rank} \,}
\newcommand{\re}{\mathrm{Re} \,}
\renewcommand{\Re}{\mathrm{Re} \,}
\renewcommand{\dim}{\mathrm{dim} \,}
\newcommand{\fr}{\frac{1}{2}}
\newcommand{\txi}{(\tau,\xi)}

\newcommand{\txim}{(\tau,\xi\imj{R})}
\newcommand{\ltx}{L\txi}
\newcommand{\lutx}{L_1\txi}

\newcommand{\kl}{\ker \ltx}
\newcommand{\klu}{\ker \lutx}

\newcommand{\lo}{L_1(0,\xi)}
\newcommand{\klo}{\ker \lo}
\newcommand{\tia}{\widetilde{A}}
\newcommand{\tiu}{\widetilde{U}}
\newcommand{\til}{\widetilde{L}}
\newcommand{\tilu}{\widetilde{L}_1}
\newcommand{\tiltx}{\til\txi}
\newcommand{\tilutx}{\til_1\txi}

\newcommand{\tiklu}{\ker \tilutx}

\newcommand{\tilo}{\til_1(0,\xi)}
\newcommand{\tiklo}{\ker \tilo}
\newcommand{\imj}[1]{^{\mbox{\tiny{$#1$}}}}
\newcommand{\lmj}[1]{_{\mbox{\tiny{$#1$}}}}
\newcommand{\cl}{{\cal L}}

\newcommand{\cau}{{\cal U}}
\newcommand{\cav}{{\cal V}}
\newcommand{\caw}{{\cal W}}
\newcommand{\caB}{{\cal B}}
\newcommand{\trace}[1]{_{\,\mid#1}}

\newcommand{\ip}{(I-\Pi_{L})\,}
\newcommand{\ipr}{(I-\Pi_{L}\imj{R})\,}
\newcommand{\ipc}{(I-\Pi_\cl)\,}
\newcommand{\ipcr}{(I-\Pi_\cl\imj{R})\,}

\def\ueta{{\underline \eta}}
\def\utau{{\underline \tau}}
\def\uxi{{\underline \xi}}
\def\ux{{\underline x}}
\def\ut{{\underline t}}

\newcommand{\sig}{\sigma}

\newcommand{\eps}{\varepsilon}
\renewcommand{\phi}{\varphi}

\def\bfe{{\bf e}}

\def\bfv{{\bf v}}

\def\caS{{\cal S}}
\def\caN{{\cal N}}

\def\caA{{\cal A}}
\def\caO{{\cal O}}

\begin{document}

\markboth{L. Halpern, S. Petit-Bergez, J. Rauch}{The Analysis of Matched Layers}

%%%%%%%%%%%%%%%%%%% Publisher's Area please ignore %%%%%%%%%%%%%%%%%%%%%%
\catchline{}{}{}{}{}
%%%%%%%%%%%%%%%%%%%%%%%%%%%%%%%%%%%%%%%%%%%%%%%%%%%%%%%%%%%%%%%%%%%%%%%%%

\title{The Analysis of Matched Layers}

\author{ L. Halpern}

\address{LAGA,Institut Galil\'ee,
  Universit\'e Paris XIII, 93430 Villetaneuse, FRANCE\\
halpern@math.univ-paris13.fr}

\author{S. Petit-Bergez}

\address{LAGA,Institut Galil\'ee,
  Universit\'e Paris XIII, 93430 Villetaneuse, FRANCE\\
sabrina.bergez@laposte.net}

\author{J. Rauch }

\address{Department of Mathematics,
University of Michigan, Ann Arbor  48109  MI, USA.
\footnote{Research partially supported by the National Science
Foundation under grant  NSF DMS 0405899.}}

\maketitle

\begin{history}
\received{Day Month Year}
\revised{Day Month Year}
%\accepted{(Day Month Year)}
\end{history}

\maketitle

\begin{abstract}
A systematic analysis of matched layers is undertaken with special attention to better understand the remarkable method of B\'erenger. We prove that the
B\'erenger and closely related layers define well posed  transmission problems in great generality.
When the B\'erenger method or one of  its close relatives
is well posed, perfect matching is proved.
 The proofs use the energy method,
 Fourier-Laplace transform,
 and real coordinate changes for Laplace transformed equations.
  It is proved that the loss of derivatives associated with  the
  B\'erenger method does not occur for elliptic generators.
  More generally,  an essentially
   necessary and sufficient condition for loss of derivatives
  in B\'erenger's method is proved. The sufficiency
   relies on  the energy method with  pseudodifferential multiplier.
  Amplifying and nonamplifying layers are identified by a geometric optics computation.
  Among the various flavors of B\'erenger's algorithm for Maxwell's
   equations our favorite choice
   leads to a strongly well posed augmented system and is both
   perfect and nonamplifying in great generality.     We construct by an extrapolation argument an alternative matched layer method which preserves the strong hyperbolicity of the original problem and though not perfectly matched has {\it leading} reflection coefficient equal to zero at all angles of incidence.
Open problems are indicated throughout.\\
\end{abstract}


\vskip.3cm
\centerline{\bf Dedication}
\vskip.1cm

\noindent
{\it En \'ecrivant ce papier, les auteurs ont toujours pr\'esente
\`a l'esprit leur amiti\'e pour Michelle Schatzman. Comment la garder vivante sinon en manifestant chaque jour la curiosit\'e, l'exigence scientifique et le plaisir du partage qui \'etaient les siens.
}

\vskip.5cm

\tableofcontents
%%%%%%%%%%%%%%%%%%%%%%%%%%%
%  New section
%%%%%%%%%%%%%%%%%%%%%%%%%%
\section{Introduction}\label{sec:introduction}

This paper analyses  absorbing layer methods for calculating approximations to the solution, $U$, of  first order systems of hyperbolic partial differential equations,
\begin{equation}\label{hyp}
L(\partial _t,\partial_x)\,U
\ :=\
\partial _tU
\ +\
\sum_{l=1}^d\, A_l\, \partial _{l}U
\ =\ F\,,
\qquad
(t,x)\ \in \
\RR^{1+d}\,,
\qquad
U(t,x)\in \CC^N\,.
\end{equation}
Approximate values are sought on a finite domain. The source term $F$  and/or initial condition is compactly supported in the domain. The absorbing layer strategy surrounds the domain with a layer of finite thickness intended to be absorbing and weakly reflective.

%We lay foundations sufficient  to give a good (but not entirely complete)
% understanding of the remarkable method of
% B\'erenger, and its close relatives, which have no reflection at all from a planar % interface even in the presence of  a discontinuous coefficient. We propose a
% new extrapolation method, which we call the Harmoniously Matched Layer. It % is reflectionless to leading order and inherits   the $L^2$ estimates of
% symmetric hyperbolic systems. Numerical tests confirm that  it is competitive
% with the B\'erenger algorithms for Maxwell's equations.

The simplest case is dimension $d=1$ with computational domain
$x_1<0$ and absorbing layer in $x_1>0$.  For the  first example
consider inhomogeneous initial data and zero right hand side. The simplest absorbing layers add a lower order term $\sigma {\bf 1}_{x_1>0}\, C\, U$ where $\bf 1$ denotes the characteristic function, for example,
$$
\partial_t U +
\begin{pmatrix}
1 & 0\\
0& -1
\end{pmatrix}
\partial_1 U
 \ +\
 \sigma\,
 {\bf 1}_{x_1>0}
\begin{pmatrix}
1 & 0\\
c & b
\end{pmatrix}
U
\ =\ 0.
$$
To get a feeling for the reflections, consider the solution $U(t,x_1)$ so that,
$$
 {\rm for}\ \  t<0\,,\qquad
U\ =\
\big(
\delta(x_1-t)\,,\, 0\big)\,.
$$
Then
$$
U_1 = \delta(x_1-t)\, e^{-\sigma x_1\, {\bf 1}_{x_1>0}},\qquad
\big(\partial_t - \partial_1 + b\sigma
{\bf 1}_{x_1>0}\ \big)U_2
\ =\
 -\sigma \,
 {\bf 1}_{x_1>0}\,
c\, U_1.
$$
If $c\ne 0$, then $\nabla_{t,x_1}U_2$ is discontinuous across the ray $\{x_1=-t\}$.  From the perspective of a numerical method, such a reflected singularity is undesirable.

The reflected singularity from a discontinuous lower order term is weaker than the singularity of the  incident wave. For the  equation
$$
\partial_t U \ +\ A_1 \partial_1 U \ +\
\sigma \,
 {\bf 1}_{x_1>0}
 \, CU\ =\ 0,
$$
if $C$ is diagonal in a basis diagonalising $A_1$, the reflections are avoided. The ease of eliminating reflections for this problem with $d=1$ is deceptive. No such simple remedy exists in dimensions $d>1$. For symmetric hyperbolic systems $A_1=A_1^*$, it is wise to choose $C=C^*\ge 0$ so that the absorption term is  dissipative in the $L^2(\RR^d)$ norm.

Consider next the wave equation with friction $\partial_{tt}v-\partial_{11}v+ 2 \sigma\, {\bf 1}_{x_1>0} \,\partial_t v=0$ written in characteristic coordinates $(U_1,U_2)=(\partial_tv-\partial_1v,\partial_tv+\partial_1v)$ with absorption $B=\sig C$:
$$
\partial_t U \ +\
\begin{pmatrix}
1&0\\
0&-1
\end{pmatrix}
\partial_1U
\ +\
\sigma \,
 {\bf 1}_{x_1>0}
 \, C\,U
\ =\
0\,,
\qquad
C\ =\
\begin{pmatrix}
1&1\\1&1
\end{pmatrix}
.$$
The absorption matrix $C$ is symmetric and nonnegative but does not commute with $A_1$. It produces unacceptably strong reflections. The absorption from  Israeli and Orszag \cite{Israeli:1981:ARB}, $\partial_{tt}v -\partial_{11}v + \sigma (\partial_tv+\partial_1v)=0,$ absorbs only rightward waves and corresponds to
$$
C=
\begin{pmatrix}
1&0\\0&0
\end{pmatrix}
\ =\
\pi_+(A_1),
$$
introducing the notation  $\pi_+(A_1)$ for the spectral projector on the eigenspace corresponding to strictly positive eigenvalues of $A_1$. The general nonnegative symmetric choice commuting with $A_1$ is a  positive multiple of
\begin{equation}
\label{eq:smartB}
C\ =\  \pi_+(A_1) + \nu\, \pi_-(A_1),
\qquad
\nu \ge 0.
\end{equation}
We call these {\it smart layers}. They dissipate the $L^2$ norm. As observed by Israeli and Orszag, the numerical performance of the smart layers is not as good as one would hope. One reduces reflections by choosing $\sigma(x)\ge 0$ vanishing to order $k\ge 0$ at the origin.  That reduces the rate of absorption and thereby increases the width of the layer required. The leading reflection by such  smart layers of incoming  wave packets of amplitude $O(1)$ and wavelength $\eps$ is $O(\eps^{k+1})$.  The leading  reflection   is linear in $\sigma$. In section \ref{sec:hml}, we introduce the method of Harmoniously Matched Layers which remove the leading order reflections (at all angles of incidence) by an extrapolation.\\

\noindent
{\bf Open problem.}   {\sl Repeated extrapolation further reduces the order of reflection. It is easy to program and it is  possible that an optimization could pay  dividends.}\\

Elaborate absorbing layer strategies, like B\'erenger's PML
introduce operators related to but often more complicated than the original operator $L$.  The operators in the absorbing layer and in the domain of interest may not be the same.    For the case of a layer in $\{x_1>0\}$, absorbing layer algorithms solve a transmission problem for an unknown $(V,W)$ where $V$ is a $\CC^N$ valued function on $x_1<  0$ and $W$ is a function on $x_1>  0$. The equations in  $x_1>0$ are chosen to be absorbing and the transmission problem weakly reflective. The ingenious innovation of B\'erenger was to realize that the operator $R$ in the layer can differ  substantially from  $L$.  He increased the number of unknown functions in the layer. So $W$ is $\CC^M$ valued with  $M>N$.

The pair $(V,W)$ is determined by a well posed transmission problem,
 \begin{equation}
 \label{eq:LR}
 LV\ =\ F
\quad
{\rm on}
\quad
\RR^{1+d}_-
\ :=\ \{(t,x)\,:\,x_1<0\}
\qquad
R\, W \ =\ 0
\quad
{\rm on}
\quad
\RR^{1+d}_+\,,
\end{equation}
with the homogeneous transmission condition
\begin{equation}
\label{eq:caN}
(V,W)\ \in \
\caN
\quad
{\rm on}
\quad
\{x_1= 0\}.
\end{equation}
Here $\caN\subset\CC^N\times \CC^M$ is a linear subspace.
\hskip -.15cm
\footnote{Transmission conditions which involve derivatives
can also be treated.  The algorithms of B\'erenger and  our
HML do not require that generality.}
The choice of the hyperbolic operator $R$ and transmission condition $\caN$ is made with three goals,
\begin{itemize}
\item[$\bullet$]
 The transmission problem is well posed,
and  not hard to approximate numerically.
\item[$\bullet$]
Waves from the left are  at most weakly
reflected at $x_1=0$.
\item[$\bullet$]
Waves moving rightward decay rapidly in $x_1>0$
so that the layer can be chosen thin.
\end{itemize}

The criterion for perfection that we adopt is that
of Appelo, Hagstr\"om and Kreiss \cite{Appelo:2006:PML}
. In the case of one absorption, it is formulated as follows.

\begin{definition}
\label{def:perfection}
A well posed  transmission problem is {\bf perfectly matched}
when for all $F$ supported in $x_1<0\,,\, t\ge 0$,
the solution supported in $t\ge 0$ satisfies $V=U\big|_{x_1< 0}$.
\end{definition}

\noindent
We prove
in \S \ref{subsec:perfectionberenger}
 that B\'erenger's method with
one discontinuous absorption $\sigma_1$  is perfect in this sense.

\vskip.2cm

In practice one does not absorb in only one
direction and the computational
domain is rectangular.
We give in \S \ref{subsec:changevariable} a definition with absorptions
in more than one direction and a proof of perfection.

\vskip.2cm
The strategy of B\'erenger
 is quite ingenious.  For an artificial boundary in two dimensions at
 $\{x_1=r\}$ and  domain of interest $\{x_1<r\}$ it consists of two steps.  The first is a doubling
 of the system and the second is insertion of an absorption
 term in $\{x_1>r\}$.  The doubled system involves the
 unknown $\widetilde U:=(U^1,U^2)\in
 \CC^N\times\CC^N$.  When $F=0$, the doubled equation without dissipation is
 \[
 \partial_t U^ j \ +\
 A_j\partial_j (
  U^1
 +
   U^2
    )
   \ = \ 0\,,
   \qquad
   j=1,2\,.
\]
 The system with  damping in $x_1$
 changes the $j=1$ equation to
\[
\partial_t  U^1 \ +\
A_1\partial_1( U^1 + U^2 ) \ +\ \sigma(x_1)\, U^1
\ = \ 0\,, \quad {\rm supp}\, \sigma \ \subset \{x_1\ge r\}\,.
\]
 Then $U:=\sum_j  U^j$ satisfies $L(\partial)U=0$
 in $x_1 < r$.
 In practice it is the restriction of $U$ to $x_1< r$ that is of
 interest.
 There are three distinct ways to view this.  One can think of
 the unknowns as $U$ defined in $x_1< r$ and $\widetilde U$
 in $x_1\ge r$ with the transmission condition that
 $A_1 U=A_1(  U^1 +  U^2)$ on $x_1=r$.
 One is given initial values of $U$ and takes initial values
 of $\widetilde U$ vanishing.  This is the most natural choice
 and the one presented by B\'erenger.

From the computational point of view it is simpler to have the same unknowns throughout.  The simplification is greater when one passes from the half space case to a computational domain equal to a rectangular domain in $\RR^d$. One introduces $\widetilde U$ everywhere with transmission condition $[A_1(  U^1+  U^2)]=0$ where $[ * ]$ denotes the jump at $x_1=r$.  The transmission condition is then equivalent to the validity of the differential equation satisfied by $\widetilde U$ in all of $\RR^d$. When one uses $\widetilde U$ everywhere, the initial values of $\widetilde U$ are taken equal to zero outside the computational domain.  The initial values are constrained to satisfy $U=\sum_j  U^j$ within the computational domain.   The choice is otherwise arbitrary. For the case of the doubling above the choice $U^j(0,x) = U(0,x)/2$ for $j=1,2$ is common.

  If the domain of interest is $|x_1|\le r$ one would
 choose $\sigma>0$ on $|x_1|>r$ and vanishing for
 $|x_1|<r$.   The transmission condition is $[A_1( U^1
 +  U^2)]=0$ with the jump   at $x_1=r$ and
 also at $x_1=-r$.

In a rectangular geometry in $\RR^d$ introduce $\tiu:=(U^1,\dotsc,U^d)$, where $U^l \in \CC^N$ for $1 \le l \le d$. Then $\widetilde U$ with values in $\CC^{Nd}$ is required to satisfy (in the case $F=0$),
\begin{equation}\label{PMLgen}
  (\til(\partial _t,\partial_x)\,\tiu)_l
  \ :=\
  \partial _tU^l
  \ +\
  \ds A_l \partial _{l}(\sum_{j=1}^d U^j)
  \ +\ \sigma_l(x_l) U^l
  \ =\
  0, \qquad  1\leq l \leq d  .
\end{equation}
Each absorption coefficient $\sigma_l(x_l)\ge 0$ depends on only one variable.  It is strictly  positive between the inside rectangle and a larger outside rectangle.  In the layer between the rectangles the solution is expected to decay. If $\tiu$ solves \eqref{PMLgen}, then $U= \sum_{j=1}^d U^j$ solves \eqref{hyp} on the set $\{x\,:\, \sigma_l(x_l) =0\ \mbox{ for }\ 1 \le l \le d\}$ including the inner rectangle. In the case considered by B\'erenger the $\sigma$ were discontinuous and the equations \eqref{PMLgen} are equivalent to transmission problems where on the discontinuity surface of $\sigma_j$ one imposes the transmission condition of continuity of $A_j \sum_\ell U^\ell$.

Our first  technique is the energy method. In \S \ref{subsec:elliptic} we show that if $(\xi_1\,,\, \dots\,,\, \xi_d)=0$ does not meet the characteristic variety of $L$ then the B\'erenger method is well posed {\it without loss of derivatives}.  This applies in particular to linearized elasticity and suggests that in some ways the B\'erenger method is better adapted to that situation than the Maxwell equations for which it was intended. In \S \ref{subsec:maxwell2DTE} we give a nontrivial extension of the method of M\'etral and Vacus to show that B\'erenger's method for the Maxwell equations in dimension $d=2$ (resp. $d=3$) is well posed provided that $\sigma_j(x_j)\in W^{1,\infty}(\RR_{x_j})$ (resp.  $\sigma_j(x_j)
\in W^{2,\infty}(\RR_{x_j})$).  The method introduces a norm that is the sum of $L^2(\RR^d_{t,x})$ norms of suitable differential operators $P_\alpha(D)$ applied to $U$. It has  the property that the norm at time $t_1$ is estimated in terms of the norm at time $t_2$. If one introduces the vector of unknowns $P_\alpha(D) U$ this  shows that the B\'erenger problem becomes strongly well posed without loss of derivatives.   Such transformations are typical of weakly well posed problems.
(see the Dominics' proof
of Theorem 1.1 in \S IV.1
of
\cite{Taylor:1981:PO}).

When such an estimate
is known, we  prove sharp finite speed in \S \ref{sec:finitespeed}  and perfection in \S \ref{subsec:changevariable} and \S\ref{subsub:PMRB}, the latter
concerned with
several variants of the Berenger strategy.
The perfection proof passes by  a study of the Laplace transform
on $\{ {\rm Im}\,\tau =0\}$.
The transformed problem is conjugated to
the problem without absorption by a $\tau$ dependent change
of independent variable $x$, an idea inspired by
\cite{Diaz:2006:TDA}.
\\

Our second method is the Fourier-Laplace method.
B\'erenger introduced his PML for Maxwell's equations with $\sigma$ piecewise constant.  Using a computation which resembles plane wave analysis of reflections for problems without lower order terms, B\'erenger argued that the layers were perfectly matched for all wave numbers and all angles of incidence. Using variants of the  same approach  other closely related PML were constructed afterward. Performance is observed to be enhanced using $\sigma$ which are not discontinuous. Twice differentiable cubic functions are the most common. The B\'erenger method is a very good method for Maxwell's equations.
The Fourier-Laplace method gives a framework for understanding
the computations of B\'erenger.  In addition, it is the only
method we know for proving well posedness of B\'erenger's
PML with discontinuous $\sigma$
for Maxwell's equations.

Plane wave  analysis is  sufficient to study reflection and transmission for linear constant coefficient operators without lower order terms.
Problems with lower order terms require other tools as it is no longer true that the plane waves generate all solutions. The first level of generalization is to use the Fourier-Laplace transform for problems where an absorbing layer occupies $x_1\ge 0$ and both $L$ and $R$ have constant coefficients.
Hersh \cite{Hersh:1963:MPS} found necessary and sufficient conditions for  (weak) well posedness of transmission problems.  We recall those ideas in \S \ref{sec:Hershcond}
including the modifications needed for characteristic interfaces,
 and verify in \S\ref{subsec:wellposed} that  the condition is satisfied for the
 B\'erenger splitting of general systems with one discontinuous absorption coefficient. To  our knowledge this is the first proof that the B\'erenger split transmission problem with discontinuous $\sigma(x_1)$ is well posed.

We give necessary and sufficient conditions for perfection at a planar
interface.  In \S \ref{subsec:perfectionberenger}
 we verity that the condition is satisfied for the B\'erenger
splitting.
%Our current favorite is that
%of Ziolkowski, Abrarbanel and Gottlieb which we call the
%Z-A-G method.

In \S \ref{sec:continuation} we prove  that
in the case of Maxwell's equations (and not in general)
 the perfection criterion follows by analytic continuation
 from the plane wave identities established by B\'erenger.

In \S \ref{subsec:FLAV} we  prove using the Fourier-Laplace method that B\'erenger's method with  one coefficient $\sigma_1(x_1)\in {\rm Lip}(\RR_{x_1})$ is well posed and perfectly matched. In our use of the Fourier-Laplace method, including this one, a central role is played by the Seidenberg-Tarski Theorem estimating the  asymptotic behavior of functions defined by real polynomial equations and inequalities. The Fourier-Laplace method is limited to coefficients that depend only on $x_1$.\\

Our third method of analysis is to study the behavior of short wavelength  asymptotic solutions. For such solutions we examine in \S\ref{sec:amplification} the decay in the absorbing layers, and reflections at discontinuities of $\sigma_j(x_j)$ or its derivatives when smoother transitions are used.
For problems other than Maxwell, Hu \cite{HU:1996:ABC} and B\'ecache, Fauqueux , and Joly \cite{Becache:2003:SPM} have already shown that the supposedly absorbing layers may in fact lead to growth. The study of short wavelength solutions in the layer yields precise and clear criteria, also valid for variable coefficients, explaining the phenomenon.

The analysis of the reflection of short wavelength wave packets at the interface with the layer also  leads us to propose in \S \ref{sec:hml}, a new absorbing layer strategy which we call Harmoniously Matched Layers. The method starts with a  smart layer for a symmetric hyperbolic system.  Then for wavelength $\eps$ asymptotic solutions of amplitude $O(1)$ and discontinuous $\sigma$, the leading order reflected wave at nonnormal incidence typically has amplitude proportional to $\sigma \eps$.  Therefore an extrapolation using computations with two values of $\sigma$ eliminates the reflections proportional to $\sigma$. This yields a method with leading order reflection $O(\eps^2)$ at all angles of incidence. The resulting method inherits the simple $L^2$ estimates of the symmetric systems. More generally if the first discontinuous derivative of the absorption coefficient is the $k^{\rm th}$ then the reflection is $O([D^k\sigma]\eps^{k+1})$ and the same extrapolation removes the leading order reflection.
  In \S  \ref{subsec:numerics}  we investigate several implementations of this idea and show that the method with cubic $\sigma$ is  competitive with that of B\'erenger with the same $\sigma$. On short wavelengths or random data it performs better than the B\'erenger method. On long wavelengths B\'erenger performs
better.\\

Though we provide satisfactory answers to a wide range of questions about absorbing layers, there is a notable gap.\\

\noindent
{\bf Open problem.} {\sl
For  the original strategy of B\'erenger for Maxwell's equations
with discontinuous absorptions in more than one
direction we do not know if the resulting problem is
well posed.
}{\bf Discussion.}
{\bf 1.}  In \S \ref{subsub:PMRB} we prove
 well posedness and perfection for a closely related method.
{\bf 2.}~In practice discontinuous
$\sigma$ have been abandoned, but it is striking that this
problem remains open.
{\bf 3.} Once well posedness is proved, perfection follows by the proof in \S \ref{subsec:changevariable}.

%%%%%%%%%%%%%%%%%%%%%%%%%%%
%  New section
%%%%%%%%%%%%%%%%%%%%%%%%%%
\section{Well posed first order Cauchy problems }\label{sec:preliminaires}

\subsection{Basic definitions}\label{subsec:cauchy}
Consider a  first order system of partial differential equations for $\CC^N$ valued functions on $\RR^{1+d}$,
\begin{equation} \label{gen}
\cl(x,\partial_t,\partial_x)\, U
\ :=\
\partial _tU
\ +\
\sum_{l=1}^d \,{\cal A}_l \,\partial _{l}U
\ +\
 {\cal B}(x) \,U
 \ =\
 0.
\end{equation}
The principal part of $\cl$, denoted $\cl_1$,
\[
\cl_1(\partial_t,\partial_x)\ :=\ \partial _t \ +\
\sum_{l=1}^d {\cal A}_l\, \partial _{l}\,,
\]
has constant matrix coefficients   $ {\cal A}_l$.
In the B\'erenger strategy, the operators $\widetilde L$
are the centerpieces and they differ from $L$.   It is for
this reason that we introduce $\mathcal L$ that can be
 $L$
or $\widetilde L$.

\begin{definition}
The {\bf characteristic variety}
 $ {\rm Char}(\cl)\subset \CC^{1+d}\setminus \{0\}$ of $\cl$ is the
 set of $(\tau,\xi)$ such that $\det \cl_1(\tau,\xi)=0$.
\end{definition}

\begin{definition}
The {\bf smooth variety hypothesis} is satisfied at $(\utau,\uxi)\in {\rm Char}(\cl)$ if there is a conic neighborhood $\Omega$ of $\uxi\in \RR^{d}\setminus \{0\}$ and a $C^\infty$ function $\xi \mapsto\tau(\xi)$ on $\Omega$ so that on a neighbourhood of $(\utau,\uxi)$, the characteristic variety has equation $\tau=\tau(\xi)$.  At such a point the associated {\bf group velocity} is defined to be $\bfv:=-\nabla_\xi\tau(\xi)$.
\end{definition}

\begin{example}
This hypothesis holds if an only if for $\xi$ near $\uxi$ the spectrum of $\cl(0,\xi)$ near $-\utau$ consists of a single point with multiplicity independent of $\xi$.
For the polynomial $(\tau+\xi_1)(\tau^2-|\xi|^2)$ with $d>1$ the hypothesis
fails at and only at $\tau+\xi_1=0$ where two  sheets of the variety are
tangent.  Replacing the first factor by $\tau+c\,\xi_1$ with $ c>1$
the hypothesis fails where the two sheets cross transversally.
For $0\le c<1$ the hypothesis holds everywhere.
\end{example}

The Cauchy problem for $\cal L$ is to find a solution $U$ defined on $[0,\infty[\times \RR^d $ satisfying  \eqref{gen} with prescribed initial data $U(0,\cdot)$.
\begin{definition} \label{faible}
The Cauchy problem for $\cl$ is {\bf weakly well posed} if there exist $q>0$, $K>0$, and $\alpha \in \RR$ so that for any initial values in $H^{q}(\RR^d)$,
there is  a unique solution $U \in
  \mathcal{C}^0([0,+\infty[\,;\,L^{2}(\RR^d))$ with
\begin{equation}\label{stab}
  \forall t \geq 0,
  \qquad
  \| U(t,\cdot)\|_{L^2(\RR^d)}\leq K e^{\alpha t}  \|U(0,\cdot)\|_{H^{q}(\RR^d)}.
\end{equation}
When the conclusion holds with $q=0$, the Cauchy problem is called {\bf strongly well posed}.
\end{definition}

\begin{theorem} \label{fort}

\begin{enumerate}[label={\it(\roman{*})}, ref={\it(\roman{*})},leftmargin=0.7 cm]
\item\label{1}The Cauchy problem  for $\cl_1$ is weakly well posed if and only if for each $\xi \in \RR^d$, the eigenvalues of  $\cl_1(0,\xi)$ are real.
        \item\label{2}The Cauchy problem for $\cl_1$ is strongly well posed if and only if for each $\xi \in \RR^d$, the eigenvalues of  $\cl_1(0,\xi)$ are real and $\cl_1(0,\xi)$ is uniformly diagonalisable,
        there is an invertible $S(\xi)$ satisfying,
            $$
            S(\xi)^{-1}\cl_1(0,\xi)\,S(\xi) = {\rm diagonal},
            \qquad
             S\,,\,S^{-1} \in L^\infty(\RR^d_\xi )\,.
             $$
        \item\label{3}If ${\cal B}$ has constant coefficients, then the Cauchy problem for $\cl$ is weakly well posed if and only if there existes $M \ge 0$ such that for any $\xi \in \RR^d$,
        $\det \cl(\tau,\xi)=0 \implies |\Im \tau| \le M$.
\end{enumerate}
\end{theorem}

\begin{remark}
 \label{rq:semisimple}
\newline
{\bf 1.}
The algebraic conditions in $\ref{1}$ and $\ref{3}$
express \textit{weak hyperbolicity}, in the sense of
  G{\aa}rding. The necessity of uniform diagonalisability
in  $\ref{2}$ expressing \textit{strong hyperbolicity }is  due to Kreiss \cite{Chazarain:1981:ITE}, \cite{Kreiss:1989:IBV}.

\noindent
{\bf 2}.
An application of Gr\"{o}nwall's inequality shows that
if $\cl_1$ satisfies the condition of Theorem
\ref{fort}, $(ii)$, then for all $B(x)\in L^\infty(\RR^d\,;\,{\rm Hom}(\CC^N))$,
the Cauchy problem for $\cl_1 + B$ is
strongly well posed.

 \noindent
{\bf 3.}  By property $\ref{2}$, if $\cl$ is strongly hyperbolic, then every eigenvalue $-\tau$ of $\cl_1(0,\xi)$ is semi-simple. Equivalently, for any $(\tau,\xi) \in {\rm Char}(\cl)$ the eigenvalue $0$ of $\cl_1(\tau,\xi)$ is semi-simple, \textit{i.e.} its geometric multiplicity is equal to its algebraic multiplicity. It is equivalent to saying that $\ker \cl_1(\tau,\xi)= \ker(\cl_1(\tau,\xi))^2$, or that $\CC^N=\ker \cl_1(\tau,\xi)\oplus \im \cl_1(\tau,\xi)$.
\end{remark}
\subsection{Characteristic variety and projectors for B\'erenger's
$\widetilde L$}\label{subsec:charac}

To study the Cauchy problem
 for B\'erenger's
split operators $\widetilde L$ one starts with a study
of the characteristic variety.
The  coefficients of  B\'erenger's operator $\widetilde L$
are the $dN\times dN$ matrices,
\begin{equation}
\label{eq:Ltildecoeff}
\tia_l := \begin{pmatrix}
  0 & \dots & \dots &\dots & 0 \\
  \vdots  & & & & \vdots \\
  A_l & \dots &\dots & \dots & A_l \\
  \vdots & & & & \vdots \\
  0& \dots & \dots & \dots & 0
\end{pmatrix} ,\qquad
B(x):=\begin{pmatrix}
  \sigma_1(x_1) I_N&\dots & 0 \\
  \vdots & \ddots & \vdots \\
  0 & \dots &\sigma_d(x_d) I_N
  \end{pmatrix}.
\end{equation}
The principal symbol of $\til$ is
\[
\til_1(\tau,\xi)=
    \begin{pmatrix}
\xi_1 A_1 +\tau I_N& \xi_1 A_1 &  \dots & \xi_1 A_1 \\
\xi_2 A_2 & \xi_2 A_2 +\tau I_N&  \dots & \xi_2 A_2 \\
\vdots & \vdots & \ddots &\vdots \\
\xi_d A_d & \xi_d A_d &  \dots & \xi_d A_d +\tau I_N
    \end{pmatrix}.
\]
\begin{theorem}
\label{prop:charvar}
\begin{enumerate}
[label={\it(\roman{*})}, ref={\it(\roman{*})},leftmargin=0.7cm]
\item  The characteristic polynomial of $\widetilde L$ is
\begin{equation}
\label{eq:det1}
\det \til_1(\tau,\xi)
\ =\ \tau ^{N(d-1)} \,
\det L(\tau,\xi)\,.
\end{equation}
The
polynomial associated to the full symbol including the absorption
is
\begin{equation}
\label{eq:det}
\det \widetilde L(\tau,\xi)
\ =\
\ds\det L\Big(
\prod_{j=1}^d (\tau +\sigma_j)
\,,\,
\xi_1
\prod_{j\ne 1}(\tau +\sigma_j)
\,,\,
\xi_2
\prod_{j\ne 2}(\tau +\sigma_j)
\,,\,
\dots
\,,\,
\xi_d
\prod_{j\ne d}(\tau +\sigma_j)
\Big).
\end{equation}
\end{enumerate}
\noindent
If $\txi\in {\rm Char}\,L$ with $\tau\ne 0$, the following properties hold.
\begin{enumerate}
[label={\it(\roman{*})}, ref={\it(\roman{*})},leftmargin=0.7cm,resume]
\item
The mapping
\[
{\cal S}:\ \widetilde{\Phi}= (\Phi_1,\cdots,\Phi_d)
 \quad \mapsto\quad
 -\,\sum_{j=1}^d\, \Phi_j
 \]
is a linear bijection from $\tiklu$ onto $\klu$. Its inverse is given by
\[
\Phi
\quad \mapsto\quad
 \Big(\frac{\xi_1}{\tau} A_1\Phi\,, \,\dots\,,\frac{\xi_d}{\tau} A_d\Phi\Big).
\]
\item
The kernel of the adjoint $\til_1(\tau,\xi)^*$
     is equal to the set of vectors $\widetilde{\Phi}=(\Phi,\cdots,\Phi)$ such that $\Phi\in \ker L_1(\tau,\xi)^*$. The range of $\til_1(\tau,\xi)$ is equal to
    the set of vectors $\widetilde{\Psi}=\bigl(\Psi_1,\cdots,\Psi_d\bigr)$ such that
    $(\sum_{j=1}^d\Psi_j,\Phi)=0$ for all $\Phi \in \ker\,L_1(\tau,\xi)^*$.
\item
If moreover the eigenvalue $0$ of $L_1(\tau,\xi)$ with $\tau\ne 0$ is semi-simple, the eigenvalue $0$ of $\widetilde L_1(\tau,\xi)$ is semi-simple.  Equivalently,
    $$
    \ker\til_1(\tau,\xi) \ \oplus\
    {\rm Range} \,\til_1(\tau,\xi) \ =\ \CC^{dN}\,.
    $$

\end{enumerate}
\end{theorem}
\begin{proof}
$(i)$  Adding the sum of the other rows   to the first row in the determinant of $\til_1(\tau,\xi)$ yields,
    \[
        \det\til_1(\tau,\xi) =
        \begin{vmatrix}
            L(\tau,\xi) & \dots & \dots & L(\tau,\xi)  \\
            \xi_2 A_2 &  \xi_2 A_2+\tau I_N & \dots &\xi_2 A_2 \\
            \vdots & & \ddots & \vdots \\
            \xi_d A_d & \dots &\dots &  \xi_d A_d+\tau I_N
        \end{vmatrix}
    \]
   Substracting the first column from the others yields,
    \[
        \det\til_1(\tau,\xi) \ =\
       \begin{vmatrix}
            L(\tau,\xi) & 0 & \dots & 0 \\
            \xi_2 A_2 &  \tau\,I_N & \dots &0 \\
            \vdots & & \ddots & \vdots \\
            \xi_d A_d & 0 &\dots & \tau I_N
        \end{vmatrix}
        \]
 The first result follows.
For the second write,
    \[
        \det\til(\tau,\xi) \ =\
        \begin{vmatrix}
            \xi_1 A_1+(\tau+\sigma_1)I_N  & \xi_1 A_1 & \dots & \xi_1 A_1 \\
            \xi_2 A_2 &  \xi_2 A_2+(\tau+\sigma_2) I_N & \dots &\xi_2 A_2 \\
            \vdots & & \ddots & \vdots \\
            \xi_d A_d & \dots &\dots &  \xi_d A_d+(\tau+\sigma_d) I_N
        \end{vmatrix}
    \]
    For each $i$ divide the $i^{\rm th}$ row by $\tau+\sigma_i$ to find,
   \[
        \det\til(\tau,\xi)
        \ =\
         \ds\prod_{j=1}^d \ (\tau +\sigma_j)^N\,
       \det\til_1(1,\frac{\xi_1}{\tau+\xi_1},\cdots,\frac{\xi_d}{\tau+\xi_d})\,.
    \]
    By formula \eqref{eq:det1} this implies
   \[
        \det\til(\tau,\xi) \ =\  \ds\prod_{j=1}^d (\tau +\sigma_j)^N
       \det L(1,\frac{\xi_1}{\tau+\xi_1},\cdots,\frac{\xi_d}{\tau+\xi_d}),
    \]
which is equivalent to \eqref{eq:det}.\\

\noindent $(ii)$ Suppose that
  $0\ne {\widetilde \Phi}=( \Phi_1,\cdots, \Phi_d)\in \ker\,\til_1(\tau,\xi)$.
    Then, for any $l$,
        \begin{equation}
        \label{eq:phi_l}
        \tau \Phi_l
        \ +\ \xi_lA_l\sum_{j=1}^{d} \Phi_j
        \ =\
        0,
        \end{equation}
    Add to find
        $$
        L_1(\tau,\xi)
        \,   \sum_{j=1}^{d}\Phi_j\ =\ 0.
        $$
        Therefore the map $\widetilde \Phi\mapsto -\sum_j \Phi_j$
        maps $\ker\,\widetilde L_1\txi$ to $\ker\,L_1\txi$.

    If $\sum_{j=1}^{d}\Phi_j=0$, equation \eqref{eq:phi_l} implies
    that all the $\Phi_j$ vanish since $\tau \ne 0$. Therefore the mapping
    is injective.

Let $\Phi \in \klo$. Define
    \begin{equation}
    \label{eq:inverse}
        \Phi_j\ =\ \frac{\xi_j}{\tau}A_j \Phi\,.
    \end{equation}
    This defines an element $\widetilde \Phi=(\Phi_1,\cdots,\Phi_d)$ in $\tiklo$ with   ${\cal S}\widetilde \Phi= \Phi$, so the mapping is surjective with inverse given by \eqref{eq:inverse}.\\

\noindent$(iii)$        Since
        $
            \til_1(\tau,\xi)^*\widetilde{\Phi}=
            (L_1(\tau,\xi)^*\Phi,\ldots ,L_1(\tau,\xi)^*\Phi)
        $
        it follows that the set of $\widetilde \Phi$ is included in the kernel.
        Since the matrices are square,   $\klu$ and $\klu^*$ have the same dimension. The set of $\widetilde \Phi$ has dimension equal to this common dimension which by $(ii)$ is equal to the dimension of $ \ker \, \widetilde L_1\txi $  proving that they exhaust the kernel.
    The last property follows directly from the fact that  $\im \tilutx$ is the orthogonal of $\klu^*$.\\

\noindent $(iv)$ It suffices to show that the intersection of these spaces consists of the zero vector.  Equivalently, it suffices to show that there is no  $\Phi\ne 0$ in   $\klu$ such that
 \[\forall \Psi \in \klu^*, \qquad \ds (\sum_{j=1}^d \frac{\xi_j}{\tau}A_j\Phi, \Psi)=0\,.
  \]
The quantity above is equal to $-(\Phi,\Psi)$, and $\Phi$ would belong to $ (\klu^*)^\bot= \im \lutx$. Since $\tau\ne 0$ and $\klu\cap \im\lutx={0}$, this would imply that $\Phi = 0$, leading to a contradiction.\qquad
 \end{proof}

Denote by  $\Pi_{L}(\tau,\xi)$ (resp.  $\Pi_{\til}(\tau,\xi)$)
the spectral projector onto the kernel of  $L_1(\tau,\xi)$ (resp.  $\til_1(\tau,\xi)$) along its range.  For $L$ it is given by
$$
\Pi_{L}(\tau ,\xi)=\frac{1}{2\pi i}
\oint_{|z|=\rho}
\big(z\,I-L_1(\tau ,\xi)\big)^{-1}\
dz
$$
with $\rho$ small. Like the characteristic variety, $\Pi_{L}$ depends
only on the principal symbol $L_1$.  It is characterized by,
\begin{equation}
\label{eq:proj}
\Pi_{L}^2=\Pi_{L},\qquad
\Pi_{L}\,L_1(\tau,\xi)=0,\qquad
L_1(\tau,\xi)\,\Pi_{L}=0,\qquad
\rk \Pi_{L}= {\rm dim}\,\ker L_1(\tau,\xi),
\end{equation}
where the $\tau,\xi$ dependence of $\Pi_{L}$ is suppressed for ease of reading.
The first three conditions assert that $\Pi_{L}(\tau,\xi)$ is a projector
annihilating  ${\rm Range}\, L_1(\tau,\xi)$ and projecting onto a subspace
of ${\rm Ker}\,L_1(\tau,\xi)$.  That it maps onto the kernel is implied by
the last equality.

\begin{prop} \label{prop:projection}
The matrix
 $\Pi_{\til}(\tau,\xi)$
is given by
$$
\Pi_{\til}(\tau,\xi)
\ =\
-\begin{pmatrix}
  \ds\frac{\xi_1A_1}{\tau}\,\Pi_{L}(\tau,\xi) & \dots &
       \ds\frac{\xi_1A_1}{\tau }\,\Pi_{L}(\tau,\xi)\\
  \vdots & & \vdots \\
  \ds\frac{\xi_dA_d }{\tau }\,\Pi_{L}(\tau,\xi) & \dots &
       \ds\frac{\xi_dA_d }{\tau }\,\Pi_{L}(\tau,\xi)
                       \end{pmatrix} \,.
   $$
\end{prop}

\noindent
\begin{proof}
Call  the matrix on the right $M(\tau,\xi)$.
The properties of the projectors associated to $L$
yield formulas for the  $(i,j)$ block of the products
 $$
 (M(\tau,\xi)\til_1(\tau,\xi))_{i,j}
 \ =\
  -\frac{\xi_iA_i}{\tau}\,\Pi_{L} L_1
  \ =\
  0,
 \qquad
 (\til_1(\tau,\xi)M(\tau,\xi))_{i,j}= -\frac{\xi_iA_i}{\tau}\,L_1\Pi_{L} L_1 =0,
\quad
 {\rm and,}
$$
$$
(M(\tau,\xi)M(\tau,\xi))_{i,j}\ =\
\frac{\xi_iA_i}{\tau^2}\Pi_{L}(L_1-\tau I)\Pi_{L}
\  =\
-\frac{\xi_iA_i}{\tau}\Pi_{L}^2
  \ =\
  (M(\tau,\xi))_{i,j}\,.
$$
This proves the first three equalities of \eqref{eq:proj}.
Since $M$ projects onto a subspace of $\ker \widetilde L_1$,
$\rk M \le {\rm dim}\,{\rm Ker}\, \widetilde L_1$.
Apply  $M$ to a vector $(\Psi, 0,\dots,0)$ and compare
with
part $(ii)$  of
Theorem \ref{prop:charvar} to see that the range of $M$
contains $\ker \widetilde L_1(\tau,\xi)$ so
$\rk M \ge {\rm dim}\,\ker {\widetilde L_1}$.
This proves the last equality of   \eqref{eq:proj}.
\end{proof}

\begin{remark}\newline
{\bf 1.}
The characteristic varieties of
$L$ and $\til$ are identical
in $\tau\ne 0$.

\noindent{\bf 2.}
In particular,
 the smooth variety hypothesis is satisfied at $(\tau,\xi)$
 with  $\tau\ne 0$ for one if and only if it holds for both,
 and the varieties have the same equations and the same
 group velocities.

\noindent{\bf 3.}  When the smooth variety
hypothesis is satisfied, the spectral projection
$
\Pi_{\til}(\tau(\uxi),\uxi)$
is analytic in $\xi$, hence of constant rank.  It follows that
$0$ is a semi-simple eigenvalue of $\til(\tau(\xi), \xi)$
on a conic neighborhood of $\uxi$.
\end{remark}

If the eigenvalue $0$ of $L_1(\tau,\xi)$ is semi-simple,   the kernel and the range of $ L_1(\tau,\xi)$
are complementary subspaces as mentioned in Remark \ref{rq:semisimple}
{\bf 3.}, and the partial inverse
$Q_L(\tau,\xi)$ of $ L_1(\tau,\xi)$ is uniquely determined by
\begin{equation}\label{eq:inversedef}
Q_L(\tau,\xi)\,
\Pi_{L}(\tau,\xi)\ =\ 0,
\qquad
Q_L(\tau,\xi)\,
L_1(\tau,\xi)
\ =\
I-\Pi_{L}(\tau,\xi)\,.
\end{equation}
The partial inverse $Q_{\tilde L}(\tau,\xi)$ is defined
in the same way from ${\tilde L}_1(\tau,\xi)$.

\subsection{The Cauchy problem for B\'erenger's split operators}\label{subsec:cauchysplit}
Part $(i)$ of Theorem \ref{prop:charvar} proves the following.
\begin{corollary}\label{cor:weak}
If the Cauchy problem for $L$ is weakly well posed, then so is the Cauchy problem for the principal part  $\til_1$.
\end{corollary}
An important observation is that though the
Cauchy problem for $\widetilde L_1$ is at least
weakly well posed, the root $\tau=0$ is for
all $\xi$ a multiple root.  When there
are such multiple roots it is possible that
order zero perturbations of $\tilde L_1$
may  lead to ill posed Cauchy problems.
The next example shows that this
phenomenon  occurs for the
B\'erenger split operators with constant
absorption
$\sigma_j$. Theorem \ref{th:thetheorem}
 shows that
when $\tau=0$ is a root of constant multiplicity
of $\det L_1(\tau,\xi)=0$, the constant coefficient
B\'erenger operators
have well posed Cauchy problems.   Cases
where the problem are strongly well posed are
identified.  In the latter cases, arbitrary bounded
zero order perturbations do not destroy the strong
well posedness.

\begin{example}\newline
{\bf 1.}   For
$
L := \partial_t + \partial_1+\partial_2,
$
$
\det L(\tau,\xi) \ =\ \tau + \xi_1 +\xi_2$.
Therefore $\tau=0$ is a root if and only if $\xi_1 +\xi_2=0$. The doubled system with absorption $\sigma=1$ in $x_1$ is
$$
\widetilde L
\ :=\
\partial_t
+
\left(
\begin{matrix}
1 & 1\cr
0 & 0
\end{matrix}
\right)
\partial_1
\ +\
\left(
\begin{matrix}
0 & 0\cr
1 & 1
\end{matrix}
\right)
\partial_2
\ +\
\left(
\begin{matrix}
1 & 0\cr
0 & 0
\end{matrix}
\right)\,,
\qquad {\rm and},
$$
$$
\det \widetilde L(\tau,\xi,-\xi)
\ =\
\det
\left(
\begin{matrix}
\tau +\xi +1 & \xi\cr
-\xi&\tau - \xi
\end{matrix}
\right)
\ =\
(\tau+\xi + 1)(\tau-\xi) + \xi^2
\ =\
\tau^2 + \tau -\xi\,.
$$
The roots of $\det \widetilde L(\tau,\xi,-\xi)=0$ are $ \tau \ =\ (-1\pm \sqrt{1+4\xi})/2$.  Taking $\xi \to -\infty$ shows that the Cauchy problem for $\widetilde L$ is not weakly well posed by Theorem \ref{fort}, \ref{3}.

\noindent
{\bf 2.}  More generally if $\tau + \xi_1 +\xi_2$ is a factor of
$\det L_1(\tau,\xi)$ then for $\sigma\ne 0$ the
operator $\widetilde L$ is not even weakly hyperbolic.
In this case \eqref{eq:det} implies that
$\tau^2 + \tau -\xi$ is a factor of
$\det\widetilde L(\tau,\xi,-\xi)$.

\noindent
{\bf 3.} This generalizes to linear hyperbolic factors
 in arbitrary dimension.
\end{example}

A key  tool is  the following special case of Theorem A.2.5 in
\cite{Hormander:2005:II}.

\begin{theorem}[Seidenberg-Tarski Theorem]\label{thm:seidenber}
  If $Q(\rho,\zeta)$,
$R(\rho, \zeta)$, and $S(\rho, \zeta)$ are polynomials with
real coefficients in the $n+1$ real variables $(\rho, \zeta_1,\dots,\zeta_n)$
and the set
$$
M(\rho) \ :=\ \big\{\zeta\ :\ R(\rho,\zeta) =0, \ S(\rho,\zeta)\le 0\big\}
$$
is nonempty when $\rho$ is sufficiently large,   define
$$
\mu(\rho) \ :=\
\sup_{\zeta\in M(\rho)} \ Q(\rho,\zeta).
$$
Then either $\mu(\rho)=+\infty$ for $\rho$ large,
or
there  are $a\in \QQ$ and $A\ne 0$ so that
$$
\mu(\rho) \ =\
A\, \rho^a\big(1 + o(1)\big),
\qquad
\rho \ \to\ \infty\,.
$$
\end{theorem}

\begin{theorem}
\label{th:thetheorem}
Suppose that $\tau=0$ is
an isolated root of constant multiplicity $m$ of $\,\det L_1(\tau,\xi)=0$.
\begin{enumerate}
[label={\it(\roman{*})}, ref={\it(\roman{*})},leftmargin=0.7cm]
\item  If the Cauchy problem for $L_1$ is strongly well posed,
then for arbitrary constant absorptions $\sigma_j\in \CC$,
the Cauchy problem for $\til_1+B$ is
weakly well posed.
\item  If the Cauchy problem for $L_1$ is strongly well posed, and if there is a $\xi \ne 0$ such that $\ker L(0,\xi)\ne \underset{\xi_j\ne 0}{\cap} \ker A_j$, then $\til_1(0,\xi)$ is not diagonalizable. Therefore the Cauchy problem for $\til$ is not strongly well posed.
\item   If the Cauchy problem for $L$ is
    strongly well posed and for all  $\xi$, $\klo = \underset{\xi_j\ne 0}{\cap} \ker A_j$, then the Cauchy problem for $\widetilde L$ is strongly well posed.
    This condition holds if $L_1(0,\partial_x)$ is elliptic, that is
    $\det L_1(0,\xi)\ne 0$ for all real $\xi$.
\end{enumerate}
\end{theorem}

\begin{remark}\newline
{\bf 1.}   Part $( i)$  is a generalisation of results in \cite{HU:2001:SPM} and Theorem 1 in \cite{Becache:2003:SPM}.  In the latter paper, B\'ecache {\it et al.}   treated the case $N=2$ assuming that the nonzero eigenvalues of $L_1(0,\xi)$ are of multiplicity one.  They conjectured that the result was true more generally. Like them  we treat the roots near zero differently from those that are  far from zero.  The treatment of each of these cases is different from theirs. The tricky part is the roots near zero.  We replace their use of Puiseux series by the related Seidenberg-Tarski Theorem \ref{thm:seidenber}.

\noindent
{\bf 2.}   Arbarbanel and Gottlieb \cite{Abarbanel:1997:AMA} proved $( ii)$   in the special case of Maxwell's equations.  The general argument below is simpler and yields a necessary and sufficient condition for loss of derivatives when the eigenvalue 0 of $L(0,\xi)$ is of constant multiplicity.

\noindent
{\bf 3.}   Part $( iii)$  is new, extending a result in
the thesis of S. Petit-Bergez \cite{Petit:2006:PFB}.
\end{remark}

\begin{proof}
\newline\noindent$( i)$
For $\xi\in \RR^d\setminus 0$, define for $\rho\in \RR_+$,
\[
E(\rho)
\ :=\
\max \Big\{\Im (\tau)\ :\
\det \widetilde L(\tau,\xi)=0\,,
\quad
\xi\in \RR^d, \quad |\xi|^2= \rho^2
\Big\}\,.
\]
Apply the Seidenberg-Tarski Theorem \ref{thm:seidenber} with real variables $\rho$, $\zeta = (\re \,\tau,\Im \,\tau, \xi)$ and polynomials $R(\rho,\zeta)=|\det \widetilde L\txi|^2 + (|\xi|^2-\rho^2)^2$, $S=0$ and $Q(\rho,\zeta)=\Im\tau$. Conclude that there is an $\alpha \ne 0$ and a rational $r$ so that
\[
E(\rho) \ =\ \alpha\,\rho^r\big(1+o(1)\big),\qquad
\rho \to \infty\,.
\]
To prove the result it suffices to prove that $\Im \tau$ is bounded, \textit{i.e.} to show that $r\le 0$.  Suppose on the contrary that $r>0$.

Given $\tau,\xi$
define  $k\in S^{d-1}$, $\rho\in \RR_+$, and, $\theta$ by
\[
k\ :=\
 \frac{\xi}{|\xi|}\,,
\qquad
\xi \ =\ \rho\,k,
\qquad
\theta \ :=\
\frac{\tau}{\rho}
\ .
\]
Choose sequences
$\tau(n)$, and $\xi(n)$
so that for $n\to \infty$,
\begin{equation}
\label{eq:badxi}
\det
\widetilde L(\tau(n),\xi(n))=0,
\qquad
\Im \tau(n) = \alpha \big(\rho(n)\big)^r(1+o(1))
\,.
\end{equation}

Write
\[
\tiltx
\ =\
\tilutx + B
\ =\
\rho\,
\Big(
\tilu(\theta,k) + \frac{1}{\rho}\,B
\Big)
\ =\
\rho
\Big(
\theta I_{2N\times 2N} +
\tilu(0,k) + \frac{1}{\rho}\,B
\Big)
.
\]
The matrix $\widetilde L(\tau,\xi)
$ is singular if and  only if $-\theta$ is an
eigenvalue of $\widetilde L_1(0,k) + \rho^{-1}B$.

For large $\rho$ this is a small perturbation of
$\tilu(0,k)$.  Choose $\mu>0$ so that
for $|k|=1$, the only eigenvalue of $\tilu(0,k)$
in the disk $|\theta|\le 2\mu$ is $\theta=0$.

Because of the strong well posedness of $L$,
there is a  uniformly independent basis of unit
eigenvectors for the eigenvalues
of $L_1(0,k)$  in $|\theta|\ge \mu$.
By part  $(iv)$ of Theorem \ref{prop:charvar}
there is a  uniformly independent basis of unit
eigenvectors for the eigenvalues of $\tilu(0,k)$ in $|\theta|\ge \mu$.

It follows that there is a $C_0$ so that for $\rho>C_0$ the eigenvalues of $\widetilde L_1(0,k) + \rho^{-1}B$ in $|\theta|>\mu$ differ from the corresponding eigenvalues of $\widetilde L_1(0,k)$ by no more than $C_0/\rho$.  In particular their imaginary parts are no larger than $C_0/\rho$. Therefore,  the corresponding eigenvalues $\tau=\rho\,\theta$ have bounded imaginary parts. Thus for $n$ large, $E(\rho(n))$ can be reached only for the eigenvalues $-\theta(n)$ which are perturbations of the eigenvalue $0$ of $\widetilde L_1(0,k(n))$.

Perturbation by  $O(1/\rho)$ of
the uniformly bounded family of
$dN\times dN$ matrices, $\widetilde L_1(0,k)$,
can move the eigenvalues by no more
than $O(\rho^{-\frac{1}{dN}})$.
Since the unperturbed eigenvalue is 0,
$|\theta(n)|\le C\,\rho(n)^{-1/dN}$, so
\[
|\tau(n)|\ \le\  C\rho(n)^{ 1 -\frac{1}{dN} },
\qquad
{\rm Im}\, \tau(n)
\ =\
\alpha\, \rho(n)^r(1+o(1))
\,,
\quad
\alpha\ne 0\,.
\]
Therefore $r\le 1-1/dN<1$ and
\[
\ds\prod_{j=1}^d (\tau(n) +\sigma_j)
\ =\
\tau(n)^d(1+ o(1)),
%\quad
%{\rm and,}
\quad
\xi_\ell(n)
\ds\prod_{j\ne \ell}(\tau(n) +\sigma_j)
=
\xi_\ell(n)\,\tau(n)^{d-1}(1+
o(1))
\,.
\]
Insert in identity
 \eqref{eq:det} to find
\[
\det L_1\Big(\tau(n)^d(1+o(1))\, ,\, \xi(n)\,\tau(n)^{d-1}(1+o(1))
\Big) = 0.
\]
Divide the argument by  $\rho(n)\,\tau(n)^{d-1}$
and use homogeneity to find
\[
\det L_1\bigg(
\frac{\tau(n)}
{\rho(n)}\big(1+o(1)\big)\,,\, k(n)(1+o(1))
\bigg) \ =\ 0.
\]

The constant multiplicity hypothesis shows that
\begin{equation}
\label{eq:F1}
\det L_1(\tau,\xi) \ =\ \tau^m\, F_1(\tau,\xi)\,,
\qquad
{\rm and}
\qquad
\forall \xi\in \RR^d,\quad
F_1(0,\xi) \ne 0.
\end{equation}
Since for $n$ large $(\tau(n)/\rho(n))(1+o(1))\ne 0$ we have
\[
F_1\bigg(
\frac{\tau(n)}
{\rho(n)}\big(1+o(1)\big)\,,\, k(n)(1+o(1))
\bigg) \ =\ 0.
\]
Passing to a subsequence we may suppose
that the bounded sequence $k(n)\to k$.
In addition, $\tau(n)/\rho(n)\to 0$ so
passing to the limit yields
$
F_1\big(0,k)= 0
$
contradicting \eqref{eq:F1}.
This contradiction proves $( i)$.\\

\noindent $( ii)$
Part $( i)$ of Theorem \ref{prop:charvar}
shows that  $0$ is an eigenvalue of $\tilo$ with algebraic
multiplicity equal to  $N(d-1)+m$. It remains to see that with the assumption, the dimension of $ \tiklo$ is strictly smaller than $N(d-1)+m$.
By definition
\[
\tiklo \ =\
\Big\{
\widetilde{\Phi}=(\Phi_1,\cdots,\Phi_d)\ :\  \ \ds \sum_{j=1}^d \Phi_j \in \cap_p \ker(\xi_pA_p)
\Big\}\,.
\]

Define
$$
{\cal E}_1\ :=\
\Big\{
\widetilde \Phi \ =\ \big(
\Phi_1,\cdots, \Phi_d\big)
\ :\
\sum_{j=1}^d \Phi_j = 0
\Big\}\,.
$$
Then ${\cal E}_1\subset  \tiklo  $ and ${\rm dim}\, {\cal E}_1=N(d-1)$.

 Define
 $$
 {\cal E}_2
 \ :=\
 \klo \otimes O^{d-1},
 \qquad
 {\rm dim}\, {\cal E}_2 \ =\ m.
 $$

 If $\widetilde \Phi\in \tiklo$, $\sum_{j=1}^d \Phi_j \in\klo $ , write
 \[
\tilde{\Phi}
\ =\
\Big(\sum \Phi_j\,,\,0\,,\,\dots\,,\,0\Big)
\ -\
W,
\qquad
\Big(\sum  \Phi_j\,,\,0\,,\,\dots\,,\,0\Big)\in
{\cal E}_2, \quad W\in {\cal E}_1.
\]
Thus,
$\tiklo \subset {\cal E}_1 \oplus {\cal E}_2$.

Pick $V$ in $\klo$, but not in $\underset{\xi_j\ne 0}{\cap} \ker A_j$. Then
$$
\widetilde V =(V,0,\cdots,0)\ \in\
{\cal E}_2
\quad
{\rm and}
\quad
V\notin\tiklo\,.
$$
This proves that $\tiklo$ is
a proper subset of
${\cal E}_1 \oplus {\cal E}_2$, so
$$
{\rm dim }\, \big(\tiklo\big)
\  <\
{\rm dim}\,
{\cal E}_1
\ +\
{\rm dim}\,
{\cal E}_2
\ =\
N(d-1) +m\,.
$$
Thus  the geometric multiplicity of the eigenvalue
$0$
is strictly less than its algebraic multiplicity.
Therefore,
$\widetilde L_1(0,\xi)$ is not diagonalizable.
This proves $( ii)$.\\

\noindent $(iii)$
To prove that the split problem is strongly well posed
it suffices to consider the principal part.
Suppose  $L(0,\xi)$ is uniformly
diagonalisable on a conic neighbourhood of
$\underline\xi \in \RR^N\setminus 0$.
For $\tiu=(U^1,\cdots,U^d)$, introduce
\begin{equation}
\label{eq:deftildeV}
\widetilde{V}= (V^1,\cdots,V^d)\quad {\rm with}\quad
  V^1:=\sum_{j=1}^d U^j, \quad{\rm and}\quad
  V^l:=U^l\quad {\rm for}\quad 2\leq l \leq d\,.
  \end{equation}
   Then
    \[
  \tilde{L}_1(\partial_t ,\partial_x)\tiu=0 \iff \partial_t \widetilde{V}+\tilde{Q}(\partial_x)\widetilde{V}=0,\quad
  \mbox{ with, }\quad
   \widetilde{Q}(\xi):=
  \begin{pmatrix}
     \lo \quad& 0 & \dots & 0 \\
      \xi_2 A_2 \quad&0&\dots&0 \\
      \vdots \quad& \vdots & \ddots & \vdots\\
      \xi_d A_d \quad&0&\dots&0
    \end{pmatrix} .
  \]

The eigenvalues of $\widetilde{Q}(\xi)$ are those of $L_1(0,\xi)$, therefore real. It suffices to diagonalize uniformly $\widetilde{Q}(\xi)$ on the conic neighbourhood of $\underline\xi$.  By homogeneity it suffices to consider $\xi$ with $|\xi|=1$.
By hypothesis there exist a real diagonal matrix
$D(\xi)$ and an invertible matrix
$S(\xi)$ so that
\[
  L_1(0,\xi)=S(\xi) D(\xi) S^{-1}(\xi)\,, \quad \textrm{ and, }\quad \exists K >0, \  \forall \xi \in \RR^d,\quad \|S(\xi)\|+\|S^{-1}(\xi)\| \leq K.
\]
%  $D(\xi)$ has the form
%  \[
% D(\xi)=\begin{pmatrix}
%      0      & 0           & \dots    & 0 \\
%      0      &\lambda_1(\xi)\,I & \dots    &0 \\
%      \vdots &  \vdots     & \ddots   & \vdots\\
%      0      &    0        & \dots    & \lambda_{m_1}(\xi)\,I
%        \end{pmatrix} .
% \]
%
Seek a diagonalization of $ \widetilde{Q}(\xi)$ on $|\xi|=1$ in the form,
\begin{equation}
\label{eq:deftildeS}
  \widetilde{S}(\xi)\ =\
        \begin{pmatrix}
            S(\xi) \quad                 &   0       & \dots     & 0 \\
            \xi_2 A_2 Q(\xi) S(\xi) \quad&   Id      &\dots      & 0 \\
            \vdots     \quad             &   \vdots  & \ddots    & \vdots \\
            \xi_d A_d Q(\xi) S(\xi) \quad&   0       &\dots      & Id
        \end{pmatrix} \quad
        {\rm so},
        \quad
        (\widetilde{S}(\xi))^{-1}\ =\
        \begin{pmatrix}
            (S(\xi) )^{-1}   \quad  &   0       & \dots     & 0 \\
            -\xi_2 A_2 Q(\xi) \quad &   Id      &\dots      & 0 \\
            \vdots            \quad &   \vdots  & \ddots    & \vdots \\
            -\xi_d A_d Q(\xi) \quad &   0       &\dots      & Id
        \end{pmatrix}
\end{equation}
%
%with
%\[
%Q(\xi)= S(\xi)
%    \begin{pmatrix}
%      0      & 0           & \dots    & 0 \\
%      0      &1/\lambda_1D(\xi)\,I & \dots    &0 \\
%      \vdots &  \vdots     & \ddots   & \vdots\\
%      0      &    0        & \dots    & 1/\lambda_{m_1}D(\xi)\,I
%    \end{pmatrix} .
%    S^{-1}(\xi)
%    \]
%
with $Q(\xi)$  to be determined. Then,
\[
\widetilde{S}^{-1}(\xi) \widetilde{Q}(\xi) \widetilde{S}(\xi)
\ =\
        \begin{pmatrix}
            D(\xi)          \quad        &   0       & \dots     & 0 \\
           \xi_2 A_2(I- Q(\xi)\lo)S(\xi) \quad &   0      &\dots      & 0 \\
            \vdots                 \quad &   \vdots  & \ddots    & \vdots \\
            \xi_d A_d(I- Q(\xi)\lo)S(\xi)\quad&   0       &\dots      & 0
        \end{pmatrix}
\]
and
\begin{equation}
\label{eq:qp}
\widetilde{S}^{-1}(\xi) \widetilde{Q}(\xi) \widetilde{S}(\xi)\text{\ \  is diagonal}
\quad
\Longleftrightarrow
\quad
 \xi_j A_j (I-Q(\xi)L_1(0,\xi))=0,
 \ \
 2\le j\le d
 \,.
 \end{equation}

 From the strong well posedness of $L_1$ it follows that uniformly in $\xi$
 one has $\ker\lo \oplus \im\lo = \CC^N$.
Choose $Q$ equal to the left inverse of $\lo$ defined in \eqref{eq:inversedef}. Since $\klo=\cap\, \ker \xi_jA_j$,
the condition  on the right in
\eqref{eq:qp} holds so $\widetilde S(\xi)$ diagonalizes $\widetilde Q(\xi)$.
Since
$S(\xi)$ and $S(\xi)^{-1}$ are bounded on a conic neighborhood,
it follows that
$\tilde{S}(\xi)$ and $\tilde{S}(\xi)^{-1}$ are  bounded on
a neighborhood of $\uxi$ in $|\xi|=1$.
A finite cover of the sphere, completes the proof.
\end{proof}

\begin{remark}
\label{rmk:conservednorm}
Denote by $\widetilde S(\xi)$ the function homogeneous
of degree zero given by
\eqref{eq:deftildeS}
for $|\xi|=1$ with  $Q$ constructed in the proof.
Then $\|\widetilde S(D)\widetilde V(t)\|_{L^2(\RR^d))}$
with $\widetilde V$ from \eqref{eq:deftildeV}
is a norm  equivalent to
$\|\widetilde U(t)\|_{L^2(\RR^d)}$
and is conserved for solutions of $\widetilde L_1(\partial)\widetilde U = 0$.
Those solutions yield a unitary group with respect to  the norm
$\|\widetilde S(D) \widetilde V\|_{L^2(\RR^d)}$.
\end{remark}

%%%%%%%%%%%%%%%%%%%%%%%%%%%
%  New section
%%%%%%%%%%%%%%%%%%%%%%%%%%
\section{Analysis  of the B\'erenger's PML by energy methods}
\label{sec:stability}
This section contains results
proving that the
initial  value problems
so defined  are
well posed.  We begin with the case of Gevrey
absorptions, then $W^{2,\infty}$, and finally the case
of the Heaviside function.

In Section
\ref{subsec:Gevrey}
we prove that when $\widetilde L$ is only weakly well posed,
Gevrey regular $\sigma_l$ lead to well
posed initial value problems in Gevrey classes.
Commonly used $\sigma$ are not this smooth.

The strongest result,
from Section
\ref{subsec:elliptic},
applies when
$L(0,\partial_x)$ is elliptic.
Important cases are the wave equation
and linearized elasticity. In these cases
the operator
$\widetilde L_1$ is strongly hyperbolic so
remains strongly hyperbolic even with general
bounded zeroth order perturbations.   Thus for
bounded $\sigma_l(x_l)$ the initial value problem
is strongly well posed.

In Section
\ref{subsec:maxwell2DTE}
we  analyse the case of
$\widetilde L$ associated to Maxwell's
equations with  finitely smooth $\sigma$.
  We follow
the lead of
\cite{Petit:2006:PFB}
and
extend the analysis of \cite{Metral:1999:CBP}
to several absorptions $\sigma_l$ and to higher
dimensions.
Related estimates for the linearized Euler
equation have been studied by
L. M\'etivier
\cite{Metivier:2009:UEE}.

The results of this section do not treat the case
of
B\'erenger's method for  Maxwell's equations with discontinuous
$\sigma_j$.   The case of one absorption is treated in
\S  4.  A closely related method is treated by an
energy
method in \S \ref{subsub:PMRB}.

\subsection{General operators and Gevrey absorption}
\label{subsec:Gevrey}
The next result
is implied by Bronstein's Theorem
\cite{Bronshtein:1979:SRP}
\cite{Bronshtein:1980:CPH} \cite{Nishitani:1983:EIN}\ \cite{Nishitani:1983:SEH}.
It
shows
that when $\lo$ has only real eigenvalues and the
$\sigma_j$ belong to the appropriate Gevrey class,
then the Cauchy problem for $\widetilde L$ is solvable for
Gevrey data.

\begin{definition}  For $1\le s<\infty$,
$f\in {\cal S}^\prime(\RR^d)$ belongs to the Gevrey
class $G^s(\RR^d)$ when
$$
\exists C,M, \ \
\forall\alpha\in \NN^d,
\ \
\big\|\partial^\alpha f\big\|_{L^2(\RR^d)}
\ \le\
M \,
\alpha !
\,
C^{|\alpha|}.
$$
\end{definition}
Then $G^s\subset\cap_\sigma H^\sigma(\RR^d)\subset C^\infty(\RR^d)$.
For $s>1$ the compactly supported elements of $G^s$ are dense.  If
$|\hat f(\xi)|\le C\,e^{-|\xi|^a}$ with $0<a<1$, then $u\in G^{1/a}$.

\begin{theorem}
If the principal part $L_1$ is weakly hyperbolic, and $\sigma_j\in G^{N/N+1}(\RR^d)$
then for arbitrary
$f\in G^{N/(N+1)}(\RR^d)$ there is one and only one solution
$u\in C^\infty(\RR^{1+d})$ to
$$
\widetilde L\, u \ =\ 0,
\qquad
u (0,\cdot)\ =\ f\,.
$$
The solution depends continuously on $f$.
\end{theorem}

\subsection{Strong hyperbolicity when
$L(0,\partial)$ is elliptic}
\label{subsec:elliptic}

\begin{theorem}
\label{thm:elliptic}
If $L$ is strongly well posed and $L(0,\partial)$ is elliptic, then $\til$ is strongly well posed for  any absorption $(\sigma_1(x_1), \dots,
  \sigma_d(x_d))$ in $(L^{\infty}(\RR))^d$.
\end{theorem}
\begin{proof}
Kreiss' theorem \ref{fort} asserts that an operator with  constant coefficient principal part is uniformly well posed if and only if the principal part is uniformly diagonalisable on a conic neighborhood of each $\uxi\ne 0$.  Therefore the Corollary  follows from the third part of Theorem \ref{th:thetheorem}.
\end{proof}

\begin{example}
This result implies that the PML model for the elastodynamic system is strongly well posed.
The system is written in the velocity-stress $(v,\Sigma)$ formulation,
\[
\rho\,\partial_t v - div \,\Sigma =0\,,\qquad
\partial_t \Sigma -C\eps(v)=0\,,
\qquad
\eps_{ij}(v):=(\partial_iv_j+\partial_j v_i)\,,
\]
with  positive definite elasticity tensor $C$ and  $\Sigma:=C\eps$. See \cite{Becache:2003:SPM}, where the authors showed that such  layers may be amplifying (see Section \ref{sec:amplification}).
\end{example}

\subsection{The method of M\'etral-Vacus extended to  the  $3d$ PML Maxwell system}
\label{subsec:maxwell2DTE}
M\'etral and Vacus proved in \cite{Metral:1999:CBP} a stability estimate for B\'erenger's two dimensional  PML Maxwell system with one absorption $\sigma_1(x_1)\in W^{1,\infty}(\RR)$ and $x=(x_1,x_2)\in \RR^2$. There are two crucial elements in their method.   First following B\'erenger, they do not split all variables in all directions. This section begins by showing that the partially split model is equivalent to the fully split model restricted to functions $\widetilde U$ some of whose components vanish.  The $\widetilde L$ evolution leaves this space invariant and its evolution on that subspace determines its behavior everywhere.

The second element is that on
the  partially split subspace there is an {\it a priori} estimate bounding the norm at time $t$ by the same norm at time $0$.  This looks inconsistent with the fact that the Cauchy problem is only weakly well posed. However the norm is not homogeneous.  Certain linear combinations of components have more derivatives estimated than others.  The observation of \cite{Metral:1999:CBP} is that the system satisfied by the fields and certain combinations of the fields and their derivatives, yields a large but symmetrizable first order system. These estimates have been obtained, and extended in Sabrina Petit's thesis \cite{Petit:2006:PFB} in the $2d$ case with two coefficients, and in the $3d$ case for an absorption in only one direction.

In this section,
motivated in part by the
clarification of the role of symmetrizers in
the work of
 L. M\'etivier \cite{Metivier:2009:UEE} for the $2d$ variable coefficient
 Euler equations in geophysics,
we construct analogous more elaborate functionals
which suffice for the  general case of three absorptions in three dimensions.
They require $\sigma_j\in W^{2,  \infty}(\RR)$.

Maxwell's equations for $\partial_t E_1$ and $\partial_t B_1$ contain only partial derivatives with respect to $x_2,x_3$ and not $x_1$.  In such a situation B\'erenger splits the corresponding equations in directions $x_2,x_3$
 but not in direction $x_1$.   To see why this is a special case of the general splitting algorithm \eqref{PMLgen} reason as follows. If the equation for $\partial_t U_j$ from $L$ does not contain any terms in $\partial_k$, that is the $j^{\rm th}$ row of   $A_k$ vanishes, then the equation for the $j^{\rm th}$ component of the unknown $U^k$ corresponding to the splitting for the $k^{\rm th}$ space variable is,

\begin{equation}
\label{eq:silentsplit}
\partial_tU^k_j\ +\
\sigma_k(x_k)\,U^k_j\ =\ 0\,,
\qquad
U^k_j=e^{-\sigma_k(x_k)t}\,U^k_j(0,x).
\end{equation}
Plugging this into the other equations reduces
the number of unknowns by one.
The simplest strategy is to take initial data
$U^k_j(0,x)=0$ which yields the operator $\widetilde L$
restricted to the invariant subspace of functions
so that
$U^k_j=0$.
Conversely if one knows how to solve that restricted
system then the full system can be reduced to the
restricted system with an extra source term from
\eqref{eq:silentsplit}.

\vskip.1cm
\noindent
{\bf Summary.}  {\sl To study the fully split system it is
sufficient to study the system restricted to
$\{U^k_j=0\}$. Performing this reduction for
each missing spatial derivative, corresponds
to splitting  equations only along directions
containing the corresponding spatial derivatives.}
\vskip.1cm

\noindent
An extreme case of this reduction occurs if an equation
contains no spatial derivatives, that unknown is eliminated
entirely.  For the Maxwell system which is the subject of this
section this does not occur.  The use of unsplit variables
\begin{itemize}
  \item reduces the size of $\widetilde U$ reducing computational
cost,
  \item corresponds to B\'erenger's original algorithm,
\item is important for the method of M\'etral-Vacus which takes advantage of the vanishing components $U^k_k$.
\end{itemize}

Consider the $3d$  Maxwell equations,
\[
\partial_tE-\nabla \times H=0, \qquad
\partial_tH+\nabla \times E=0.
\]
Defining $U=E+iH$, they take the symmetric hyperbolic form \eqref{hyp} with hermitian matrices
\begin{equation}
\label{matricemaxwell2D}
A_1
\ =\
\begin{pmatrix}
 0&0&0\\0&0&-i\\0&i&0
 \end{pmatrix}
,\quad
 A_2
\ =\
\begin{pmatrix}
0&0&i\\0&0&0\\-i&0&0
 \end{pmatrix}
 \quad\mbox{and}\quad
 A_3
\ =\
\begin{pmatrix}
0&-i&0\\i&0&0\\0&0&0
 \end{pmatrix}\,.
 \end{equation}

Introduce the splitting \eqref{PMLgen}
with some components unsplit. 
Define the subspace $\cal H$ with vanishing components
corresponding  to the unsplit components
\[
{\cal H}\ :=\
\Big\{\widetilde U=(U^1,U^2,U^3) \in H^2(\RR^3\,;\,\CC^3)^3\big\}\ :\ \
 U^1_1=0,
\quad U^2_2=0
,
\quad U^3_3=0
\Big\}.
\]
For $\widetilde U=(U^1,U^2,U^3)$ in ${\cal H}$, define
\begin{equation}\label{eq:defV}
\begin{array}{c}
U\ :=\
 U^1+U^2+U^3, \qquad V^j\ :=\ \pd j U,
 \qquad V^{i,j}\ : =\
 \partial_{ij} U,
 \qquad
 W\ :=\
  \sum_k \sig_k(x_k)U^k,
 \\[2mm]
 W^j\ :=\  \pd jW,
 \qquad
 Z\ :=\
 \sum_k\pd k(W_k+\sig_k(x_k)U_k),\qquad
 Z^{j}\ :=\
 \pd j Z,
 \\[2mm]
 \mathbb{V}
 \ :=\
 \Big(U \,,\, V^i,V^{i,j} \,,\, W^j \,,\, U^j \,,\, W \,,\, Z^{j}
 \Big) \ \in\
  \CC^{54}.
\end{array}
\end{equation}
The function $Z$ and therefore $Z^j$ are $\CC$ valued.  The other slots
in $\VV$
are $\CC^3$ valued.
The second derivatives $V^{i,j}$ of $U$ are ordered as
$V^{1,1}, V^{2,1},V^{3,1},V^{2,2},V^{3,2},V^{3,3}$.
This convention is important when the equations for $\VV$
are written in matrix form.  Computing in turn $W^j, Z, Z^j$
requires two derivatives of $\sigma_j$.

The unknown in \eqref{PMLgen} is $\widetilde U = (U^1,U^2,U^3)$.
The $U^j$ appear in the fifth slot  of $\VV$.  Therefore,
\[
\| \VV(t,\cdot)\|_{(L^2(\RR^3))^{54}} \ge   \| \widetilde U(t,\cdot) \|_{(L^2(\RR^3))^{9}}.
\]
 For the Cauchy problem the initial data is $\widetilde U_0= (U^1_0,U^2_0,U^3_0)$, from which $\VV_0$ is deduced by the derivations  above, and
\[
\| \VV_0\|_{(L^2(\RR^3))^{54}}
\le  C
\| \widetilde U_0  \|_{(H^2(\RR^3))^{9}}.
\]

\begin{theorem}
\label{th:vacus}
If $\sig_j$, for $j=1,2,3$, belong to $ W^{2,\infty}(\RR)$,
then for any $\widetilde U_0=(U^1_0,U_0^2,U^3_0)$ in ${\cal H}$
there is a unique solution $\widetilde U$ in $L^2(0,T;{\cal H})$ of the
split Cauchy problem \eqref{PMLgen} with initial value $\widetilde U_0$. Furthermore
there is a $C_1>0$ independent of $\widetilde U_0$ so that
for all positive time $t$,
\begin{equation}\label{eq:estMV}
  \|\widetilde U(t,\cdot) \|_{(L^2(\RR^3))^{9}}
\    \leq\
C_1e^{C_1\,t}
\,
\big\|\widetilde U_0\, \big\|_{(H^2(\RR^3))^{9}}\,
.
\end{equation}
\end{theorem}
%\begin{theorem}
%\label{th:vacus}
%If $\sig_j$, for $j=1,2,3$, belong to $ W^{2,\infty}(\RR^2)$, there exists
%$C$, such that, for any $\dot{U}=(\dot{U}^1,\dot{U}^2,\dot{U}^3)$ in ${\cal H}$
%there is a unique solution of the doubled Cauchy problem with initial value $\dot{U}$ \textbf{***discuter la formulation avec la donnée initiale ***} satisfying
%for all positive time $t$,
%%
%\[
%  \|\mathbb{V}(t,.)\|_{L^2}
%\    \leq\
%\|\mathbb{V}(0,.)\|_{L^2}\,
%e^{C\,t} \,.
%\]
%\end{theorem}
%
\begin{proof}   The main step is to derive a system of equations
satisfied by $\VV(t,x)$ together with a symmetrizer $S(D)$.  These
imply an estimate for $t\ge 0$,
\begin{equation}
\label{eq:VVest}
\big\|
\VV(t)
\big\|_{L^2(\RR^3)}
\ \le \
C_2\,e^{C_2t} \,
\big\|
\VV(0)
\big\|_{L^2(\RR^3)}
\,.
\end{equation}
From this estimate it easily follows that the Cauchy problem for the
$\VV$-equations is uniquely solvable.
It is  true but not immediate that if the initial values of $\VV$
are computed from those of $\widetilde U$ then the solution $\VV$
comes from a solution $\widetilde U$ of the B\'erenger system.
 The strategy
has three steps,

\vskip.1cm

$\bullet$  Discretize the B\'erenger system in $x$ only.

$\bullet$ Derive an estimate analogous to \eqref{eq:VVest} for the
semidiscrete problem.   The estimate is uniform as the discretization parameter
tends to zero.   The proof is a semidiscrete analogue
of \eqref{eq:estMV}.

$\bullet$ Solve the semidiscrete problem and pass
to   the limit to prove existence.

\vskip.1cm
\noindent
This is done for the case $d=2$ in
\cite{Petit:2006:PFB}  to which we refer for details.
Uniqueness of the solutions
to the $\VV$-system and therefore $\widetilde U$ is simpler and classical
and is also in \cite{Petit:2006:PFB}.

Equation \eqref{PMLgen} yields,

\begin{equation}\label{eq:Uj}
\pd t U^j+ \sum_k A_k V^j+\sig_j U^j =0.
\end{equation}
Summing on $j$ yields
\begin{equation}\label{eq:U}
\pd t U+L(0,\partial)U+W=0.
\end{equation}
Differentiate in direction $x_j$  to find,
\begin{equation}\label{eq:Vj}
\pd t V^j+L(0,\partial) V^j+W^j=0.
\end{equation}
Differentiate once more to get
\begin{equation}\label{eq:Vij}
\pd t V^{i,j}+L(0,\partial) V^{i,j}+\pd iW^j=0.
\end{equation}
The quantity  $\partial_iW^j$
on the left  is replaced using the next lemma.

\begin{lemma}\label{lem:intermWj}
\begin{equation}\label{eq:intermWj}
\begin{array}{lcl}
\partial_j W &=& L(0,\partial)A_jW+Ze_j -\sum_kE_{jk}( \sig_k' U+ \sig_k V^k),
\end{array}
\end{equation}
where $e_j$ is the $j^{\rm th}$ vector of the standard basis, and $E_{ij}$ is the $3\times 3$ matrix all of whose entries vanish except the  $(i,j)$ element that  is equal to $1$.
\end{lemma}
\begin{proof}
%\begin{proof}[Proof of Lemma \ref{lem:intermWj}]

First evaluate $L(0,\partial)A_jW$ to find
\[
L(0,\partial)A_jW= \sum_k A_kA_j \pd k W.
\]
The matrices in  Maxwell's equations
 satisfy $A_j A_k = -E_{jk}$ for $j \ne k$, and
$A_j^2= \sum_{k\ne j}E_{kk}=I-E_{jj}$. This yields
\[
L(0,\partial)A_jW = -\sum_{k\ne j}E_{jk} \pd k W +(I-E_{jj})\pd j W
=\pd j W -\sum_{k}E_{jk} \pd k W
=\pd j W -div(W)\,e_j.
\]

Introduce the definition of $Z$ to find
\[
\begin{array}{lcl}
W^j
\ :=\
\partial_j W &=& L(0,\partial)A_jW+div(W)e_j
\\
             &=& L(0,\partial)A_jW+Ze_j -(\sum_k\pd k (\sig_k U_k))\, e_j.
             \end{array}
\]
Compute
\[
\begin{array}{lcl}
\big( \sum_k\pd k (\sig_k U_k) \big)\, e_j&=&
      \sum_kE_{jk}\pd k(\sig_k U)\\
&=&
\sum_kE_{jk}\sig_k' U
+ \sum_kE_{jk}\sig_k V^k,
\\
\end{array}
\]
which proves \eqref{eq:intermWj}.  The proof of the lemma is complete.
\end{proof}

Differentiate \eqref{eq:intermWj} in space to obtain
\[
\begin{array}{lcl}
\pd i W^j &=& L(0,\partial)A_jW^i +\pd iZe_j -
\sum_kE_{jk}\pd i\big(  \sig_k U+ \sig_k'V^k  \big).\\
\end{array}
\]
Inserting into \eqref{eq:Vij} yields
\[
\partial_t V^{i,j}+L(0,\partial) V^{i,j}+
L(0,\partial)A_jW^i +Z^{ij} -
\sum_kE_{jk}\pd i\big( \sig_k U+ \sig_k'V^k  \big)=0.
\]
This is equivalent to
\begin{equation}\label{eq:Vijbis}
\begin{array}{l}
\partial_t V^{i,j}
\ +\
L(0,\partial) V^{i,j}
\ +\
L(0,\partial)A_jW^i\ +
 \\
\hspace{15mm}
Z^{i}e_j \ -\
\sig_i'E_{ji}U
\ -\
(\sig_i"E_{ji}
\ +\
 \sum_k\sig_k E_{jk}) V^i
\ -\
\sum_k \sig_k' E_{jk} V^{i,k}
\ =\ 0.
\end{array}
\end{equation}
To close the system it remains to evaluate the time derivatives of $W, W^j$ and $\pd jZ$.
\[
\pd t W
\ =\
\sum_k \sig_k \pd t U^k
       \ =\
       -\sum_k \sig_k (A_k U+\sig_k U^k)
       \ =\
       -\sum_k \sig_k A_k V^k -\sum_k \sig_k^2 U^k.
       \]
Using the particular form of the equations yields
\[
\begin{array}{rcl}
\sum_k \sig_k^2 U^k&=&
(\sum_k \sig_k)(\sum_k \sig_k U^k)-
\sum_k\sig_k\left(\sum_{l\ne k}\sig_l U^l
\right)\\
&=&
(\sum_k\sig_k) W
\ -\ diag(\sig_2 \sig_3,\sig_1\sig_3,\sig_1\sig_2)U
\ -\ diag(\sig_1,\sig_2,\sig_3)W\,.
\end{array}
\]
Therefore
\begin{equation}\label{eq:W}
\pd t W
+\sum_k\sig_k A_k V^k
+(\sum_k\sig_k) W
-diag(\sig_2 \sig_3,\sig_1\sig_3,\sig_1\sig_2)U
-diag(\sig_1,\sig_2,\sig_3)W=0\,.
\end{equation}
Differentiate  in $x_i$ to find
\begin{equation}\label{eq:Wj}
\begin{array}{l}
\pd t W^i+
\sum_k(\pd i(\sig_k) A_k V^k+\sig_k A_k V^{ki} )
+\pd i(\sum_k\sig_k) W+ \sum_k\sig_k  W^i\\[2mm]
\hspace{20mm}
-\pd i\big(diag(\sig_2\sig_3,\sig_1\sig_3,\sig_1\sig_2)\big)U
-diag(\sig_2\sig_3,\sig_1\sig_3,\sig_1\sig_2)V^i\\[2mm]
\hspace{35mm}
-\pd i\big(diag(\sig_1,\sig_2,\sig_3)\big)W
- diag(\sig_1,\sig_2,\sig_3)W^i
=0
\end{array}
\end{equation}
Next  compute
\[
\pd t Z=
\sum_i \pd i\pd t
 (W_i+\sig_iU_i)\,.
\]
Consider the pair of  equations
\[
\pd t U^2_1 +i \pd 2 U_3 +\sig_2 U^2_1 =0,
\qquad
{\rm and}
\qquad
\pd t U^3_1 -i \pd 3 U_2 + \sig_3 U^3_1 =0.
\]
Add the two equations.
Also add  $\sig_2$  times the first to $\sig_3$ times the second.
This yields two equations,
\[
\pd t U_1 +i(\pd 2 U_3 - \pd 3 U_2) + W_1 =0,
\quad
{\rm and}
\quad
\pd t W_1 +i(\sig_2\pd 2 U_3 - \sig_3\pd 3 U_2) +
\sig_2^2 U^2_1 +\sig_3^2 U^3_1 =0.
\]
Rewrite the last term as
$
\sig_2^2 U^2_1 +\sig_3^2 U^3_1=(\sig_2+\sig_3)W_1-\sig_2\sig_3U_1,
$
to find
\[
\pd t U_1 +i(\pd 2 U_3 - \pd 3 U_2) + W_1 =0,
\quad
{\rm and,}
\quad
\pd t W_1 +i(\sig_2\pd 2 U_3 - \sig_3\pd 3 U_2) +
(\sig_2+\sig_3)W_1-\sig_2\sig_3U_1 =0.
\]
Multiply the first equation by $\sig_1$ and add the second to obtain
\[
\begin{array}{l}
\pd t (W_1+\sig_1U_1)
+i\left((\sig_1+\sig_2)\pd 2 U_3 - (\sig_1+\sig_3)\pd 3 U_2\right)
+(\sum_k\sig_k )W_1 -\sig_2\sig_3U_1=0\,.
\end{array}
\]
The other indices follow by permutation.  Differentiate in $x_k$ and add
to find
\[
\begin{array}{l}
\hspace{-10mm}
\sum_k\pd k \pd t (W_k+\sig_kU_k)
+i\sum_k\pd k\left((\sig_k+\sig_{k+1})\pd {k+1} U_{k+2}
- (\sig_k+\sig_{k+2})\pd {k+2} U_{k+1}\right)\\[2mm]
\hspace{30mm}
+\sum_i\pd i\left((\sum_k\sig_k )W_i\right)
-\sum_k\pd k(\sig_{k+1}\sig_{k+2} U_k)=0.
\end{array}
\]
The terms with two spatial derivatives cancel.
This leaves
\[
\begin{array}{l}
%\hspace{-10mm}
 \pd t Z%\\
%\hspace{15mm}
+i\sum_k\sig_k'( \pd {k+1} U_{k+2}- \pd {k+2} U_{k+1})
\\[3mm]
\hspace{20mm}
+(\sum_k\sig_k )(\sum_kW_k^k)
+\sum_k \sig_k'W_k
-\sum_k \sig_{k+1}\sig_{k+2} V^k_k=0.
\end{array}
\]
Since
\[
 Z =\sum_k\pd k(W_k+\sig_k(x_k)U_k)
 = \sum_k (W_k^k+\sig_k'U_k+\sig_kV_k^k),
\]
we can replace $(\sum_k\sig_k )(\sum_kW_k^k)$ in the previous equation by $(\sum_k\sig_k )(Z-\sum_k (\sig_k'U_k+\sig_k\,V_k^k))$
so
\[
\begin{array}{l}
%\hspace{-10mm}
 \pd t Z +(\sum_k\sig_k )Z%\\
%\hspace{15mm}
+i\sum_k\sig_k'( V^{k+1}_{k+2}- V^{k+2}_{k+1})
\\[3mm]
\hspace{20mm}
-(\sum_k\sig_k )(\sum_k (\sig_k'U_k+\sig_k\,V_k^k))
+\sum_k \sig_k'W_k
-\sum_k \sig_{k+1}\sig_{k+2} V^k_k=0.
\end{array}
\]
Differentiating in $x_j$ yields
\[
\begin{array}{l}
%\hspace{-10mm}
 \pd t  Z^j +(\sum_k\sig_k )Z^j+ \sig_j'Z
 %\\
%\hspace{15mm}
+i \sig_j"( V^{j+1}_{j+2}- V^{j+2}_{j+1})
+i\sum_k\sig_k'( V^{k+1,j}_{k+2}- V^{k+2,j}_{k+1})
\\[3mm]
\hspace{2mm}
- \sig_j' (\sum_k (\sig_k'U_k+\sig_k\,V_k^k))
-(\sum_k\sig_k )( \sig_j"U_j+\sig_j'\,V_j^j)
-(\sum_k\sig_k )(\sum_k (\sig_k'V^j_k+\sig_k\,V_{j,k}^k))\\[3mm]
\hspace{6mm}
+ \sig_j'W_j
+\sum_k \sig_k'W^j_k
-\sum_k \pd j(\sig_{k+1}\sig_{k+2}) V^k_k
-\sum_k \pd j(\sig_{k+1}\sig_{k+2}) V^{j,k}_k
=0.
\end{array}
\]
Replace $Z$ by $\sum_k(W^k_k+\sig_kV^k_k)$ to  end up with
\begin{equation}\label{eq:Zj}
\begin{array}{l}
%\hspace{-10mm}
 \pd t Z^j +(\sum_k\sig_k )Z^j+ \sig_j'\sum_k(W^k_k+\sig_kV^k_k)
 %\\
%\hspace{15mm}
+i \sig_j"( V^{j+1}_{j+2}- V^{j+2}_{j+1})
+i\sum_k\sig_k'( V^{k+1,j}_{k+2}- V^{k+2,j}_{k+1})
\\[3mm]
\hspace{2mm}
- \sig_j' (\sum_k (\sig_k'U_k+\sig_k\,V_k^k))
-(\sum_k\sig_k )( \sig_j"U_j+\sig_j'\,V_j^j)
-(\sum_k\sig_k )(\sum_k (\sig_k'V^j_k+\sig_k\,V_{j,k}^k))
+ \sig_j'W_j\\[3mm]
\hspace{4mm}
+\sum_k \sig_k'W^j_k
-\sum_k \pd j(\sig_{k+1}\sig_{k+2}) V^k_k
-\sum_k \pd j(\sig_{k+1}\sig_{k+2}) V^{j,k}_k
=0.
\end{array}
\end{equation}

Summarizing, $\mathbb{V}$ is solution of a first order
system,
$\partial_t \mathbb{V} + P(\partial_x)\mathbb{V} + B(x) \mathbb{V}=0$, whose principal symbol is given by
\[
P(\partial)=
\begin{pmatrix}
I_4\otimes L(0,\partial) & 0_{4,6}\otimes 0_{3,3}& 0_{4,3}\otimes 0_{3,3}  & 0_{4,3}\otimes 0_{3,3} & 0_{4,4}\otimes 0_{3,3}  \\[2mm]
0_{6,4}\otimes 0_{3,3} & I_6\otimes L(0,\partial)& (I_6\otimes L(0,\partial)) M & 0_{6,3}\otimes 0_{3,3}& 0_{6,4}\otimes 0_{3,3}\\[2mm]
0_{3,4}\otimes 0_{3,3} & 0_{3,6}\otimes 0_{3,3}& 0_{3,3}\otimes 0_{3,3}& 0_{3,3}\otimes 0_{3,3} & 0_{3,4}\otimes 0_{3,3}  \\[2mm]
0_{3,4}\otimes 0_{3,3} & 0_{3,6}\otimes 0_{3,3}& 0_{3,3}\otimes 0_{3,3}& 0_{3,3}\otimes 0_{3,3}& 0_{3,4}\otimes 0_{3,3}  \\[2mm]
0_{4,4}\otimes 0_{3,3} & 0_{4,6}\otimes 0_{3,3}& 0_{4,3}\otimes 0_{3,3}& 0_{4,3}\otimes 0_{3,3}& 0_{4,4}\otimes 0_{3,3}  \\[2mm]
\end{pmatrix}\,.
\]
Here the $V^{i,j}$ are ordered as indicated before the theorem and,
\[
M
\ :=\
\begin{pmatrix}
A_1 & 0   & 0   \\
0   & A_1 & 0 \\
0   & 0   & A_1 \\
0   & A_2 & 0   \\
0   & 0   & A_2 \\
0   & 0   & A_3 \\
\end{pmatrix}\,.
\]

To symmetrize it suffices to construct  a symmetrizer for the upper left hand block
\[
Q(\partial)\ :=\
\begin{pmatrix}
I_4\otimes L(0,\partial) & 0_{4,6}\otimes 0_{3,3}& 0_{4,3}\otimes 0_{3,3}  \\[2mm]
0_{6,4}\otimes 0_{3,3} & I_6\otimes L(0,\partial)& (I_6\otimes L(0,\partial)) M \\[2mm]
0_{3,4}\otimes 0_{3,3} & 0_{3,6}\otimes 0_{3,3}& 0_{3,3}\otimes 0_{3,3} \end{pmatrix}
\,.
\]
We verify that
\[
\tilde{S}\ :=\
\begin{pmatrix}
I_4\otimes I_3 & 0_{4,6}\otimes 0_{3,3}& 0_{4,3}\otimes 0_{3,3} \\
0_{6,4}\otimes 0_{3,3} & I_6\otimes I_3& (I_6\otimes I_3) M \\
0_{3,4}\otimes 0_{3,3} & 0_{3,6}\otimes 0_{3,3}& I_3\otimes I_3
\end{pmatrix}
\quad {\rm with}\quad
\tilde{S}^{-1}=
\begin{pmatrix}
I_4\otimes I_3 & 0_{4,6}\otimes 0_{3,3}& 0_{4,3}\otimes 0_{3,3}   \\
0_{6,4}\otimes 0_{3,3} & I_6\otimes I_3& - M(I_3\otimes I_3)  \\
0_{3,4}\otimes 0_{3,3} & 0_{3,6}\otimes 0_{3,3}& I_3\otimes I_3\\
\end{pmatrix}\
\]
 is a symmetrizer for $Q(i\xi)$. Compute
\[
\tilde{S}Q\tilde{S}^{-1}\ =\
\begin{pmatrix}
I_4\otimes L(0,\cdot) & 0_{4,6}\otimes 0_{3,3}& 0_{4,3}\otimes 0_{3,3}  \\
0_{6,4}\otimes 0_{3,3} & I_6\otimes L(0,\cdot)& 0_{6,4}\otimes 0_{3,3}  \\
0_{3,4}\otimes 0_{3,3} & 0_{3,6}\otimes 0_{3,3}& 0_{3,3}\otimes 0_{3,3}   \\
\end{pmatrix}
\]
which is symmetric since $L(0,\cdot)$ is.

Therefore $P(\xi)$ is
symmetrizable by a matrix independent of $\xi$.
Hence,  the Cauchy problem for \eqref{eq:U}, \eqref{eq:Vj}, \eqref{eq:Vijbis}, \eqref{eq:W}, \eqref{eq:Wj}, \eqref{eq:Zj} is
strongly well posed.
The norm of the zero order  terms depends on the coefficients $\sig_j$ and their derivatives up to order 2.  The estimate  of the Theorem follows.

%The existence and uniqueness of solutions satisfying the estimate
%can be proved by replacing $\partial_k$ in the equations
%by $j_\epsilon * \partial_k * j_\eps *$ with $j_\eps$ a symmetric
%approximation of $\delta$.  The solutions of the regularized equations
%satisfy analogous  estimates uniform in $0<\eps<1$.  One concludes on passing
%to the limit.
% Uniqueness follows follows from uniqueness
% of solutions of the Cauchy problem for the
% extended system.
 \end{proof}
\begin{remark}
We recall a computation from \cite{Petit:2006:PFB},
showing that when there  are only 2 coefficients $\sigma_1$ and $\sigma_2$
only one derivative of $\sigma_j$ is needed.
This is always the case in dimension $d=2$.
When  $\sig_3 \equiv 0$, split $W$ as
\[
W=E_{33}W + \begin{pmatrix} W_1\\W_2\\0\end{pmatrix}\,.
\]
Then
\[
\begin{pmatrix} W_1\\W_2\\0\end{pmatrix}=
\begin{pmatrix} \sig_2U^2_1\\\sig_1U^1_2\\0\end{pmatrix}=
diag(\sig_2,\sig_1,0)U.
\]
Rewrite \eqref{eq:U} as
\begin{equation}\label{eq:Ubis}
\pd t U+L(0,\partial)U+E_{33}W  + diag(\sig_2,\sig_1,0)U=0.
\end{equation}
Differentiate with respect to $x_1$ and $x_2$ to obtain
\begin{equation}\label{eq:Vjbis}
\pd t V^j+L(0,\partial) V^j+
E_{33}\pd j W
+\pd j(diag(\sig_2,\sig_1,0))U +
diag(\sig_2,\sig_1,0)W
=0.
\end{equation}
To find an equation on $W$, proceed as in the $3d$ proof to  get,
\begin{equation}\label{eq:Wbis}
\pd t W
+\sum_k\sig_k A_k V^k
+(\sum_k\sig_k) W
- \sig_1\sig_2 U=0
\end{equation}
Therefore $\mathbb{V}$ is solution of a first order
system, whose principal symbol is given by
\[
P(\partial)=
\begin{pmatrix}
L(0,\partial)& 0  & 0  & 0 \\
 0  & L(0,\partial)& 0  & E_{33}\pd 1 \\
 0  & 0  &L(0,\partial)& E_{33} \pd 2\\
 0  & 0  & 0 & 0
\end{pmatrix}
\ .
\]
A symmetrizer is given by
\[
\tilde{S}\ =\
\begin{pmatrix}
\ I\     & 0  & 0  & 0 \\
 0  & \ I\   & 0  & i E_{23} \\
 0  & 0  &\ I\    & i  E_{13}\\
 0  & 0  & 0 & I
\end{pmatrix}\,,
\qquad
{\rm with}
\qquad
\tilde{S}^{-1}\ =\
\begin{pmatrix}
\ I\    & \ 0\   & \ 0\   & 0 \\
 0  & I   & 0  & -i  E_{23} \\
 0  & 0  & I   & -i  E_{13}\\
 0  & 0  & 0  & I_3
\end{pmatrix}
\,.
\]

\end{remark}

\subsection{Sharp Finite Speed for B\'erenger's PML}\label{sec:finitespeed}

Recall some notions associated with estimates
on the domains of influence and determinacy
for a hyperbolic operator $L$ (see \cite{Joly:2005:HDD}).
The {\it timelike cones} are the connected components
of $(1,0,\dots,0)$ in the complement of
the characteristic variety. The {\it forward propagation cone}
is dual to the  timelike cone.  A lipschitzean curve
$[a,b]\ni t\mapsto (t,\gamma(t)$ is an influence curve
when $\gamma^\prime$ belongs to the propagation
cone for Lebesgue almost all $t\in [a,b]$.

\begin{theorem}
\label{thm:finitespeed}
  Suppose that $L$ defines a strongly well
posed Cauchy problem and that the multiplicity of
$\tau=0$ as a root  of $\det L_1(\tau,\xi)=0$
is independent of $\xi\in \RR^d\setminus 0$.  The support of the solution of the B\'erenger
transmission problem is contained in the union of the propagation
curves of $L$ starting in the support of the source terms when
either of the following conditions is satisfied.

$(i)$
$
\forall \xi\in \RR^d\setminus 0,
\quad
{\rm Ker}\,L_1(0,\xi) \ =\ \cap_j\,{\rm Ker}\, \xi_jA_j
$, and $\forall j,\ \sigma_j\in L^\infty$.

$(ii)$   $L_1$ is Maxwell's equation and $\forall j,\ \sigma_j\in W^{2,\infty}$.
\end{theorem}

\begin{proof}  The characteristic varieties of $L$ and $\widetilde L$ satisfy
${\rm Char}\,\widetilde L={\rm Char}\,L\cup\{\tau=0\}$.
 When $\{\tau=0\}$ has
multiplicity as a root of $\det L_1(\tau,\xi)=0$
independent of $\xi\in \RR^d\setminus 0$ the timelike
cones of $L$ and $\widetilde L$ coincide.
Therefore the propagation cones and influence curves
coincide too.

{\bf Case $(i)$.}
Part $(ii)$ of Theorem \ref{th:thetheorem}   proves that $\widetilde L_1$  defines a strongly
well posed Cauchy problem.
It follows that the sharp propagation conclusion of the Theorem is
valid for $\widetilde L_1 + B(t,x)$ for any bounded $B(t,x)$.  This
follows on remarking that the solution of $(\widetilde L +B)\widetilde U=0$
with initial data $\widetilde U_0$ is the limit at $\nu\to \infty$ of Picard iterates
$\widetilde U^\nu$.  The first, $\widetilde U^1$,
is defined as the solution of the Cauchy problem without
$B$.  For $\nu>1$ the iterates are defined by,
$$
\widetilde L_1 \widetilde U^{\nu+1} \ +\
B(t,x)\,\widetilde U^{\nu}\ =\ 0,
\qquad
\widetilde U^{\nu+1}(0,\cdot)\ =\ \widetilde U_0\,.
$$
Since $\widetilde L_1$ has constant coefficients, sharp finite
speed is classical for that operator.
An  induction proves that each iterate is supported in the union of influence
curves starting in the support of $\widetilde U_0$.

{\bf Case $(ii)$.}    Reason as above constructing by Picard
iteration approximations $\VV^\nu$ converging to the solution
 $\VV$ from \eqref{eq:defV}.  Since the equation satisfied by $\VV$ is strongly well
 posed the iterates converge.
An induction shows that they
are supported in the set of influence curves starting in the support of
$U_0$.
\end{proof}

\subsection{Proof of perfection for B\'erenger's PML by a change of variables}
\label{subsec:changevariable}
This section continues the analysis of B\'erenger's method
when the hypotheses of Theorem \ref{thm:finitespeed} are satisfied.
In those cases well posedness is proved by an energy
method. In addition suppose that
\begin{equation}
\label{eq:Lj}
\forall j,\ \ \exists L_j>0, \qquad
\sigma_j=0 \quad {\rm when}\quad |x_j|\le L_j\,.
\end{equation}
Denote by $R:=\Pi_j]-L_j,L_j[$.

\begin{definition}   In this setting the method is perfectly matched
when for arbitrary ${\widetilde F}\in C^\infty_0\big(]0,\infty[\times R\big)$
the unique solutions $\widetilde V$ and $\widetilde U$ of
\begin{equation}
\label{eq:2sols}
\widetilde L \widetilde V \ =\ {\widetilde F},
\quad
\widetilde V\big|_{t\le 0} \ =\ 0,
\qquad
\widetilde L_1 \widetilde U \ =\ {\widetilde F},
\quad
\widetilde U\big|_{t\le 0} \ =\ 0,
\end{equation}
with $\widetilde L$ as in \eqref{PMLgen}
satisfies
\begin{equation}
\label{eq:rectangularperfection}
\widetilde V\big |_{\RR\times R}
\ =\
 \widetilde U \big |_{\RR\times R}\,.
 \end{equation}
 \end{definition}

 \begin{theorem}
 \label{thm:perfectionCOV}
 With the assumptions of Theorem \ref{thm:finitespeed}
 and $\sigma_j$ as above, B\'erenger's method is
 perfectly matched.
 \end{theorem}

 \begin{proof} Taking the Laplace transform of the $\widetilde V$
 equation in \eqref{eq:2sols} yields a transform holomorphic
 in ${\rm Re}\,\tau>\tau_0$ with values in $L^2(\RR^d)$
 satisfying for $1\le j\le d$,
 \begin{equation}
 \label{eq:jtransform}
\widehat {V^j}
  \ +\
  \ds (\tau + \sigma_j(x_j) )^{-1} A_j \partial _{j}
  \widehat{V}
  \ =\
  \widehat{  F^j},
\quad
{\rm with}
\quad   V:=\sum_j   V^j,
\quad
  F:= \sum_j  F^j .
\end{equation}
Multiply by $\tau$ and sum on $j$ to obtain
\begin{equation}\label{eq:sumed1}
\tau\, \widehat V
  \ +\
  \ds \sum_j
  \frac{\tau}{\tau+ \sigma_j(x_j)} \ A_j\, \partial _{j}\widehat V
  \ =\
\tau\,  \widehat F\, .
\end{equation}

When $\tau$ is {\it fixed real and positive} this equation can be
transformed to the corresponding equation without
the $\sigma_j$ by a change of variables.   The change of
variables depends on $\tau$.   The resulting
equation is exactly that determining
$\widehat U:= \sum_j\widehat{  U^j}$.   In this way
we find that $\widehat V$ is obtained from $\widehat U$ by
this change of variables. This idea is inspired by Diaz and Joly in \cite{Diaz:2006:TDA}.

For real $\tau >0$ define
$d$ bilipschitzean homeomorphisms  $X_j(x_j)$ of $\RR$ to itself by
$$
\frac{dX_j(x_j)}{dx_j}
\ =\
\frac{\tau+\sigma_j(x_j)}{\tau}
\,,
\qquad
X_j(0)\ =\
0\,.
$$
Then,
$$
\frac{\partial}{\partial x_j}  \ =\
\frac{\partial X_j}{\partial x_j}\,
\frac{\partial}{\partial X_j}
\ =\
\frac{\tau + \sigma_j(x_j)}{\tau}\
\frac{\partial}{\partial X_j},
\qquad
  \frac{\tau}{\tau+ \sigma_j(x_j)} \, \frac{\partial} {\partial  x_{j}}
  \ =\
  \frac{\partial}{\partial X_j}
  \,.
$$
Therefore if $\widehat U(X)$ is the solution of
\begin{equation}\label{eq:sumed2}
\tau\, \widehat U(X)
  \ +\
  \ds \sum_j
   \,A_j\frac{\partial}{ \partial X_{j}}\widehat U
  \ =\
 \widehat F(X)\,,
\end{equation}
then the solution $\hat V$ of \eqref{eq:sumed1} is given by
$\hat V(x):= \hat U(X(x))$
since the latter function of $x$ satisfies the equation
determining $\hat V$.

Since $X(x)=x$ for
$x\in R$ this proves that the transforms of $\widetilde U$
and $\widetilde V$
satisfy for real $\tau>\tau_0$
\begin{equation}
\label{eq:almostperfect}
\sum_j \widehat{ V^j}(\tau,x)
\ =\
\sum_j \widehat{ U^j}(\tau,x)\,,
\qquad
x\in R\,.
\end{equation}
Since both sides of
the identity in  \eqref{eq:almostperfect}
are holomorphic in ${\rm Re}\,\tau>\tau_0$
it follows
that  the identity extends to that
domain by analytic continuation.

Equation \eqref{eq:jtransform} and its analogue for
$\widetilde U$ then imply
that for all $j$, $\widehat{V^j}\big|_R=\widehat{ U^j}\big|_R$.
Uniqueness of the Laplace transform implies
$ V^j\big|_R= U^j\big|_R$ for all $t$
proving perfection.
\end{proof}

\begin{remark}
The proof  is very general.  It shows that once
the initial value problem defined by $\widetilde L$ is
well posed there is perfect matching.  The proof works
more generally
for at least weakly well posed  methods for which the Laplace transform can be
reduced to \eqref{eq:sumed1} for real $\tau$.   Our
favorite version of the B\'erenger algorithm is analysed this way in
\S \ref{subsub:PMRB}.
\end{remark}

\subsection{Perfection for methods related to B\'erenger's PML}
\label{subsub:PMRB}

Consider
\eqref{eq:sumed1} with $F=0$.
This equation is the starting point for many authors to
construct well posed PML.
It has been viewed as a complex stretching of coordinates (see \cite{Rappaport:1995:IIA}, \cite{Chew:1994:IIA}, \cite{Petropoulos:2000:RSL},  \cite{Hesthaven:1998:ACP}). This idea, for $\tau$ real, becomes an honest change of variables as in \cite{Diaz:2006:TDA}, that is at the heart of the proof in \S \ref{subsec:changevariable}. In the case of Maxwell system, it can be viewed as a system with modified constitutive equations (a lossy medium \cite{Petropoulos:1998:AMD},\cite{Abarbanel:1998:WPM}), or recovered as above from the B\'erenger's system.
The system \eqref{eq:sumed1}
is not differential because of the division by
$\tau +\sigma_j (x_j)$.  In order to recover a hyperbolic system, a change of unknowns is performed. We adopt  the approach in \cite{Mazet:1998:IDM} for the Maxwell system.

\begin{lemma}\label{lem:decompmatrix}
 With matrices given in \eqref{matricemaxwell2D},
define $S_j:= (\tau+ \sigma_j(x_j))/\tau$. There exists a pair of invertible
matrices $M,N$, unique up to a multiplication by the same constant,
such that
 \begin{equation}\label{eq:decompmatrix}
S_j^{-1}N A_j
\ =\ A_jM,\qquad  j=1,2,3.
\end{equation}
They are given by
\begin{equation}\label{eq:decompmatrices}
M=\gamma \begin{pmatrix}
S_1  &  0    &  0  \\
0    &  S_2  &  0  \\
0    &  0    &  S_3
\end{pmatrix},\qquad
N= \gamma
\begin{pmatrix}
S_2S_3  &  0    &  0  \\
0    &  S_1S_3  &  0  \\
0    &  0    &  S_1S_2
\end{pmatrix}\,,
\qquad
\gamma\in \CC\setminus 0\,.
\end{equation}
\end{lemma}
\begin{proof}
Since
\[
A_j e_j=0,\qquad
A_je_{j+1}=-ie_{j+2},\qquad
A_je_{j+2}=ie_{j+1},
\]
it is easy to see by applying \eqref{eq:decompmatrix} to $e_j$ that $M$ is necessarily diagonal,
$M=\mbox{diag}(m_1,m_2,m_3)$. Applying \eqref{eq:decompmatrix} to $e_{j+1} $ and $e_{j+2} $
shows that for any $j$,
\[
Ne_{j+1}=m_{j+2}S_{j} e_{j+1},\qquad Ne_{j+2}=m_{j+1}S_{j} e_{j+2}\,.
\]
This implies  that $N$ is also diagonal, equal to
$\mbox{diag}(m_2S_3,m_3S_1,m_1S_2)$, and
\[
m_1S_3=m_3S_1,\quad
m_2S_1=m_1S_2,\quad
m_3S_2=m_2S_3
\]
This leaves no choice but to choose \eqref{eq:decompmatrices}.
\end{proof}

\vskip.2cm

In the rest of the analysis take $\gamma=1$.
In \eqref{eq:sumed1} with $F=0$ replace
$\widehat V$ replace by $U$.
Insert
\eqref{eq:decompmatrix}
  to obtain

\begin{equation}\label{eq:sumedmod}
\tau N U
  \ +\
  \ds \sum_j A_j M \partial _{j}U
  \ =\
  0   .
\end{equation}

The fact
that  $\sig_j$ depends only on $x_j$
and the form of the matrices guarantees $A_j\partial_jM=0$.
This yields
\[
A_j \,M\, \partial _{j}U
\ \equiv\
 A_j \,\partial _{j}\,(M\,U)
 \,.
\]
Define a new unknown  $V:=MU$ to find
\begin{equation}\label{eq:sumedmod2}
N M^{-1}\tau V
  \ +\
  \ds \sum_j A_j \partial _{j}V
  \ =\
  0   .
\end{equation}
$N M^{-1}=\ds\mbox{diag}
(S_1^{-1}S_2S_3,S_2^{-1}S_3S_1,S_3^{-1}S_1S_2)$.
Next compute a rational fraction expansion of
$\tau S_1^{-1}S_2S_3$ as
\[
\ds\frac{(\tau+\sig_2)(\tau+\sig_3)}{\tau+\sig_1}\ =\
\tau \ +\
 (\sig_2+\sig_3-\sig_1)\  +\
\frac{ \sig_1^2+\sig_2 \sig_3-\sig_1(\sig_2+\sig_3)}{\tau+\sig_1}\,.
\]
 Introduce a new unknown $W$ by
\[
\tau N M^{-1} V \ =\
\tau V + \Sigma^1 V + \Sigma^2 W,\qquad
{\rm equivalently}
\qquad
W_j
 \ =\
 \ds \frac{1}{\tau + \sig_j(x_j)} V_j=
 \frac{1}{\tau} U_j,
\]
with
\[
\begin{array}{l}
\Sigma\ :=\
\mbox{diag}\big(\sig_1,\sig_2,\sig_3\big),
\\
\Sigma^{(1)}\  :=\
\mbox{diag}\big(
\sig_2+\sig_3-\sig_1,
\sig_3+\sig_1-\sig_2,
\sig_1+\sig_2-\sig_3 \big)
,
\
\\
\Sigma^{(2)}\  :=\
\mbox{diag}\big(
(\sig_1-\sig_2)(\sig_1-\sig_3),
(\sig_2-\sig_1)(\sig_2-\sig_3),
(\sig_3-\sig_1)(\sig_3-\sig_2)
\big).
\end{array}
\]
This leads to a system in the unknowns $V$ and $W$
\begin{equation}
	\label{eq:rahmouni}
L(\partial_t, \partial_x) V
\ +\  \Sigma^{(1)}V \ +\ \Sigma^{(2)} W\ =\  0\,,
\qquad
\partial_t W\ +\ \Sigma W \ -\
V\ =\  0\,.
\end{equation}
Finally,  $U$ is recovered from
$$
U\ =\ \partial_t W\ =\
V\ -\ \Sigma W\,.
$$
The system of equations for $V,W$  is strongly well posed since $L$ is symmetric
hyperbolic.
In the case of a single layer in the $x_1$ direction,
there is  only one coefficient $\sig$ and
therefore a single complex valued suplementary variable.
The equations for the magnetic and electric fields are
\[
\begin{array}{l}
\partial_t E_1 -(\nabla \wedge H)_1 -\sig E_1 +\sig^2W_1=0,\\
\partial_t E_2 -(\nabla \wedge H)_2 +\sig E_2=0,\\
\partial_t E_3 -(\nabla \wedge H)_3 +\sig E_3=0,\\
\partial_t W_1 +\sig W_1 -E_1=0,
\end{array}
\quad\ \
\begin{array}{l}
\partial_t H_1 +(\nabla \wedge E)_1 -\sig H_1 +\sig^2W_2=0,\\
\partial_t H_2 +(\nabla \wedge E)_2 +\sig H_2=0,\\
\partial_t H_3 +(\nabla \wedge E)_3 +\sig H_3=0,\\
\partial_t W_2 +\sig W_2 -H_1=0,
\end{array}
\]
In $2d$
this is identical to the layers in \cite{Petropoulos:1998:AMD} and equivalent to
those in \cite{Abarbanel:1998:WPM}.

%\subsection{Abarbanel-Gottlieb}
%We analyse the 2-D Maxwell equation.  The
% analysis in 3-D should be done as well. Begin
%  with Abarbanel-Gottlieb (Appl Num Maths 1998)
%They refer to Ziolkowski: introduce a current $J$ in the equation, call $P=J+\sig E_x$, and obtain
%\[
%\left\{
%\begin{array}{l}
%\pdt{E_x}-\pd yH \textcolor[rgb]{1.00,0.00,0.00}{-\sig E_x +P}=0,\\
%\pdt{E_y}+\pd xH \textcolor[rgb]{1.00,0.00,0.00}{+\sig E_y}=0,\\
%\pdt{H}+\pd xE_y -\pd yE_x  \textcolor[rgb]{1.00,0.00,0.00}{+\sig H}=0,\\
%\textcolor[rgb]{1.00,0.00,0.00}{\pdt{P}+\sig P -\sig^2 E_x=0}.
%\end{array}
%\right.
%\]
%\textcolor[rgb]{0.98,0.00,0.00}{Z-A-G is strongly well posed, perfect and absorbing for Maxwell 2D.
%}
The principal symbol and lower terms are
\[
R_1=
\begin{pmatrix}
L& 0 \\
0& I_{3}\partial_t
\end{pmatrix}
\qquad
{\rm and}
\qquad
B=
\begin{pmatrix}
\Sigma^{(1)}& \Sigma^{(2)} \\
-I_{3} & \Sigma
\end{pmatrix}\,.
\]
Reversing
the computation %backward
shows that $(V,W) \in \ker R(\tau,\xi)$ if and
only if $V=MU$, $W=1/\tau U$, and
$L_1(\tau,\frac{\xi_1\tau}{\tau+\sig_1}, \cdots,\frac{\xi_3\tau}{\tau+\sig_3})U=0$. The characteristic polynomial is therefore the same as for B\'erenger's layer.   Thus, by Theorem \ref{prop:charvar}
\[
\det R(\tau,\xi)
\ =\
\tau^2-\sum \frac{\xi_j^2\tau^2}{(\tau+\sig_j)^2}.
\]

\begin{theorem}  If $\sigma_j(x_j)\in L^\infty(\RR)$ and vanish
for $|x_j|\le L_j$ then the system \eqref{eq:rahmouni} for
$V:=MU$ and $W_j:=V_j/(\tau+\sigma_j(x_j))$ is strongly
well posed in $L^2(\RR^d)$ and perfectly matched in the sense
that for sources supported in $R:=\Pi_j]\{|x_j|\le L_j\}$
the function $U$ computed from $V,W$ agrees in
$\RR_t\times R$ with the
solution of Maxwell's equation with corresponding sources.
\end{theorem}

\begin{proof}  The proof of Theorem \ref{thm:finitespeed}
applies with only minor
modifications.
\end{proof}

\begin{remark}  Note the ease with which strong well posedness is established
and the lack of regularity required
of the functions $\sigma_j$.
\end{remark}

\section{Analysis of layers with only one absorption by Fourier-Laplace transform}

There are cases where the energy method presented
above does not prove well posedness.   This is the case
for the B\'erenger algorithm
when the ellipticity assumption is not satisfied and the
absorptions are not regular.  Notably for the Maxwell
system and discontinuous absorptions.
In this section we present a systematic analysis by
Fourier-Laplace transformation of transmission problems
with absorption in only one direction.

\subsection{Fourier analysis of piecewise
constant coefficient transmission problems}
\label{subsec:piecewiseconstant}

Return to the situation
of \eqref{eq:LR} with operators $L$ and $R$
on the left and right half spaces and transmission
condition \eqref{eq:caN}.  Suppose that both
$L$ and $R$ are weakly hyperbolic in the sense
of G\aa rding.
An example is the classical method of B\'erenger
with one absorption.   Among other things we will
prove that the method is well posed and perfect.
Note the open problem at the end of the introduction
emphasizing that we do not know if the classic
algorithm with two discontinuous absorptions is well
posed.
In addition, we show,  by a non trivial analytic continuation
argument
in \S \ref{sec:continuation},
 that the perfection of B\'erenger's method can
be verified using the modified plane wave solutions
from his original paper.
It is our hope that the analysis may help in  the
construction
of new perfectly matched layers.

\subsubsection{Hersh's condition for transmission problems}
\label{sec:Hershcond}

This section takes up the analysis of  mixed problems
following  Hersh in \cite{Hersh:1963:MPS}.
In the present context we treat transmission problems
which are essentially equivalent.
The analysis of Hersh supposed the interface
is noncharacteristic which is never the case
for Maxwell's equations.
 We
address the changes that are needed to treat
problems with characteristic  interfaces.

First analyse the solution of the constant coefficient
pure initial value problem
$\cl U=F$ on $\RR^{1+d}$ by Laplace transform in time and Fourier
transform in $x'=(x_2,\cdots,x_d)$.
The  transform
$$
\widehat U(\tau, x_1, \eta)
\ :=\
\int
\int_0^\infty\
e^{-\tau t}
\
(2\pi)^{-d/2}\,
e^{-ix'\cdot\eta}
\
U(t, x')\
dt\, dx'
$$
 decays as $|x_1|\to \infty$ and satisfies
$$
\cl (\tau,d/dx_1,i\eta)\,\widehat U
\ =\
\widehat F
\quad
{\rm in}
\quad
\RR\,.
$$
When $A_1$ is invertible this is a standard ordinary differential
equation in $x_1$.  When $A_1$ is singular,
 the analysis
requires care.
The homogeneous equation $
\cl (\tau,d/dx_1,i\eta)\,\widehat U
=0$ has purely exponential solutions
$e^{\rho x_1}$ corresponding to the roots
$\rho$ of the equation
\begin{equation}
\label{eq:rho}
\det \cl (\tau\,,\,\rho\,,\, i\eta)\ =\ 0\,.
\end{equation}
Hyperbolicity of $\cl$
guarantees that for ${\rm Re}\,\tau>\tau_0$
and $\eta\in \RR^{d-1}$ this equation has no purely imaginary roots.

The number of boundary conditions at $x_1=0$ for the boundary value problem in the right half space is chosen equal to the number of roots with negative real part (see also Remark \ref{rmk:bcnum}). That  integer must be independent of $\tau,\eta$.  Since roots cannot cross the imaginary axis, the only way the integer can change is if roots escape to infinity. That can happen when the coefficient of the highest power of $\rho$ vanishes.  The next hypothesis rules that out.

\vskip.2cm

\begin{definition}   A hyperbolic operator $\cl(\partial_t,\partial_x)$
is {\bf nondegenerate with respect to $x_1$}
when there is a $\tau_1>0$ so that
the degree in $\rho$ of the polynomial
$
\det \cl(\tau,\rho,i\eta)
$
is  independent of $(\tau,\eta)$
for ${\rm Re}\,\tau>\tau_1\,,\, \eta\in \RR^{d-1}$.
\end{definition}

\begin{example}
\label{ex:nondeg}
\newline
{\bf 1.} In the noncharacteristic
case, $\det A_1\ne 0$,
the condition is satisfied and  the degree
with respect to $\rho$ is equal to $N$.

\noindent
{\bf 2.}   For Maxwell's  equations
written in the real $6\times 6$ form,
the degree with respect to $\rho$ is
equal to 4. If written in the complex form \eqref{matricemaxwell2D}, the degree is 2.

\noindent
{\bf 3.}   The formula for the characteristic polynomial
in Theorem \ref{prop:charvar} shows that
if $L$ is nondegenerate then so is the
B\'erenger doubled operator $\widetilde L$ with
one
absorption $\sigma_1$ in $x_1>0$.
The degree in $\rho$ is the same for $L$ and
$\widetilde L$.

\noindent
{\bf 4.}  If $\cl=\cl_1+B$ is nondegenerate with respect to $x_1$
then so is the operator
$P:=\cl_1(\partial) + a^{-1}B = a^{-1}\cl(a\partial)$ for any $a>0$.
If the degree for $\cl$ is constant in ${\rm Re}\,\tau>\tau_1$,
then the degree for $P$ is constant for
${\rm Re}\,\tau> a^{-1}\tau_1$.
\end{example}

For the lemmas to follow it is useful to transform so that ${\cal A}_1$ has block form.

\begin{lemma}
\label{lem:rankA}
If $\cl$ in \eqref{gen} is nondegenerate with respect to $x_1$ then for
${\rm Re}\,\tau>\tau_1\,,\, \eta\in \RR^{d-1}$,
\begin{enumerate}
[label={\it(\roman{*})}, ref={\it(\roman{*})},leftmargin=0.7cm]
\item
the degree in $\rho$ of the polynomial
$\det \cl(\tau,\rho,i\eta)$
is equal to $\,{\rk}\,{\cal A}_1$,
\item the  number of roots
$\rho$ with positive real part is equal to the number
of negative eigenvalues of ${\cal A}_1$.
\end{enumerate}
\end{lemma}

\begin{proof}
\newline
Since $\cl$ is nondegenerate, it suffices to study the case $\eta=0$. \\
\noindent $(i)$
Choose invertible $K$ so that
$$
K^{-1}\, {\cal A}_1\, K\ =\
\begin{pmatrix}
\caA & 0
\cr
0 & 0
\end{pmatrix}
,
\qquad
\caA \ \text{an invertible square matrix of size} \ \rk {\caA}_1\,.
$$
Then
$$
K^{-1} \cl(\tau,\rho, 0)\,K
\ =\
\begin{pmatrix}
\tau I + \caA \,\rho& 0
\cr
0 & \tau\,I
\end{pmatrix}
\ +\
{\rm matrix  \ independent \ of \ }\tau,\rho\,.
$$
It follows that the degree in $\rho$ is no larger than
${\rk}\, {\cal A}_1$.

The coefficient of
$\rho^{ {\rk}\caA_1}$
in $\det \cl (\tau,\rho,0)$  is a polynomial in $\tau$ of degree
$\,\le N-{\rk}\, \caA_1$.
For large $\tau$ the coefficient  is equal to
$$
\big(\det \caA\big)\ \big(\tau^{N-{\rk}\,{\cal A}_1}\big) + {\rm lower \ order\ in \ }\tau
\,.
$$
Thus the degree in $\rho$ is $\,{\rk}\,{\cal A}_1$ for such $\tau$
proving the result.\\

\noindent $(ii)$
\[
\det \cl(\tau,\rho,0)=
\det\left(
\begin{pmatrix}
\tau I+\rho \caA & 0\\
0 & \tau I
\end{pmatrix}
+B
\right)\,.
\]
For fixed $\tau$, sufficiently large, $\rho\ne 0$, and
$\rho/\tau$ is a root of the polynomial $p(x,1/\tau)$
of degree $\rk \caA$:
\[
p(x,\eps)=
\begin{vmatrix}
 I+x\caA+\eps B_{11} & \quad \eps B_{12}\\
\eps B_{21} & \quad  I+\eps B_{22}
\end{vmatrix}
\]
This polynomial has exactly $\rk \caA$ roots. By Rouch\'{e}'s theorem,
\[
\frac{\rho_j}{\tau}\sim  -\frac{1}{\lambda_j}, \quad \tau \gg 1,
\]
 where the $\lambda_j$ are the $\rk\caA$ eigenvalues of
 $\caA$ repeated according to their algebraic multiplicity.
 Since the eigenvalues of $\caA$ and the nonzero eigenvalues
 of ${\cal A}_1$ are the same,
 this completes the proof.

\end{proof}

\begin{remark}
\label{rmk:bcnum}
For the transformed one dimensional hyperbolic operator
$\cl(\partial_t, \partial_{1}, i\eta)$ the number of incoming characteristics at the boundary $x_1=0$ in the right half space is equal to the number of strictly positive eigenvalues of ${\cal A}_1$.  The second part of the Lemma shows that this is equal to the number of roots with negative real part.  The two natural ways to compute the number of necessary boundary conditions yield the same answer.
\end{remark}

The next lemma shows that for nondegenerate operators, the characteristic case can be transformed to a standard ordinary differential equation.

\begin{lemma}
\label{lem:blochform} Suppose
that $A,M\,\in\,{\rm Hom}\,(\CC^N)$
and  the equation
$\det(A\rho +M)=0$ has degree in $\rho$ equal to
$\rk\, A$ and no purely imaginary roots.
Then,
\begin{enumerate}
[label={\it(\roman{*})}, ref={\it(\roman{*})},leftmargin=0.7cm]
\item   The matrix $M$ is invertible
and all solutions of the homogeneous
equation
\begin{equation}
\label{eq:singode}
A\,
\frac{dU}{dx_1}
\ +\
M\,U\ =\ 0
\end{equation}
take values in the space
$\GG:= M^{-1}(\im\,A)$ satisfying
 $\, {\dim}\,\GG={\rk}\,A$.
\item  There is a $\widetilde M \in {\rm Hom}\,\GG$
so that a function $U$ satisfies
\eqref{eq:singode} if and only if $U$ is $\GG$ valued and
satisfies
\begin{equation}
\label{eq:Mtilde}
\frac{dU}{dx_1}
\ +\
\widetilde M\, U \ =\ 0\,.
\end{equation}
\item    The vector space $\mathbb U$ of solutions
 of
\eqref{eq:singode}
is a  linear subspace
 of $C^\infty(\RR)$ with dimension equal to
 ${\rk}\,A$.
The Cauchy problem with data in $\GG$ is
well posed.
\end{enumerate}
\end{lemma}

\noindent
\begin{proof}
\newline\noindent $(i)$  Since $\rho=0$ is not a root,  $M$
is invertible.   The equation
$
U\ = \
  M^{-1}A \,
dU/dx_1
$
shows that continuously differentiable solutions
$U$ takes values in $\GG$.
More generally, if
$U$ is a distribution solution and
$\psi\in C^\infty_0(\RR)$
takes values in the annihilator, $\GG^\perp$  of $\GG$, then
$$
\big\langle
U\,,\, \psi
\big\rangle
\ =\
\big\langle
M^{-1}\, A \, dU/dx_1
\,,\,
\psi
\big\rangle
\ =\
\big\langle
dU/dx_1
\,,\,
(M^{-1}\,A)^*
\psi
\big\rangle
$$
But
$\GG = {\rm range}\, M^{-1}A$
so
$\GG^\perp= {\rm ker}\, (M^{-1}\,A)^*$.
Therefore $(M^{-1}\,A)^*\psi=0$
so
$\langle U,\psi\rangle=0$ which is the desired
conclusion.\\

\noindent $(ii)$
Multiplying the equation  by an invertible
$P$ and making the change of variable
$U=KV$ transforms  the equation
to the equivalent form
$$
P\, A\,K\,
\frac{dV}{dx_1}
\ +\
P\,M\,K\,V
\ =\ 0\,.
$$
Choose invertible $P,K$ so that
$PAK$ has block form
$$
P\, A\, K\,
\ =\
\begin{pmatrix}
I & 0\cr 0 & 0
\end{pmatrix},
$$
where $I$ is the ${\rk}\, A\,\times\, {\rk}\, A$ identity
matrix.
With $V=(V_1,V_2)$ one has  the block forms
$$
P\, M\, K\,
\ =\
\begin{pmatrix}
H_{11}&H_{12}\cr
H_{21}&H_{22}
\end{pmatrix}
\qquad
{\rm and}
\qquad
\begin{pmatrix}
I & 0\cr 0 & 0
\end{pmatrix}
\frac{dV}{dx_1}
\ +\
\begin{pmatrix}
H_{11}&H_{12}\cr
H_{21}&H_{22}
\end{pmatrix}
V\ =\ 0\,.
$$
One has.
$$
\det (A\rho+M) \ =\
\det P^{-1}\
\det
\begin{pmatrix}
\rho I + H_{11}& H_{12}\cr
H_{21} & H_{22}
\end{pmatrix}
\
\det K^{-1}.
$$
The first part of the preceding lemma implies that the determinant on the left is a polynomical of degree ${\rk}\, A$ in $\rho$. It follows that $H_{22}$ is invertible.

The solutions $V$ satisfy $H_{21}V_1+H_{22}V_2=0$ so take values in $\VV:=\big\{V_2=-H_{22}^{-1}H_{21}V_1\big\}$. The function $V$ is a solution if and only if it takes values in $\VV$ and
$$
\frac{dV_1}{dx_1}
\ +\
RV_1\ =\ 0,
\qquad
R
\ :=\
H_{11}  - H_{12} H_{22}^{-1}H_{21}
\,.
$$
If $N:\VV\to \VV$ is the map,
$$
\big(V_1,V_2\big)
\ \mapsto\
\Big(
RV_1 \,,\ -H_{22}^{-1}H_{21}RV_1
\Big)\,,
$$
then $V$ is a solution if and only if it is $\VV$ valued
and satisfies $dV/dx_1=NV$.
Writing $V=K^{-1}U$
and $\widetilde M = -K\,N$ implies $(ii)$.\\

\noindent $(iii)$  Follows from $(ii)$.
\end{proof}

\begin{lemma}
\label{lem:L_1}
 If $\cl$ is hyperbolic and
nondegenerate with respect to $x_1$,
then
its principal part
$\cl_1(\partial_t,\partial_x)$
is also nondegenerate with respect to $x_1$.
The degree in $\rho$ of
$\det\,\cl_1(\tau,\rho,i\eta)$ is constant for ${\rm Re}\,\tau>0$
and $\eta\in \RR^{d-1}$.
\end{lemma}

\noindent
\begin{proof}
With notation from the preceding proof,
\begin{equation}
\label{eq:conjugate}
\cl(\partial_t, \partial_{x_1},\partial_{x^\prime})
\ =\
P^{-1}
\bigg(
\begin{pmatrix}
I\partial_{x_1}  & 0
\cr
0 & 0
\end{pmatrix}
\ +\
\begin{pmatrix}
H_{11}
(\partial_t, \partial_{x^\prime})
 &
 H_{12}
 (\partial_t, \partial_{x^\prime})
 \cr
 H_{21}
(\partial_t, \partial_{x^\prime})
 &
 H_{22}
 (\partial_t, \partial_{x^\prime})
\end{pmatrix}
\bigg)
K^{-1}\,.
\end{equation}
The proof of the last lemma showed that
for ${\rm Re}\,\tau > \tau _1$ and
$\eta\in \RR^{d-1}$,
$H_{22}(\tau,\eta)$
is invertible.

The computation in
Lemma \ref{lem:rankA}   shows that for $\eta=0$
and $\RR\ni\tau\to \infty$ the coefficient
of $\rho^{{\rk}\,{\cal A}_1}$ has
modulus $\ge c\tau^{N-{\rk}\,{\cal A}_1}$
with $c>0$.
This implies that  $\eta=0$ is noncharacteristic for $H_{22}$.
Therefore
 $H_{22}(\partial_t,\partial_{x^\prime})$
is hyperbolic.

Replacing $\cl$ by its principal part $\cl_1$ has the effect of
replacing each operator $H_{ij}(\partial)$ by its principal
part.  This yields identity \eqref{eq:conjugate} with
$\cl$ and the $H_{ij}$ replaced by their principal parts.

Since the principal part of a hyperbolic operator
is hyperbolic, it follows that
$(H_{22})_1(\partial_t,\partial_{x^\prime})$ is
a homogeneous  hyperbolic operator.  Therefore
$(H_{22})_1(\tau,i\eta)$ is invertible for $\eta\in \RR^{d-1}$
and ${\rm Re}\,\tau\ne 0$.   Thus, the
coefficient of $\rho^{{\rk}\,{\cal A}_1}$ in
$\det \cl_1(\tau, \rho, i\eta)$ is
nonzero
for $\eta\in \RR^{d-1}$
and ${\rm Re}\,\tau\ne 0$.
\end{proof}

\begin{lemma}
\label{lem:transdata}
Suppose that the ordinary differential equation
\eqref{eq:singode}
satisfies  the hypotheses of Lemma
\ref{lem:blochform}.
 Denote by $E^\pm$ the linear space of solutions
 which tend exponentially to zero as $x_1\to \pm \infty$
 and by $\dot E^\pm$ their traces at $x_1=0$.
 Then
\begin{enumerate}
[label={\it(\roman{*})}, ref={\it(\roman{*})},leftmargin=0.7cm]
\item   $ \dot E^\pm \cap {\rm ker}\,A\ =\
\{0\}\,,$
\item  ${\rm dim}\, A E^\pm \ =\
 {\rm dim}\,E^\pm\,,$
\item   The map $U\mapsto U(0)$ is an
 isomorphism from $E^\pm$ to $\dot E^\pm$\,,
\item    $A \dot E^+ \cap A\dot E^-\ =\
 \{0\},$
\item   $A \dot E^+ \oplus A\dot E^-\ =\
{\im}\,A$.
\end{enumerate}
\end{lemma}

\begin{proof}
\newline\noindent $(i)$
The absence of purely imaginary roots shows that every solution is uniquely the sum of two solutions. One  grows exponentially at $+\infty$ and decays exponentially at $x_1=-\infty$. The second grows at $-\infty$ and decays at $+\infty$.   In particular the only bounded solution is the zero solution.

 If $\bfe_+\in \dot E^+\cap {\rm ker}\, A$,  denote by $U(x_1)$ the solution
with this Cauchy data.
The function that is equal to $U$ on $x_1>0$ and equal
to 0 in $x_1\le 0$ is a distribution solution of
\eqref{eq:singode}
on all of $\RR$ since $A[U]_{x_1=0}=0$.
This solution is bounded hence identically equal to zero.
Therefore $\bfe_+=0$.  The case
for $\dot E^-\cap {\rm ker}\,A$ is analogous.\\

\noindent $(ii)$
Follows from $(i)$.\\

\noindent $(iii)$
It is surjective by definition.  If it were not
injective for $E^+$,
there would be a nontrivial solution $U(x)$ exponentially
decaying as $x_1\to +\infty$ with $U(0)=0$ violating $(i)$.\\

\noindent $(iv)$
The set $A\dot E^+$ consists of the values
$AU_+(0)$ with $U_+$ satisfying \eqref{eq:singode}
and exponentially decreasing in $x_1>0$.
If the intersection were nontrivial there would
be a solutions $U_-$ decaying as
$x_1\to -\infty$ so that $A U_+(0) = AU_-(0)$.
The function $V$ equal to $U_+$
in $x_1>0$ and $U_-$ in $x_1<0$ is then a
distribution solution
for all $x_1$ exponentially decaying in both directions.
Hyperbolicity implies that $V=0$ contradicting the
nontriviality.\\

\noindent $(v)$
 Using $(ii)$ and $(iv)$,
one sees that the direct sum on the left is
a subspace of ${\im}\,A$ of full dimension.
\end{proof}

The next lemma is needed in
in \S 4.1.2.

\begin{lemma}
\label{lem:seekF}
Assume that the hypotheses and notations of
Lemma
\ref{lem:blochform}
are in force.  Then for  $K\in \dot E^+$ there
is an $F\in C^\infty_0(]-\infty,0[)$ so that the
unique solution of
\begin{equation}
\label{eq:Mtildeinhomog}
A\,\frac{dU}{dx_1}
\ +\
M\, U \ =\ F\,,
\quad
\lim_{|x_1|\to \infty} \
\|U(x_1)\| \ =\ 0\,,
\end{equation}
satifies $U(0)=K$.
\end{lemma}

\begin{proof}
  Consider first the case of invertible $A$.
A change of dependent variable yields the block form
for the new variable still denoted $U$
$$
\frac{dU}{dx_1}
\ +\
\begin{pmatrix}
M_+ & 0 \cr 0 & M_-
\end{pmatrix} U
\ =\ F,
\quad
U=(U_1,U_2),
\quad
F=(F_1,F_2)\,,
\quad
{\rm spec}\,M_\pm\subset\{\pm {\rm Re}\,z>0\}.
$$
Then $\dot E^+=\{U_2=0\}$ so $K=(K_1,0)$.  Choose
$F=(F_1,0)$.  Then $U(0)=K$ if and only if,
$$
K_1\ =\
\int_{-\infty}^0 e^{M_+ s}\, F_1(s) \ ds.
$$
This is achieved with,
$$
F_1(s)\ =\
\chi(s)\,
e^{-M_+s} K_1,
\qquad
\chi\in C^\infty_0(]-\infty,0[),
\qquad
\int \chi(s)\ ds =1\,.
$$

When $A$ is not invertible change variable as in
Lemma
\ref{lem:blochform}
to find the block form
$$
\begin{pmatrix} I & 0\cr 0&0
\end{pmatrix}
\frac{dU}{dx_1}
\ +\
\begin{pmatrix}
H_{11} & H_{12} \cr H_{21} & H_{22}
\end{pmatrix}U
\ =\ F,
$$
with invertible $H_{22}$.

Part $(i)$ of
\eqref{eq:singode}
implies that the map $\GG\ni G = (G_1,G_2)\mapsto G_1$
is an isomorphism.
Write $\GG\ni K = (K_1,K_2)$.  Choose $F=(F_1,0)$.  Then choose
a
$\GG$ valued solution $U$ defined by
$$
\frac{dU_1}{dx_1}
\ +\
H_{11} U_1 \ =\ F_1,
\qquad
U_2\ =\ -\,H_{22}^{-1} H_{21}U_1\,.
$$
One has $U(0)=K$ if and only if $U_1(0)=K_1$.
The construction in the invertible case finishes the proof.
\newline

\end{proof}

Suppose that
$$
L= \partial_t + A_1\partial_1 +\cdots
\qquad
{\rm and}
\qquad
R =  \partial_t + \caA_1\partial_1 +\cdots
$$
are nondegenerate with respect to $x_1$.
For ${\rm Re}\,\tau>\tau_0$
and $\eta\in \RR^{d-1}$, define $E^\pm_L(\tau,\eta)$ to be the set of solutions
of
$$
L(\tau, d/dx_1,i\eta)V\ =\ 0
$$
tending to zero
as $x_1\to\pm\infty$.  Denote by
$\dot E^\pm_L(\tau,\eta)\subset \CC^N$ the linear space
of traces at $x_1=0$ of solutions  in $E^\pm_L(\tau,\eta)$.

Similarly with a possibly larger value still called $\tau_0$, there are
$\dot E^\pm_R(\tau,\eta)\subset \CC^M$
so that the solutions of
$R(\tau,d/d{x_1},i\eta)Z=0$ taking values in
$\dot E^\pm_R(\tau,\eta)$ are exactly those tending to zero exponentially as
$x_1\to \pm \infty$.
The subspaces $E^\pm_L(\tau,\eta)$ and $E^\pm_R(\tau,\eta)$ depend smoothly
on $\tau,\eta$ for ${\rm Re}\,\tau>\tau_0$ and $\eta\in \RR^{d-1}$.

Lemma \ref{lem:rankA} implies that
\begin{equation}
\label{eq:dimensions}
{\rm dim}\, E^-_L(\tau,\eta) \ =\
\#\, {\rm positive\ eigenvalues\ of\ } A_1,
\qquad
{\rm dim} \, E^+_R(\tau,\eta) \ =\
\#\, {\rm negative\ eigenvalues\ of\ } \caA_1\,.
\end{equation}
Consider the inhomogeneous transmission problem,
\begin{equation}
\label{eq:LRhomogeneous}
LV=0
\quad
{\rm when}
\quad
{x_1<0},
\qquad
RW=0
\quad
{\rm when}
\quad
{x_1>0}\,,
\end{equation}
\begin{equation}
\label{eq:inhomogeneoustrans}
(V,W)-g\ \in\  \caN
\quad
{\rm when}
\quad
x_1=0\,.
\end{equation}
The problem with inhomogeneous term $F$ can be reduced to this form by subtracting on the left a solution of the hyperbolic Cauchy problem $LU=F$ on $\RR^{1+d}$ with $U|_{t<0}=0$. Denote by  $ \widehat V(\tau,x_1,\eta), \widehat W(\tau,x_1,\eta), \hat g(\tau,\eta)$ the Fourier-Laplace transforms.  The transform $\widehat U$ is defined for $x_1\in \RR$, while  $\widehat V$ (resp. $\widehat W$)  is defined for $x_1<0$ (resp. $x_1>0$). The transforms  $\widehat V,\widehat W$ decay as $|x_1|\to \infty$.  $\widehat V, \widehat W$ satisfy the ordinary differential transmission problem
\begin{equation}
\label{eq:transformed}
L(\tau,d/dx_1,i\eta)\,\widehat V = 0
\ \ {\rm in}\ \
 x_1<0,
\qquad
R(\tau,d/dx_1,i\eta)\, \widehat W = 0
\ \ {\rm in}\ \
x_1>0,
\end{equation}
\begin{equation}
\label{eq:transformed2}
(\widehat V(\tau,0,\eta)\,,\,\widehat W(\tau,0,\eta))
\ -\
\hat g(s,\eta)
\ \in\  \caN
\,.
\end{equation}

Hersh's necessary and sufficient condition for well posedness of the
transmission problem is derived as follows.
Uniqueness of solutions of
 \eqref{eq:transformed},
\eqref{eq:transformed2}
for ${\rm Re}\, \tau>\tau_0,\, \eta\in \RR^{d-1}$
is equivalent to the fact that
there are no exponentially decaying solutions of the homogeneous
transmission problem.
That is,
\begin{equation}
\label{eq:disjoint}
\caN \ \cap\
\big(\dot E^-_L(\tau,\eta)\times \dot E^+_R(\tau,\eta)\big) \ =\
\{0\}.
\end{equation}
 In order to guarantee
existence, one imposes the maximality
condition,
\begin{equation}
\label{eq:directsum}
\caN \ \oplus
\big(\dot E^-_L(\tau,\eta)\times \dot E^+_R(\tau,\eta)\big) \ =\
\CC^N\times\CC^M
.
\end{equation}
Using \eqref{eq:dimensions}, this determines
the dimension of $\caN$ from the coefficients
$A_1$ and $\caA_1$ of
$L$ and $R$ respectively.

\begin{definition}
\label{def:hersh}   If the transmission problem \eqref{eq:LRhomogeneous},
\eqref{eq:inhomogeneoustrans}
satisfies \eqref{eq:directsum}
 for all ${\rm Re}\, \tau\,>\tau_0$ and $\eta\in \RR^{d-1}$
it is said to satisfy {\bf Hersh's condition}.
\end{definition}

\begin{theorem}  Hersh's condition is satisfied if and only if
there is an $r$ and a  $\lambda_0$ so that for all
$\lambda>\lambda_0$ and
$g$ supported in $t\ge 0$ with
$e^{-\lambda t} g\in H^{s+r}(\RR^{d}_{t,x^\prime})$
with values in $\CC^N\times \CC^M$
there is a unique $V,W$ supported in $t\ge 0$ with
$$
e^{-\lambda t} V\in H^s(]-\infty,\infty[\times\{x_1<0\})
\quad
{\rm and}
\quad
e^{-\lambda t} W\in H^s(]-\infty,\infty[\times\{x_1>0\})
$$
satisfying the transmission problem
\eqref{eq:LRhomogeneous},
\eqref{eq:inhomogeneoustrans}.
\end{theorem}

\noindent
{\bf Sketch of Proof.}
 We have shown that the Hersh condition
permits one to compute a candidate Fourier-Laplace
transform.
We outline how  the  condition
 implies the desired estimate.
 The method is to
 use the Seidenberg-Tarski
 Theorem \ref{thm:seidenber} to derive a lower bound on
the real parts
of the roots $\nu$  together with
a contour integral representation.
The same elements form the
heart of \cite{Kasahara:1970:WWP},
and
\S 12.9
of \cite{Hormander:2005:II}.
In the present context we treat
a transmission problem
rather than a boundary value problem.
In addition,
one needs to use the earlier lemmas to
treat the case when $x_1=0$ is characteristic.

Choose $\Lambda>\max\{\tau_0(L),
\tau_0(R)\}$.
The
equations
$$
\det  L(\tau, \nu, i\eta) \ =\ 0,
\qquad
\det  R(\tau, \nu, i\eta) \ =\ 0
$$
with $\Re \tau\ge  \Lambda$,
$\eta\in \RR^{d-1}$
have no purely imaginary roots.
Define
$$
\begin{aligned}
\zeta(R)
\ :=\
\min\Big\{|\Re \nu| \ :\
\eta\in \RR^{d-1}, \ \
& \Re\tau \ge\Lambda,\ \
|\tau|^2 + |\eta|^2
\le R^2,
\cr
&
\big\{\det L(\tau, \nu, i\eta) = 0\ \ {\rm or}\ \
\det R(\tau, \nu, i\eta) = 0\big\}
\Big\}.
\end{aligned}
$$
The Seidenberg-Tarski Theorem
\ref{thm:seidenber} implies that there is a $\rho\in {\mathbb Q}$ and $b\ne 0$ so that
$$
\zeta(R)
\ =\
R^\rho(b+o(R)),
\qquad
{\rm as}
\qquad
R\to \infty\,.
$$
Thus, there
are $C,N$ so that, for any  $\tau$, $\eta$ with $\Re \tau \ge \Lambda$,
\begin{equation}
|\Re \nu|\ \ge\
\frac{C}{1+|(\tau,\eta)|^N}.
\label{eqn:seiden}
\end{equation}

The solutions in $E^+_{R}(\tau,\eta)$ are written using a  contour integral representation of $\widehat W$ in the block form  of Lemma \ref{lem:blochform}. Here the matrix $H_{ij}$ depend of $(\tau,\eta)$. Denote by $D=D(\tau,\eta)$ the finite union of squares with centers at the roots with $\Re \nu<0$. The side of  each square is the smaller of $1$ and  half the distance of the root to  the imaginary axis.  Then
 \begin{equation}
 \label{eqn:contour}
\widehat W_I
 =
 \frac{1}{2\pi i}
 \oint_{\partial D}
e^{\tau x_1}\,
\Big(
\tau +
\big(H_{11} + H_{12}H_{22}^{-1}H_{21}\big)\Big)^{-1}
 d\tau\,\dot{\widehat W}_I,
\quad
\widehat W_{II}=
H_{22}^{-1}H_{21}\,\widehat W_I\,.
  \end{equation}
 The Seidenberg-Tarski Theorem \ref{thm:seidenber} applied
 to
 $$
 \max
 \Big\{
 |w|^2\ :\ |z|^2=1, \ \
 H_{22}w = z,
 \ \
  \Re \tau \ge \Lambda,\ \
  |\tau|^2 +|\eta|^2\le R^2\Big\}
  $$
  proves that
  $$
  \|H_{22}(\tau,\eta)^{-1} \|\ =\
  R^\beta(a+o(1)),
 \qquad
 a\ne 0,\quad \beta \in {\mathbb Q}\,.
 $$
This estimate together with \eqref{eqn:seiden}
yields with new $C,N$,
$$
\int_0^\infty
|\widehat W(\tau,x_1,\eta)|^2\ dx_1
\ \le \
C\,(1+|(\tau,\eta)|^{2N})| {\widehat  W}_I(0) |^2.
$$
With the analogous expression for
$V$ the solution of \eqref{eq:LRhomogeneous}
satisfies
$$
\int_{-\infty}^\infty
\Big|
\widehat{V}(\tau,x_1,\eta)
\Big|^2\ dx_1
\ \le \
C\,(1+|(\tau,\eta)|^{2N})|\widehat V_I(0)|^2.
$$

The Hersh condition asserts  that
for each $(\tau, \eta)$,  ${\widehat  W}_I(0) $
and  ${\widehat  V}_I(0)$ are uniquely determined
by $\hat g(\tau,\eta)$.  Seidenberg-Tarski Theorem \ref{thm:seidenber}
yields an estimate
$$
\big\|
{\widehat  W}_I(0) \,,\,
{\widehat  V}_I(0)
\big\|
\ \le\
C\, (1+|(\tau,\eta)|)^a
\,
\|\hat g(\tau,\eta)\|^2\,.
$$
The last three estimates together with Parseval's identity proves
the desired estimate,
$$
\exists\, C,N,\ \
\forall \,g,\ \
\forall \lambda>\Lambda,
\qquad
\big\|
e^{-\lambda t} \, U
\big\|_{L^2(\RR^{1+d})}^2
\ \le \
C \,
\sum_{|\alpha|\le N}
\big\|
e^{-\lambda t}\,
\partial_{t,x}^\alpha
g\big\|_{L^2(\RR^{1+(d-1)}
)}\,.
$$
This estimate proves the existence part of the Theorem.
\qed

\subsubsection{Necessary and sufficient condition for perfection}
\label{subsec:nasc}

The Fourier-Laplace method is  used to derive
a necessary and sufficient condition for perfection of
an absorbing layer.  Begin with a closer
analysis of the transform,
$\widehat U(\tau,x_1,\eta)$, of the solution of
the basic equation \eqref{hyp}.

When $A_1$ is invertible,  $\widehat U$ is
analysed as follows.
Denote by $\Pi_\pm(\tau,\eta)$ the  projectors
associated with the direct sum decomposition
$\dot E^+_L(\tau,\eta)\oplus \dot E^-_L(\tau,\eta)=\CC^N$.
Define $S_\pm(\tau,x_1,\eta)$ as the
Hom$(\CC^N)$ valued solutions of
$$
L\big(\tau\,,\, d/dx_1\,,\,i\eta\big)\,S_\pm\ =\ 0,
\qquad
S_\pm\big|_{x_1=0}\ =\
A_1^{-1}\,\Pi_\pm\,.
$$
Then $S_\pm$ decays exponentially as $x_1\to \pm \infty$
and
$$
\chi_{]-\infty ,0[}\,S_- \ +\ \chi_{[0,\infty[}\,S_+
$$
is the unique tempered fundamental solution of $L(\tau,d/dx_1,i\eta)$.
Decompose $\widehat F=\widehat F_-+\widehat F _+$,
$\widehat U=\widehat U_+ + \widehat U_-$
according to $E^+_L(\tau,\eta)\oplus E^-_L(\tau,\eta)=\CC^N$.
Then $\widehat U_-$ is the convolution of $\hat F_-$ with
$\chi_{]-\infty, 0[}\,S_-$
and
$\widehat U_+$ is the convolution
of $\widehat F_+$ with
$\chi_{[0,\infty[}\,S_+$.
In particular $\widehat U_-(\tau,0,\eta)$ vanishes on a neighborhood of
$[0,\infty[$ so $\widehat U(\tau,0,\eta)=\widehat U_+(\tau,0,\eta)\in \dot E^+_L(\tau,\eta)$.  The value
of $\widehat U$ in $x_1\ge 0$ satisfy the
homogeneous ordinary differential equation
 $L(\tau,d/dx_1,i\eta)\widehat U=0$.
 with initial value
$\widehat U(0)\in \dot E^+_L(\tau,\eta)$.

To reach the same  conclusion when $A_1$ is singular, apply the lemmas of the preceding section to the equation
$L(\tau,d/dx_1,i\eta)Z=0$.
 Lemma \ref{lem:transdata}
applied to $A=A_1$ and $M=\tau I+i\sum_{j=2}^d A_j\eta_j$ shows
 that both $\dot E^\pm_L(\tau,\eta)$ are subspaces of
$\GG$ and that
 the space
of solutions is a direct sum $E^+_L(\tau,\eta)\oplus E^-_L(\tau,\eta)$.
It follows that
$$
\dot E^-_L
(\tau,\eta)
\oplus
\dot E^+_L
(\tau,\eta)
\ =\ \GG
(\tau,\eta)\,.
$$
Repeating the analysis in the nonsingular
case applied to
\eqref{eq:Mtilde}
shows that $\widehat U(\tau,0,\eta)\in \dot E^+_L(\tau,\eta)$.

\begin{definition}
For a transmission problem $(L,R,\caN)$ satisfying Hersh's
condition (Definition \ref{def:hersh}),
${\rm Re}\,\tau>\tau_0$ and $\eta\in \RR^{d-1}$, define
the {\bf reflection operator},
$
H(\tau,\eta)
\, :\,
\dot E^+_L(\tau,\eta)\  \rightarrow\
\dot E^-_L(\tau,\eta)
$
as follows.
Hersh's condition implies that for
each $K\in \dot E^+_L(\tau,\eta)$ there is a unique $(\dot V,\dot W)\in
\dot E^-_L(\tau,\eta)\times \dot E^+_R(\tau,\eta)$ so that
$
(K,0) \equiv
(\dot V,\dot W)
\, {\rm mod}\,
\caN.
$
Define,
$
H(\tau,\eta)\,K\, :=\,
\dot V\,.
$
\end{definition}

\begin{theorem}
\label{thm:nascperfertion}
Suppose that the transmission problem $(L,R,\caN)$
satisfies the Hersh condition.   The following are equivalent.
\begin{enumerate}
[label={\it(\roman{*})}, ref={\it(\roman{*})},leftmargin=0.7cm]
\item The
transmission problem is perfectly matched in the sense of
Definition \ref{def:perfection}.
\item There is a $\tau_0\in \RR$ so
that for all  ${\rm Re}\, \tau>\tau_0$
and $\eta\in \RR^{d-1}$,
$H(\tau,\eta) \ =\ 0$.
\item  There is a $\tau_0\in \RR$ so
that for all  ${\rm Re}\, \tau>\tau_0$
and $\eta\in \RR^{d-1}$,
\begin{equation}
\label{eq:perfectcond}
\forall K_L\in \dot E^+_L(\tau,\eta),
\ \
\exists \,! \,K_R \in \dot E^+_R(\tau,\eta),
\quad
{\rm such\ that}
\quad
(K_L\,,\,K_R)\in \caN\,.
\end{equation}
\end{enumerate}
\end{theorem}
\begin{proof}  Conditions $( ii)$   and $( iii)$
are clearly equivalent.

For the equivalence with $( i)$,
compare the values of $\widehat U$ and $\widehat V$
in $\{x_1<0\}$.
Since both satisfy $LZ=F$ and decay as $x_1\to -\infty$
it follows that
$L(\widehat V -\widehat U)=0$,
 so,
 $\widehat V-\widehat U:=\Gamma$ is an $\dot E^-_L$ valued solution
of $L\Gamma=0$.   Since $F=0$ in $x_1>0$,
$\widehat W\in E^+_R$.
The transmission condition requires that
\begin{equation}
\label{eq:decomp}
\caN\ \ni\
(\widehat V(0)\,,\,\widehat W(0))
\ =\
(\widehat U(0) +\Gamma(0)\,,\,\widehat W(0))
\ =\
(\widehat U(0)\,,\,0)
\ +\
(\Gamma(0)\,,\,\widehat W(0))\,.
\end{equation}
Since $(\Gamma(0),\widehat W(0))\in \dot E^-_L(\tau,\eta)\times \dot E^+_R(\tau,\eta)$,
\eqref{eq:decomp} expresses $(\widehat U(0),0)$  as a sum of
an element in $\caN$ and an element of $\dot E^-_L(\tau,\eta)\times \dot E^+_R(\tau,\eta)$.
The Hersh condition
\eqref{eq:directsum} asserts that such a decomposition is unique.
Therefore $(\widehat V(0),\widehat W(0))$  is uniquely determined
from $\widehat U(0)$.

The method is perfectly matched if and only if for all $F$ supported in $x_1<0, t\ge 0$
\[
V\ =\ U\big|_{x_1<0}\,.
\]
This occurs if and only $\Gamma$ vanishes for $x_1<0$ which holds
if and only if $\Gamma(0)=0$.

If the method is perfectly matched, then in the decomposition \eqref{eq:decomp} one has $\Gamma(0)=0$.  Then $(\widehat U(0), \widehat W(0))\in \caN$. Lemma \ref{lem:seekF}  asserts that for any $K\in \dot E^+_L$ there is an $F$ so that  $\widehat U(0)=K$. This proves that \eqref{eq:perfectcond} holds.

Conversely if \eqref{eq:perfectcond} holds, then in the decomposition \eqref{eq:decomp}, $\Gamma(0)=0$ so $\Gamma=0$. It follows that $U|_{x_1<0}=V$.
\end{proof}

\begin{remark}
{\bf 1.}  When \eqref{eq:perfectcond} holds,
 the decomposition of $(K,0)\in \CC^N\times\CC^M$
in the direct sum
\eqref{eq:directsum} is,
$$
(K\,,\,0)\ =\
(K\,,\,W(K))
\ -\
(0\,,\,W(K))
\ \in\
\caN
\oplus
(E^-_L\times E^+_R)
\,.
$$

\noindent
{\bf 2.}  With $K=U(0)$ as above the solution $(V,W)$
of the ordinary differential equation
transmission problem is given by
$V=U|_{x<0}$ and $W$
 is the solution of
$RZ=0$ with $Z(0)=-W(K)$.

\noindent
{\bf 3.}   In the important case where $N=M$,
invertible $A_1$ and $\caA_1$
 and  transmission
condition  $\caN =\{V=W\}$, the perfection criterion $(iii)$  asserts that
$\dot E^+_L(\tau,\eta) = \dot E^+_R(\tau,\eta)$.
\end{remark}

We present a typical example showing that
the natural absorbing layers  are virtually never
perfectly matched in dimension $d\ge 2$.

\begin{prop}  Consider  the dissipative
symmetric hyperbolic example with
$d=N=M=2$,
$$
A_1=
\begin{pmatrix}
1 & 0\cr
0 & -1
\end{pmatrix}
,
\quad
A_2=
\begin{pmatrix}
0 & 1\cr
1 & 0
\end{pmatrix}
,
\quad
R=L+P,
\quad
P=P^*\ge 0,
\quad
\caN=\{(V,W):V=W\}\,.
$$

\noindent
{\it (i)}
The transmission problem is perfectly matched if and only
if $P=0$.

\noindent
{\it (ii)}  The corresponding problem with $d=1$ is perfectly
matched if and only if $P$ is diagonal.
\end{prop}

\begin{proof}  Define
\[
M_L:=A_1^{-1}[\tau + i\eta A_2]
\,,
\qquad
M_R:=A_1^{-1}[\tau + i\eta A_2]+A_1^{-1}P\,,
\]
so that $A_1^{-1}L(\tau, \partial_1,i\eta)=\partial_1 +M_L(\tau,\eta)$
and similarly for $A_1^{-1}R(\tau,\partial_1,i\eta)$.
For ${\rm Re}\,\tau>0$ and $\eta\in \RR$,
the matrices $M_L$ and $M_R$ have one eigenvalue with positive
real part and one with negative real part.  The eigenspace
corresponding to positive (resp. negative) real part eigenvectors
is   equal to
$\dot E^+_L(\tau,\eta)$ (resp. $\dot E^+_R(\tau,\eta)$).
Therefore the necessary and sufficient condition for
perfection is that for $\Re \tau > \tau_0$ and any $\eta$,
$
\dot E^+_L(\tau,i\eta)
\ =\
\dot E^+_R(\tau,i\eta)\,.
$

Since
\[
L(\tau, \rho, i\eta)
=
\begin{pmatrix}
\tau +\rho & i\eta\cr
i\eta &\tau-\rho
\end{pmatrix}
,
\qquad
{\rm and},
\qquad
\det L(\tau,\rho,i\eta) = \tau^2 -\rho^2+\eta^2,
\]
the eigenvalue of $M_L(\tau,\eta)$ with positive
real part is $\rho=\sqrt{\tau^2+\eta^2}$.  The eigenspace
is the kernel of $L(\tau,\rho,i\eta)$.  Therefore
\begin{equation}
\label{eq:Epluswave}
 \dot E^+_L(\tau,\eta)
 \ =\
 \CC(-i\eta,\tau + \rho)\,.
 \end{equation}

Since $M_R=M_L+A_1^{-1}P$, a necessary condition
is that the family of vectors $v(\eta,\tau):=(-i\eta,\tau + \rho)$  be  eigenvectors of the constant matrix $A_1^{-1}P$, which is possible only if $A_1^{-1}P$ is a constant multiple of the identity. Therefore $P = cA_1$.  Since $P\ge 0$
and $A_1$ has eigenvalues of both signs, it follows that
$c=0$ proving $(i)$.

 In the one dimensional case there is just one eigenvector $(0,1)$ which must be an eigenvector of $A_1^{-1}P$.  Since $(0,1)$ is also an eigenvector of $A_1$ it follows that  $(0,1)$ must be an eigenvector of $P$.    Since $P=P^*$, the orthogonal vector $(1,0)$ is also an  eigenvector and $P$ is diagonal.  Conversely, if $P$ is diagonal the condition is satisfied.

\end{proof}

\begin{remark}\newline
{\bf 1.}    Examples verifying perfection
for a  family of absorbing layers related to but
not including those of B\'erenger
are presented in \cite{Appelo:2006:PML}.
To our knowlege Hersh's criterion for B\'erenger's
layers has not been verified before.

\noindent
{\bf 2.}  The perfection criterion is related to the plane wave criterion
of B\'erenger.  We examine the relation in \S 4.1.6.
\end{remark}

\subsubsection{Hersh's condition for B\'erenger's PML with piecewise constant  $\sigma_1$}
\label{subsec:wellposed}

Of our earlier results, only those of Section \ref{subsec:elliptic} apply to
discontinuous
absorptions.  So, if the generator is not elliptic,
(for example. the PML Maxwell system of B\'erenger),
the preceding results do not prove that the initial value
problem is well posed.  In this section we   prove  that
the
doubled operators of B\'erenger define a (weakly) well
posed
initial value problem provided that
\begin{equation}
\label{eq:coeffpiece}
\sigma_j\equiv 0\ \  \mbox{ for }\ \  j \ge 2
\qquad
 \mbox{ and }
 \qquad
  \sigma_1(x_1)\equiv \sigma^\pm  \mbox{ in }\RR^d_\pm
\,,
\end{equation}
and, the  constant coefficient
operators $\tilde L$ on $\RR^d_\pm$ are both
(weakly) hyperbolic.

The unknown $\widetilde U$ satisfies \eqref{PMLgen}.
Denote by $\tiu^{\pm}=\{U^{\pm}_1,\dotsc,U^{\pm}_d \}$
the restriction of the unknown $\widetilde U$ to
$\RR^d_\pm$.   They satisfy
differential equations in the half spaces $\RR^d_\pm$.

\begin{lemma}
\label{lem:trace}
For $\widetilde U$ locally square integrable on a neighborhood
of $(\ut,\ux)\in \{x_1=0\}$ the following are equivalent.
\begin{enumerate}
[label={\it(\roman{*})}, ref={\it(\roman{*})},leftmargin=0.7cm]
\item $\widetilde L \widetilde U\in L^2$  on a neighborhood
of $(\ut,\ux)$ in $\RR^{1+d}$
 in the sense of distributions.
\item There is a neighborhood $\caO$ of $(\ut,\ux)$
so that $ \widetilde L\widetilde  U^\pm$
is square integrable on
$\caO \cap \RR^d_\pm$ and
$[ \widetilde A_1\widetilde U] =0.$
\end{enumerate}
\end{lemma}

\vskip.2cm
\begin{remark}\newline
{\bf 1.}   The first hypothesis is often verified
by combining
$\widetilde L\widetilde U + \widetilde  B(x)\widetilde  U\in L^2_{loc}$, $\widetilde U\in L^2_{loc}$
and $\widetilde B\in L^\infty_{loc}$.

\noindent
{\bf 2.}    $[\widetilde A_1\widetilde  U]$ makes sense since the differential equation
implies
$$
\partial_1\big(\widetilde A_1 \widetilde U^+\big) \ \in \
L^2_{loc}\big(]0,\eps[\ ;\ H^{-1}_{loc}(\RR^d_{t,x^\prime})\big).
$$
With $\widetilde U\in L^2(]0,\eps[\,;\, H^{-1}_{loc}(\RR^d)$ this
implies that $\widetilde A_1\widetilde U^+\in C([0,\eps[\,;\, H^{-1/2}_{loc}(\RR^d))$.
An analogous result holds for $\widetilde A_1\widetilde  U^-$.  Therefore the traces
from both sides and the jump  are well defined elements of
$H^{-1/2}_{loc}$.

\noindent
{\bf 3.}     is clear on a formal level
since if $\widetilde A_1 \widetilde U$
were discontinuous there would be a $\delta(x_1)$ term from the
differential operator $\widetilde L$ applied to $\widetilde U$.

\noindent
{\bf 4.}
The standard proof based on these remarks is omitted.
\end{remark}

We have supposed that the nonzero data are
 initial values $f^\pm(x)$.   By the usual subtraction
 one can convert the problem to one with homogeneous
 initial values and right hand side and inhomogeneous
 transmission condition.
 In this way, the determination of
 $\widetilde U^\pm$
 is reduced to finding
$ \widetilde  W^\pm$ satisfying the inhomogeneous
transmission problem
\begin{equation}
\label{eq:couplage}
\til_1(\partial_t,\partial_x) \widetilde W^{\pm} +
\widetilde B^{\pm}\widetilde   W^{\pm}\ =\ 0\,,
\qquad
   \widetilde A_1[\widetilde   W]\ =\ \widetilde g\,,
\end{equation}
where
\begin{equation}\label{eq:btilde}
\widetilde B^\pm\ :=\
        \begin{pmatrix}
            \sigma^{\pm} I_N  \quad                &   0       & \dots     & 0 \\
            0 \quad&   0      &\dots      & 0 \\
            \vdots                  &   \vdots  & \ddots    & \vdots \\
            0 \quad&   0       &\dots      & 0
        \end{pmatrix}\,,
     \end{equation}
and
 $\widetilde g(t,x)$  take values
in ${\rm Range}\,\widetilde A_1$.
The unknowns
$\widetilde W$ and source $\widetilde g$ are vectors of length $dN$.

\begin{theorem}
\label{thm:Hersh}
Suppose that $L(\partial)$ is a hyperbolic operator nondegenerate with respect
to $x_1$ and that
the B\'erenger's doubled operator $\widetilde L$ is
weakly hyperbolic for $\sigma=\sigma^\pm$.
Then, the transmission problem \eqref{eq:couplage}
with absorption \eqref{eq:coeffpiece}
satisfies Hersh's condition.
\end{theorem}

\begin{proof}  Drop the tildes
on the Fourier-Laplace transforms of $\widetilde W,\tilde{g}$ for ease of reading.
The transformed problem is
 \begin{equation}
\label{jbr2}
\Big(
{\widetilde A}_1
\frac{d}{dx_1}
\ +\
{\widetilde L}(\tau,0,i\eta)
\Big)
\widehat {W^\pm}\ =\
0\,,
\qquad
\widetilde A_1\,[\widehat{W}] \ =\ \hat g.
\end{equation}
The condition of Hersh
is that for an arbitrary right hand side $\hat g$ in ${\rm range}\,\widetilde A_1$
this  transmission problem has one and only one solution.

Denote by $E^\pm_{\widetilde L}(\tau,\eta,\sigma)$ the spaces
associated to the
B\'erenger operator operator $\widetilde L$
with absorption $\sigma$.
The uniqueness of solutions of \eqref{jbr2} is equivalent
to
\begin{equation}
\label{eqn:uniqueness}
\widetilde A_1\dot E^-_{\widetilde L}(\tau,\eta, \sigma^-)\ \cap \
\widetilde A_1\dot E^+_{\widetilde L}(\tau,\eta, \sigma^+)
\ =\
\{0\}.
\end{equation}
Existence is equivalent to
\begin{equation}
\label{eqn:Hershsurjective}
\widetilde A_1\dot E^-_{\widetilde L}(\tau,\eta, \sigma^-)\ +\
\widetilde A_1\dot E^+_{\widetilde L}(\tau,\eta, \sigma^+)
\ =\
{\rm range}\, \widetilde A_1.
\end{equation}
Part $(ii)$
 of Lemma \ref{lem:transdata}
 implies that
 $$
 {\rm dim}\,\big(
 \widetilde A_1\dot E^-_{\widetilde L}(\tau,\eta,\sigma^-)\big)
 \ +\
  {\rm dim}\,\big(
  \widetilde A_1\dot E^-_{\widetilde L}(\tau,\eta,\sigma^+)
  \big)
  \ =\
  {\rm dim}\,\big({\rm range}\,\widetilde A_1\big)
  \,,
  $$
  so
  \eqref{eqn:uniqueness}
  implies \eqref{eqn:Hershsurjective}.
     It remains to prove
    \eqref{eqn:uniqueness}.

For $\sigma>0$, the split  B\'erenger operator
$\widetilde L$ is hyperbolic so for $\Re \tau>\tau_0(\sigma)$
the solutions of
\begin{equation}
\label{eqn:Lwithsigma}
\widetilde L(\tau,d/dx_1,i\eta)\,\widehat U
\ =\
0
\end{equation}
are generated by exponentially growing and exponentially  decaying
solutions.  The next lemma identifies these
solutions in terms of the corresponding solutions
of
\begin{equation}
\label{eqn:sigmaequalzero}
L(\tau,d/dx_1,i\eta)\,\widehat V=0\,.
\end{equation}
The result shows that the
traces at $x_1=0$, $\dot E^\pm_{\widetilde L}$,  are independent of
$\sigma$.

\begin{lemma}
\label{lem:rescale}
For  $\sigma>0$,
$\Re \tau>\tau_0(\sigma)$, $\eta\in \RR^{d-1}$,
\begin{enumerate}
[label={\it(\roman{*})}, ref={\it(\roman{*})},leftmargin=0.7cm]
\item  The map
$$
\widehat V(x_1)
\quad
\mapsto
\quad
\bigg(
\frac{\eta_1}{\tau}
\,
A_1\widehat V\Big((\tau+\sigma)x_1/\tau\Big)
\,,\,
\frac{\eta_2}{\tau}\,
A_2\widehat V\Big((\tau+\sigma)x_1/\tau\Big)
\,,\,
\dots\,,\,
\frac{\eta_2}{\tau}\,
A_d\widehat V\Big((\tau+\sigma)x_1/\tau\Big)
\bigg)
$$
is an isomorphism  from solutions of
\eqref{eqn:sigmaequalzero}
onto the solutions of \eqref{eqn:Lwithsigma}.
\item  $\mu$ is a root of
$\det L(\tau,\cdot,i\eta)=0$ if and only if
$\nu=(\tau + \sigma)\mu/\tau$ is a root
of $\det \widetilde L(\tau,\cdot,\eta)=0$.
\item  For the roots in $( ii)$,
the real parts of $\mu$ and $\nu$ have
the same sign.  In particular, the map
in $(i)$ is an isomorphism
$
E_L^\pm(\tau,\eta)
\
\mapsto
\
E_{\widetilde L}^\pm(\tau,\eta,\sigma).
$
\item  The map $\widehat W=(\widehat W_1,\dots, \widehat W_d)\ \mapsto\
\sum_j \widehat W_j$ is an isomorphism
$E^\pm_{\widetilde L}(\tau, \eta, \sigma) \to E^\pm_L(\tau,\eta)$.
\end{enumerate}
\end{lemma}

\begin{remark}  In $(i)$  it important to know that the solutions
$\widehat V(x_1)$ are entire analytic functions of  $x_1$ so it makes
sense to evaluate $\widehat V$ at points off the $x_1$-axis.  In
 the literature this is sometimes called a {\it complex change
of variables}.
It is only reasonable for analytic solutions.  A related
idea is used in the
Fourier-Laplace analysis for general $\sigma_1(x_1)$
presented in \S 4.2.
\end{remark}

\noindent
{\bf Proof of Lemma \ref{lem:rescale}}
\newline\noindent $( i)$    If $\widehat  U=(\widehat U_1,\dots,\widehat U_d)$ satisfies
\eqref{eqn:Lwithsigma} then with
$\widehat W := \sum_j \widehat U_j$,
\begin{equation}
\label{eqn:tildetransformed}
A_1 \frac{d \widehat W}{dx_1} \ +\
(\tau + \sigma)\,\widehat U_1\ =\ 0,
\qquad
\tau \widehat U_j \ +\ i\eta_j A_j\,\widehat W \ =\ 0, \quad j=2,\dots, d\,.
\end{equation}
Multiply the first by $\tau$ and the last $d-1$ by
$(\tau +\sigma)$.  Sum and then divide by $\tau$
to find,
\begin{equation}
\label{eqn:rescaled}
A_1
\frac{d \widehat W}{dx_1}
\ +\
\frac{\tau +\sigma}{\tau}\,
L(\tau,0,i\eta)\,\widehat W^\pm
 \ =\ 0\,.
 \end{equation}
 Conversely if $\widehat W$ satisfies \eqref{eqn:rescaled}
 and $\widehat U_j$ for $j\ge 2$
  is defined  from $\widehat W$
 using
the last equations in \eqref{eqn:tildetransformed}
and  $\widehat U_1:=\widehat W-\sum_{j\ge 2}\widehat U_j$ then
$U$ satisfies \eqref{eqn:Lwithsigma}.

 The solutions $\widehat W$ to \eqref{eqn:rescaled}
 are exactly the
 $\widehat V((\tau+\sigma)x_1/\tau)$ with
 $\widehat V$ satisfying \eqref{eqn:sigmaequalzero}.
 This proves that the mapping in $(i)$
 is surjective.

 The set of solutions $\widehat V$ of \eqref{eqn:sigmaequalzero} has dimension
 ${\rk}\,A_1$.  The set of solutions
 of \eqref{eqn:Lwithsigma} has
 dimension ${\rk}\,\widetilde A_1=
 {\rk}\,A_1$
 (see \eqref{eq:Ltildecoeff}), so surjectivity implies
 injectivity.\\

\noindent  $( ii)$ {\rm and} $( iv)$   follow from $( i)$.\\

\noindent  $(iii)$
 Denote
 by  $K$ the mapping
 from $(i)$.
 When $\tau>0$, $(\tau +\sigma)/\tau$ is
 also positive and
 real.
 Therefore $K$ maps decaying
 (resp. increasing) solutions
 to decaying (resp. increasing) solutions.
    Thus for $\tau>\tau_0$ and real,
 \begin{equation}
 \label{eqn:Kimage}
 K(E^+_L(\tau,\eta))=E^+_{\widetilde L}(\tau,\eta,\sigma)\,.
 \end{equation}
For all $\Re\,\tau>\tau_0$, $K(E^+_L(\tau,\eta))$ is a subspace
of solutions of \eqref{eqn:Lwithsigma}  with
dimension equal to $\dim E^+_L(\tau,\eta)$.
 If \eqref{eqn:Kimage} were violated,
$K(E^+(\tau,\eta))$ would contain
 exponentially growing solutions.
 If this happened at $\utau,\ueta$ with $\Re\,\utau>\tau_0$
 consider the
 values $\tau(r) = \Re\,\utau + r\Im\,\utau$ for $0\le r\le 1$.
 For $r=0$ \eqref{eqn:Kimage} is satisfied while
 for $r=1$ it is violated.
 Let
 $$
 f(r)
 \ :=\
 \max\Big\{\Re
 \frac{\tau(r)+\sigma}{\tau(r)}\,\mu
 \ :\
 \det L(\tau(r)\,,\,\mu\,,\,i\ueta)=0\,,
 \ \ {\Re}\,\mu<0
 \Big\}\,.
 $$
 Then $f(0)<0$, $f(1)>0$ and $f$ is continuous,
 so there is a $0<\underline r<1$ so that $f(\underline r)=0$.
 Then for $\tau = \Re\utau + \underline r\,\Im\utau\,$
 there is a purely imaginary root.
 This violates the hyperbolicity of $\widetilde L$
 establishing
 \eqref{eqn:Kimage}
 This proves $( iii)$ completing the proof of
 the Lemma.
 \qed

\vskip.2cm

We now finish the proof of Theorem \ref{thm:Hersh}
by proving \eqref{eqn:uniqueness}.
Lemma \ref{lem:rescale}  implies that the spaces of
Cauchy data $\dot E^\pm_{\widetilde L}$ are independent
of $\sigma$.  Therefore if \eqref{eqn:uniqueness} is violated
then also
\begin{equation*}
\widetilde A_1\dot E^-_{\widetilde L}(\tau,\eta, \sigma^+)\ \cap \
\widetilde A_1\dot E^+_{\widetilde L}(\tau,\eta, \sigma^+)
\ \ne\
\{0\}.
\end{equation*}
This contradicts part $(iv)$  of Lemma \ref{lem:transdata}
for the operator $\widetilde L$ with absorption $\sigma^+$.
The proof of Hersh's condition is complete.
\end{proof}

In these problems with only one nonzero absorption
coefficient $\sigma_1$ and $\sigma_1=0$ when
$x_1<0$ one can consider a transmission problem
which is only split in $x_1>0$.  The next result shows
that this partially split problem satisfies  Hersh's condition
 if and only
if the fully split problem does.

Introduce the partially split problem $(L,R,\caN)$ where
\begin{equation}
\label{eq:halfsplit}
L=L_1(\partial),
\quad
R = \widetilde L_1+\widetilde B^+,
\quad{\rm with}
\quad
\sigma^+>0,
\quad
\caN \ :=\
\{(V,W)
\, :\,
V-\sum_jW_j\ \in \ {\rm ker}\,A_1
\}\,,
\end{equation}
with $\widetilde B^+$ given by \eqref{eq:btilde} and the split variable on the right is $\widetilde W=(W_1,\dots,W_d)$.

\begin{corollary}  Suppose that $\sigma_j=0$ for $j\ge 2$,
  and $\sigma_1^+ >0$.  Then the partially split
B\'erenger transmission problem $(L_1, \widetilde L_1 +
\widetilde B^+, \caN)$
defined by \eqref{eq:halfsplit} satisfies Hersh's condition if and
only if the fully split problem does.
\end{corollary}

\begin{proof}  Denote by $(V, \widetilde W)= (V, W_1,\dots, W_d)$ the variables
for the partially split problem and
$\big(\widetilde U,\widetilde W\big) =
\big((U_1,\dots,U_d),(W_1,\dots,W_d)\big)$
the split variables.   If $\widehat U(\tau,x_1,\eta),\widehat W(\tau,x_1,\eta)$ is an exponentially
decaying solution of the split Laplace-Fourier transformed homogeneous transmission problem,
then $(\widehat V,\widehat W)=(\sum_j\widehat U_j,\widehat W)$ is an exponentially decaying
solution of the partially split homogeneous transmission problem.

Conversely, if $\widehat V(\tau,x_1,\eta),\widehat W(\tau,x_1,\eta)$  is a  solution of the
homogeneous partially split problem, the computation leading to
\eqref{eqn:tildetransformed} shows that
$$
\widehat U_1 := -\tau^{-1}A_1\partial_1\widehat V,
\qquad
\widehat U_j := -\tau^{-1}\, i\eta_jA_j\,\widehat V, \quad j\ge 2,
$$
is an exponentially decaying solution of the
fully split
homogeneous
transmission problem.

Therefore, if either problem  has decaying solutions for $\eta$
real and ${\rm Re}\,\tau$ arbitrarily large, then so does the other.
\end{proof}

\subsubsection{Perfection for  B\'erenger's PML  with piecewise constant
$\sigma_1 $}
\label{subsec:perfectionberenger}

\begin{theorem}
\label{thm:berengerperfect}
With the hypotheses of Theorem \ref{thm:Hersh}, the B\'erenger transmission problem is perfectly matched. The B\'erenger transmission problem that is only split on the right is also perfectly matched.
\end{theorem}
\begin{proof}  Verify condition $(ii)$  of Theorem
\ref{thm:nascperfertion}.
For $K\in \dot E^+_{\widetilde L}(\tau,\eta,\sigma^+)$
consider the unique
decomposition guaranteed by the Hersh's
condition,
\begin{equation}
\label{eq:Kdecomp}
(K,0) \ =\ (W^-,W^+) + (F^-,F^+)\,,
\end{equation}
 where,
$$
\big(W^-,W^+\big) = \Big(
(W^-_1,\dots,W^-_d)
,
(W^+_1,\dots,W^+_d)\Big)\in \caN,
\ \ \
(F^-,F^+)\in
\dot E^-_{\widetilde L}(\tau,\eta,\sigma^-)
\times
\dot E^+_{\widetilde L}(\tau,\eta,\sigma^+)
.
$$
Perfection is equivalent to $F^-=0$.

By inspection, one such decomposition
\eqref{eq:Kdecomp}  is given by
$$
(K,0) \ =\ (K,K) \ +\
(0, -K)
$$
where we use the fact from Lemma \ref{lem:rescale}  that
$$
\dot E^+_{\widetilde L}(\tau,\eta,\sigma^-)
\ =\
\dot E^+_{\widetilde L}(\tau,\eta,\sigma^+)\,.
$$
As this decomposition satisfies $F^-=0$ the proof of the
first assertion is complete.

For the partially split case,
$K\in \dot E^+_L(\tau,\eta)$  has a unique
decomposition from the Hersh's condition,
$$
(K,0) \ =\ (W^-,W^+)\ +\ (F^-,F^+),
$$
with
$$
(W^-,W^+) = \big(
W^-
\,,\,
(W^+_1,\dots,W^+_d)\big)\in \caN,
\quad
(F^-,F^+)\in
\dot E^-_{L}(\tau,\eta)
\times
\dot E^+_{\widetilde L}(\tau,\eta,\sigma^+)
.
$$

Define
$$
W^+_j \ :=\
\frac{ \eta_j}{\tau} A_jK\,.
$$
Part $(i)$ of Lemma \ref{lem:rescale}
implies that $W^+\in \dot E^+_{\widetilde L}(\tau,\eta,\sigma^+)$.
In addition,
$\sum_j W^+_j=K$ so
$(K, W^+)\in \caN.$

By inspection
$$
\big(K\,,\,0\big)\ =\
\big(
K\,,\,
W^+\big)
\ +\
\big(
0\,,\, -W^+\big)
 $$
is the unique Hersh decomposition.  Since
$F^-$ vanishes for this one the proof is complete.
\end{proof}

\subsubsection{Analytic continuation for  Maxwell like systems and
B\'erenger's plane waves}
\label{sec:continuation}

In this section we investigate B\'erenger's method for operators,
including the Maxwell system, whose
characteristic polynomial is $\tau^p(\tau^2-|\xi|^2)^q$.
For ease of exposition we treat the case $d=2$ and the
explicit operator,
\begin{equation}
\label{eq:2dsystem}
L\ =\
\partial_t
 \ +\
 \begin{pmatrix}
 1 & 0
 \cr
 0 & -1
 \end{pmatrix}
\partial_1
 \ +\
  \begin{pmatrix}
 0 & 1
 \cr
 1 & 0
 \end{pmatrix}
\partial_2 \,.
\end{equation}
Analogous results are valid for the Maxwell system with
only slightly more complicated formulas.

For ${\re}\,\tau>0$ and
$\eta\in\RR$ there is exactly one root of
$\det L_1(\tau,\rho,i\eta)=0$ with ${\rm Re}\,\rho>0$ given by
\begin{equation}
\label{eq:rootrho}
\rho \ =\ \sqrt{\tau^2+\eta^2},
\qquad
{\rm Re}\,\rho >0\,.
\end{equation}
The corresponding eigenspace
$\dot E^+_L(\tau,\eta)$  from \eqref{eq:Epluswave}
is spanned by $\Phi(\tau,\eta) = (-i\eta,\tau +\rho)$.
If $\widetilde L_1$ is the B\'erenger  operator doubled
in the $x_1$ direction one has the same roots
and $\dot E^+_{\widetilde L_1}$ is spanned by
$\big(\rho A_1  \Phi, i\eta A_2\Phi\big)$.

\begin{prop}
\label{prop:analyticity}
\begin{enumerate}
[label={\it(\roman{*})}, ref={\it(\roman{*})},leftmargin=0.7cm]
\item For each $\eta$,
$\rho(\tau ,\eta)$, $\dot E_L^-(\tau,\eta)$, and
 $\dot E_L^+(\tau,\eta)$ are
holomorphic in $\Re \tau>0$ with continuous
extension to $\Re \tau\ge 0$.
\item If $\sigma_1>0$ then for ${\rm Re}\,\tau>0$
and $\eta\in \RR$
the equation
$\det \widetilde L(\tau,\nu,i\eta) =0$
has exactly one root $\nu$ with positive real part.
It is given by $\nu=(\tau +\sigma_1)\rho/\tau$.
\item
For $\sigma_1\ge 0$, the relation \eqref{eq:directsum} with
$\widetilde L_1$ on the left and $\widetilde L$
on the right
is satisfied  on $\big\{\Re \tau\ge 0\,,\, \rho\ne 0\big\}$.
\item The mapping $H(\tau, \eta)$ is for each $\eta$
holomorphic in ${\rm Re}\, \tau >0$ with continuous
extension to $\big\{\Re \tau\ge 0\,,\, \rho\ne 0\big\}$.
\end{enumerate}
\end{prop}

 \begin{proof}
$(i)$  For ${\rm Re}\, \tau >0$
there are two roots $\pm\rho$ with $\rho$ from \eqref{eq:rootrho}.
One has strictly positive real part and the other strictly negative.
Each is holomorphic in ${\rm Re}\,\tau>0$.
Holomorphy for
$\Phi(\tau,\eta)$ follows from
its expression in
terms of $\rho$.  As $\Phi$ is a basis for $\dot E_{L}^+(\tau,\eta)$
holomorphy of the latter follows.

So long as the eigenvalues $\pm \rho$
remain apart as ${\rm Re}\,\tau\to 0$ they and their eigenspaces
will be holomorphic.    The delicate case is when $\tau^2 + \eta^2 \to 0$.
The limiting points are $(\pm i\eta, \eta)$.

If $\eta=0$, then $\rho=\tau$ and the eigenspace is $(0,1)$.
Both are continuous up to the boundary.

When $\eta\ne 0$ one has $\rho\to 0$ so $\rho$ is continuous up
to the boundary.
Then $\Phi$ is   continuous
up to the boundary and nonvanishing from its expression
in terms of $\rho$.
Therefore $ \dot E_{L}^+$ is continuous up to the boundary.

\vskip.2cm

\noindent $(ii)$ It suffices to remark that this is
an eigenvalue and then to show that the real
part is positive.  For the latter compute
$$
\frac{\partial}{\partial \sigma}
\bigg(
\frac{\tau +\sigma}{\tau}
\sqrt{\tau^2+\eta^2}
\bigg)
\ =\
\frac{
\sqrt{\tau^2+\eta^2}
}{\tau}
\ =\
\sqrt{
1 + \eta^2/\tau^2}\,.
$$
For ${\rm Re}\,\tau>0$ this has positive real part
so the real part of the eigenvalue is increasing
as a function of $\sigma$ so is positive for all
$\sigma\ge 0$.

\vskip.2cm

\noindent $(iii)$   It suffices to show that
\eqref{eq:disjoint} is valid for such $s,\eta$.
Suppose $(v,w)=(v_1,v_2,w_1,w_2)\in
\dot E_{ {\tilde L}_1}^-\times\dot E_{\tilde L}^+$.
Must show that $v_1+v_2\ne w_1+w_2$.  Since
$(w_1,w_2)\in \dot E_{\tilde L}^+$ it follows that
$w_1+w_2\in \dot E_L^+$.  Similarly
$v_1+v_2\in \dot E_{ { L}_1}^-$.
Thus it suffices to show
that $\dot E_{L_1}^-$ and $\dot E_L^+$ are uniformly transverse
as $\Re \tau\to 0$.
It suffices to show that
$(i\eta, \tau +\rho)$ and
$(i\eta, \tau -\rho)$ are uniformly independent.
This follows from $\rho\ne 0$.

\vskip.2cm

\noindent $(iv)$
The holomorphy of $H$ follows from $(i)$.
The continuous extension follows from
$(i)$ and $(iii)$ .
\end{proof}

Since the method is perfectly matched, $H=0$ for
$\Re \tau>0$.
By continuity
the map vanishes for purely imaginary
$\tau\ne 0$.  This shows that
for $\{{\rm Re}\, \tau\ge 0\}\setminus 0$,
the function equal to
$$
e^{i\tau t + \rho(\tau,\eta)  x_1 + i\eta x_2} \widetilde \Phi
\ \
{\rm for}\ \
 x_1<0,
 \quad
 {\rm  and}
 \quad
  e^{i\tau t + \rho(\tau,\eta)  x_1 + i\eta x_2}\
  e^{-\sigma \rho x_1/\tau}\  \widetilde \Phi
\ \ {\rm for}\ \
x_1>0,
$$
satisfies the B\'erenger transmission problem.
For ${\rm Re}\,\tau>0$ these solutions decay (resp. grow)
exponentially as $x_1\to \infty$ (resp. $x_1\to -\infty$).
Though such solutions serve to verify perfection they
don't look very physical in isolation.

On the other hand, when $\tau$ is purely imaginary
and not equal to zero, the solution is a bounded plane
wave in $x_1<0$ and is a plane wave modulated
by an exponentially  decaying factor in $x_1>0$.
These are the solutions which  B\'erenger
constructed to show
that the method was perfectly matched.

In the language of the analytic objects
constructed in the preceding lemma, B\'erenger's
plane wave solutions show
that $H(is, \eta)=0$ when $s$ is real
valued with $s^2>\eta^2$.
For $\eta$ fixed  the function $\tau\mapsto H(\tau,\eta)$
is holomorphic in the right half plane continuous up to the imaginary
axis punctured at $\pm i | \eta | $, and vanishes on the boundary interval
$\tau=is\in i\RR$ with $s^2>\eta^2$.  By Schwarz reflection and
analytic continuation this implies that $H$ vanishes in the right
half plane.

In summary, \textit{the computation of B\'erenger is
actually  sufficient to prove perfection for Maxwell's system given the structures
provided in this paper}.

\vskip.2cm

\begin{remark}  The perfection argument based on  plane
waves is not valid in full generality where the objects like
$\dot E$ and $H$ are analytic in ${\rm Re}\,\tau>\tau_0$
with $\tau_0>0$.   This is the case, for example, whenever
the absorbing layer is amplifying.
\end{remark}

\subsection{Fourier-Laplace analysis with variable $\sigma_1(x_1)$}
\label{subsec:FLAV}

Consider the case of only one nonzero
$\sigma_1(x_1)$.  If $\widetilde L$ is
hyperbolic for one constant value $\underline\sigma_1\ne 0$
the  scaling $(t,x)\mapsto (a t,ax)$ shows that
$\widetilde L$ is hyperbolic for $\underline\sigma_1/a$.
Therefore $\widetilde L$ is hyperbolic for all constant values $\sigma_1$.

The results of \S
 \ref{subsec:piecewiseconstant}
will be
extended to the case $\sigma_j=0$ for $j\ge 2$
and variable coefficient $\sigma_1(x_1)$.
The
Fourier-Laplace
transform $\widehat U(\tau, x_1,\eta)$
of the B\'erenger split operator
satisfies
$$
\widetilde L(\tau,d/dx_1,\eta)\widehat U\ =\ \widehat F,
\qquad
-\infty<x_1<\infty\,,
$$
with variable
coefficient $\sigma_1(x_1)$.

The first line of the proof of Lemma
\ref{lem:rescale} yields \eqref{eqn:tildetransformed}
with $\sigma=\sigma_1(x_1)$.
As in the proof of that lemma one derives
\eqref{eqn:rescaled} now with
$\sigma=\sigma_1(x_1)$.  The important
observation is that the $x_1$ dependence
of the coefficient appears only as a scalar
prefactor in \eqref{eqn:rescaled}.  Such
equations will be analysed in the same
way as the equations in Lemma \ref{lem:blochform}.

\subsubsection{Well posedness by Fourier-Laplace with variable $\sigma_1(x_1)$}
\label{sec:FLsigmavariable}

\begin{theorem}
\label{thm:variablesigma}
Suppose
$\sigma_j=0$ for $j\ge 2$ and
$\sigma_1(x_1) \in L^\infty(\RR)$
is real valued.  Suppose in addition that $L$ is nondegenerate with restect to $x_1$, and
for one value $\sigma_1\ne 0$, $\widetilde L$
is hyperbolic.
Then there is a $\tau_0>0$ and $m$ so that for
all $\lambda>\tau_0$ and
$F\in e^{\lambda t}L^2\big(\RR_t\,:\,H^m(\RR^{d}_{t,x^\prime})\big)$
there is a unique solution  solution $\widetilde U\in e^{\lambda t}L^2(\RR^{d+1})$
to the   B\'erenger split problem $\widetilde L \widetilde U=F$.  In addition, there is a constant
$C$ independent of $F,\lambda$ so that
\begin{equation}
\label{eq:expestimate}
\big\|
e^{-\lambda t}\widetilde U\big\|_{L^2(\RR^{1+d})}
\ \le \
C\,
\big\|
e^{-\lambda t} F
\big\|_{
L^2\big(\RR_t\,:\,H^m(\RR^{d}_{t,x^\prime})\big)
}\,.
\end{equation}

\end{theorem}

\begin{remark}
\newline
{\bf 1.}   The condition
$\widetilde U\in e^{\lambda t}L^2$ implies
that $\widetilde U$ tends to zero at $t\to -\infty$ as does $F$.

\noindent
{\bf 2.}     If $F$ is supported in $t\ge t_0$ it
follows from \eqref{eq:expestimate} on sending $\lambda\to \infty$
that $\widetilde U$ is supported in $t\ge t_0$.
\end{remark}

\noindent
\begin{proof}
The values of the Fourier Laplace Transform of $W=\sum U_j$ are computed from the
ordinary
differential equation
\begin{equation}
\label{eq:inhomog}
 A_1
\frac{d \widehat W}{dx_1}
\ +\
\frac{\tau +\sigma_1(x_1)}{\tau}\,
 L_1(\tau,0,i\eta)\,
\, \widehat W
\ =\
\widehat F\,.
\end{equation}

As in Lemmas \ref{lem:blochform} and \ref{lem:L_1},
 transform to the equivalent form,
$$
\begin{pmatrix}
I & 0\\
0 & 0
\end{pmatrix}
\frac{d\widehat W}{dx_1}
\ +\
\frac{\tau +\sigma_1(x_1)}{\tau}
\begin{pmatrix}
H_{11} & H_{12}
\cr
H_{21} & H_{22}
\end{pmatrix}
 \widehat W\ =\
\widehat F,
\qquad
H_{22}\ {\rm invertible}.
$$
Denote the decomposition as $W=(W_I,W_{II})$ and similarly
$F$.  The invertibility of $H_{22}$
from Lemma \ref{lem:L_1}
yields,
\begin{equation}
\label{eq:uII}
\widehat W_{II} \ =\
H_{22}^{-1}\big( \widehat F_{II}
-H_{21}\widehat W_I
\big)\,.
\end{equation}
It suffices to find $\widehat W_I$ which is determined from,
$$
\frac{d\widehat W_I}{dx_1} \ +\
\frac{\tau +\sigma_1(x_1)}{\tau}\,
M(\tau,\eta)\,
\widehat W_I
\ =\ \widehat G\,,
\qquad
M(\tau,\eta)
:=
H_{11} - H_{12}H_{22}^{-1}H_{21},
\ \
\widehat G :=
\widehat F_I + H_{22}^{-1}\widehat F_{II}\,.
$$

The hyperbolicity of $\widetilde L$ implies that
$M$ has no purely imaginary eigenvalues.
Correspondingly there is the decomposition,
into the spectral parts with positive and negative
imaginary parts,
$$
W_I\ =\ W_I^+ + W_I^-,
\qquad
G\ =\ G^+ + G^-,
\qquad
M \ =\
M^+\ \oplus \ M^-\,.
$$

For $\sigma$ constant, part $(iii)$ of
Lemma \ref{lem:rescale} (using the hyperbolicity of $\widetilde L$)
implies that
for ${\rm Re}\,\tau$ sufficiently large
(depending on $\sigma)$, one has
the spectral decomposition,
$$
\frac{\tau +\sigma}{\tau}\,
M(\tau,\eta)
\ =\
\frac{\tau +\sigma}{\tau}\,
M(\tau,\eta)^+
\ \oplus\
\frac{\tau +\sigma}{\tau}\,
M(\tau,\eta)^-
$$
corresponding to spectra with positive and negative
real parts.

\begin{lemma}
If
$g(x_1)
$
satisfies
$
dg(x_1)/dx_1\ =\
\sigma_1(x_1)
\,.
$
Then
$$
\Big(
\frac{d}{dx_1} +
M
\Big)
\Big(
e^{g(x_1)M/\tau}\widehat U_I
\Big)
\ =\
e^{g(x_1)M/\tau}
\Big(
\frac{d}{dx_1} +
\frac{\tau +\sigma_1(x_1)}{\tau}
M
\Big)
\widehat U_I\,.
$$
\end{lemma}

\noindent
{\bf Proof of  Lemma.}  Since
$
\big(
d\,e^{gM/\tau}
\widehat U_I
\big)/dx_1
=
e^{gM/\tau}
\big(
g^\prime M\widehat U_I/\tau + d \widehat U_I/dx_1
\big)
$
one has
$$
\Big(
\frac{d}{dx_1} +
M
\Big)
\Big(
e^{g(x_1)M/\tau}\widehat U_I
\Big)
=
e^{gM/\tau}
\Big(
\frac{d}{dx_1}\widehat U_I
+
\big(
\frac{dg/dx_1\, M}{\tau} +
\frac{\tau M}{\tau}
\big)
\widehat U_I
\Big)\,,
$$
proving the desired identity.
\qed
\vskip.2cm

Therefore
$$
\widehat W_I
\ =\
e^{-g(x_1)M/\tau}
\Big(
\frac{d}{dx_1} +
M
\Big)^{-1}
\Big(
e^{g(x_1)M/\tau} \widehat G
\Big)\,.
$$
The unique $L^1$ fundamental solution of $\partial_1 +
M$
is equal to,
$$
e^{-x_1M^+}\,
\chi_{[0,\infty[}(x_1)
\ +\
e^{-x_1M^-}\,
\chi_{]-\infty,0]}(x_1)\,.
$$
Therefore,
$$
e^{gM/\tau}\,
\widehat W_I^+
\ =\
\big(
e^{-x_1M^+}\
\chi_{[0,\infty[}(x_1)
\big)*(e^{gM/\tau}\widehat G\big)^+
,
$$
$$
e^{gM/\tau}\,
\widehat W_I^-
\ =\
\big(
e^{-x_1M^-}\
\chi_{]-\infty,0]}(x_1)
\big)*\big(e^{gM/\tau}\widehat G\big)^-\,.
$$
The kernel of the integral operator mapping $\widehat G^+$ to
$\widehat W_I^+$ is equal to,
\begin{equation}
\label{eq:integralkernel}
\exp
\Big(
-(x_1-y_1)
\Big[
\frac{\tau + (g(x_1)-g(y_1)/(x_1-y_1)}{\tau}M(\tau,\eta)^+
\Big]
\Big)
\
\chi_{x_1\ge y_1}\,.
\end{equation}

\begin{lemma}
$$
\exists\, \tau_0=\tau_0(\mu),\ \
\forall\,
{\rm Re}\,\tau\ge \tau_0,\ \
\forall\, \eta\in \RR^{d-1},\ \
\forall\,
\sigma\in [-\mu,\mu],
\quad
{\rm spec}\
\frac{\tau +\sigma}{\tau} M^+(\tau,\eta)
\ \subset\
\{
{\rm Re}\,z>0
\}
\,.
$$
\end{lemma}

\noindent
{\bf Proof of Lemma.}  Part $(iii)$  of Lemma
\ref{lem:rescale}  allows one
to choose $\tau_1$ so that
for $\sigma=\mu$ one has the desired conclusion
for ${\rm Re}\,\tau>\tau_1$.
Then for $\lambda\in {\rm spec}\, M(\tau,\eta)^+$
one has
$$
{\rm Re}\,\lambda > 0,
\qquad
{\rm Re}\, \big(1+ \frac{\mu}{\tau}\big)\lambda
=
{\rm Re}\,
\frac{\tau +\mu}{\tau}
>0\,.
$$

For $0\le \sigma \le \mu$ write $\sigma = a + b\mu$
with nonnegative $a,b$ summing to 1.  It follows
that ${\rm Re}\, (1+ \mu/\tau)\lambda>0$.  This
proves that $\tau_1$ suffices to treat the nonnegative
values $0\le \sigma\le \mu$.

Choosing $\tau_2$ for $\sigma=-\mu$, that value
suffices for $-\mu\le \sigma\le 0$.  Set
$\tau_0$ equal to the maximum of $\tau_1$ and
$\tau_2$.
\qed

\vskip.2cm

The Seidenberg-Tarski Theorem \ref{thm:seidenber} shows that
the absolute values of the real parts of
the eigenvalues of
$M(\tau,\eta)$ are bounded
below by
$C(|\tau| + |\eta|)^{-N}$ for some $N$.
And also that  the spectral decomposition $V\mapsto (V^+,V^-)$
and its inverse
are both bounded polynomially in $|\tau,\eta|$.
More generally
for $\tau, \eta, \mu, \sigma$ as above,
$$
{\rm spec}\
\frac{\tau +\sigma}{\tau} M^+(\tau,\eta)
\ \subset\
\Big\{
\big|
{\rm Re}\, z
\big|
 > C(|\tau|+|\eta|)^{-N}
\Big\}
\,.
$$
Taking $\mu:=\|f\|_{L^\infty}$ one finds that for all $x_1,y_1$, the matrix
\begin{equation}
\label{eq:fundsolmatrix}
\frac{\tau + (g(x_1)-g(y_1)/(x_1-y_1)}{\tau}\ M(\tau,\eta)^+
\end{equation}
has spectrum in
$$
\Big\{
{\rm Re}\, z
 > C(|\tau|+|\eta|)^{-N}
\Big\}\,,
\qquad
C>0
\,.
$$
The elements of the matrix
\eqref{eq:fundsolmatrix} are bounded above
polynomially in
$|\tau,\eta|$.
Therefore
 the kernel
\eqref{eq:integralkernel} is bounded above by
\begin{equation}
\label{eq:upper}
\
|\tau,\eta|^N\
\exp
\big(
-c(x_1-y_1)/|\tau,\eta|^N
\big)
\
\chi_{x_1\ge y_1}\,,
\qquad
c>0\,.
\end{equation}
This is proved  using Schur's Theorem to reduce
$M^\pm$ to upper triangular form by orthogonal transformations
of the spectral subspaces.  Then solve the differential equation
$X^\prime + M^+ X =0$ by back substitution
to prove
$
\|
\exp(\rho M^+)
\|
\le C|\tau,\eta|^p e^{-c \rho/|\tau,\eta|^N}
$.

The operator with kernel \eqref{eq:upper}
 is  convolution by an element
of $L^1(\RR)$ whose $L^1$ norm grows polynomially
in $|\tau,\eta|$.   By Young's theorem one concludes that
the operator with kernel \eqref{eq:integralkernel}
has norm in ${\rm Hom}(L^2(\RR))$ which grows at most
polynomially in $|\tau,\eta|$.

There is an entirely
analogous estimate for the expression for the spectrum
with  negative
real part.

Therefore,
$$
\big\|
\widehat W_I(\tau,x_1,\eta)
\big\|_{L^2(\RR)}
\ \le \
C_1\,
(1+|\tau| + |\eta|)^N\,
\big\|
\widehat G(\tau,x_1,\eta)
\big\|_{L^2(\RR)}
\ \le \
C_2\,
(1+|\tau| + |\eta|)^N\,
\big\|
\widehat F(\tau,x_1,\eta)
\big\|_{L^2(\RR)}\,.
$$
A similar estimate for $\widehat W_{II}$ follows from
\eqref{eq:uII}. Estimates for $\widehat U_j$ follow from the second equation in \eqref{eqn:tildetransformed}.
Plancherel's Theorem then implies
\eqref{eq:expestimate},
proving the existence part of well posedness.

Uniqueness is proved by a duality
argument of H\"olmgren type using existence
(backward in time)
for the adjoint differential operator (details omitted).
\end{proof}

\subsubsection{Perfection for  B\'erenger's PML with
variable coefficient
$\sigma_1(x_1)$}

\begin{lemma}
\label{lem:variableode}
Suppose that $A,M$ satisfy the hypotheses
of Lemma \ref{lem:blochform} with $\GG$ and
$\widetilde M\in {\rm Hom}\,\GG$
are from  that Lemma.
Suppose in addition that $f\in L^\infty_{loc}(\RR\,;\,\CC)$ and
$g$ is the unique solution of
$$
\frac{dg}{dx_1} =f,
\quad
g(0)=0,
\qquad
{\rm so},
\qquad
g(x_1)
\ =\
\int_{0}^{x_1} f(s)\ ds\,.
$$
Then for $\gamma\in \GG$ the unique solution of the equivalent
initial value problems for the $\GG$ valued function $U$,
$$
A\frac{dU}{dx_1}+ f(x_1)\,M\,U=0,
\quad
{\rm equivalently},
\quad
\frac{dU}{dx_1} + f(x_1)\,\widetilde M\,U \ =\ 0,
\qquad
U(0)=\gamma,
$$
is
$$
U(x_1)
\ =\
e^{-g(x_1)\,\widetilde M}\,\gamma\,.
$$
\end{lemma}

\begin{proof}  Compute using the differential equation,
$$
\frac{d}{dx_1}\Big[
e^{g(x_1)\, \widetilde M}\ U
\Big]
\ =\
e^{g(x_1)\widetilde M}
\Big[
\Big(\frac{dg}{dx_1}\Big)\widetilde M + \frac{dU}{dx_1}
\Big]
\ =\
e^{g(x_1)\widetilde M}
\Big[
f\, \widetilde M - f\widetilde M
\Big]\ =\ 0\,.
$$
The lemma follows.
\end{proof}

The next result shows that when
 the B\'erenger split problem  with
absorption
$\sigma(x_1)$ defines a stable time evolution, then the
problem is perfectly matched.
Either the split problem is ill posed, or it is well posed
and perfect.

\begin{theorem}  Suppose that $\sigma_1(x)\in L^\infty(\RR)$
has support in $[0,\rho]$  for  some $\rho>0$,  that $\sigma_j=0$ for
$j\ne 1$, and that the operator $\widetilde L$
with these absorptions is nondegenerate with respect to $x_1$ and defines a weakly well posed
time evolution.  Then, the $\widetilde L$ evolution
is perfectly matched in the sense that for
$F\in C^\infty_0(\{t>0\} \cap  \{x_1<0\})$
the solutions $\widetilde U$ and $\widetilde U^\prime$
with and without absorptions respectively,
\[
\widetilde L_1 \widetilde U \ =\
F\,,
\qquad
\widetilde U\big|_{t\le 0}\ =\ 0,
\qquad
{\rm
and,
}
\qquad
\widetilde L \widetilde U' \ =\
F\,,
\qquad
\widetilde U'\big|_{t\le 0} \ =\ 0
\]
satisfy
\[
\widetilde U\big|_{x_1<0}
\ =\
\widetilde U^\prime \big|_{x_1<0}\,.
\]
\end{theorem}

\begin{proof}
Denote by $\widehat U$ and $\widehat U'$ the Fourier-Laplace
transforms.   The functions are  characterized by
$$
\widetilde L_1(\tau, d/dx_1,\eta) \, \widehat U \ =\ \widehat F,
\qquad
{\rm
and,
}
\qquad
\widetilde L(\tau, d/dx_1,\eta) \, \widehat U' \ =\ \widehat F,
$$
both required to decay exponentially as $|x_1|\to \infty$.
The strategy is to construct a solution of the problem
defining $\widehat U'$ from the solution $\widehat U$.

The equations for $\widehat W=\sum_j \widehat U_j$ and $\widehat W'=\sum_j \widehat U'_j$ in $x_1\ge 0$
have the form
$$
 A_1
\frac{d\,\widehat W}{dx_1}
\ +\ M\, \widehat W\ = \ 0,
\qquad
A_1
\frac{d\,\widehat W'}{dx_1}
\ +\
\frac{\tau + \sigma_1(x_1)}{\tau}\,M\, \widehat W'\ = \ 0.
$$
Lemma \ref{lem:variableode}
 applies with $f(x):=(\tau+\sigma(x_1))/\tau$.

Define $g$  as in that lemma.  Set  $\widehat V=\widehat W$
in $x_1\le 0$.  For $x_1\ge 0$ define
$$
\widehat V\ :=\
e^{-g(x_1)\widetilde M} \widehat  W(\tau,0,\eta)\,.
$$
The resulting function satisfies the differential equation
required of $\widehat  W'$.  In addition since
$e^{-g(x_1)M}$ is independent of $x_1$
for $x_1\ge \rho$, $\widehat V$ decays as rapidly
as $\widehat W$.  Therefore $\widehat V$ satisfies the
conditions uniquely determining
$\widehat W'$. Therefore $\widehat V=\widehat W'$, and
$\widehat W'|_{x_1<0}=\widehat W|_{x_1<0}$.
    Use
\eqref{eqn:tildetransformed} to recover $\widehat U\,,\,\widehat U^\prime$ from
$\widehat W\,,\,\widehat W^\prime$
 shows that
 $\widehat U'|_{x_1<0}=\widehat U|_{x_1<0}$, proving
perfection.
\end{proof}

\begin{example}
\newline
{\bf 1.}
  If $L_1(0,\partial_x)$ is
elliptic then
Corollary \ref{thm:elliptic}
shows that the evolution of $\widetilde L$
is strongly well posed.  This includes the case
of anisotropic wave equations for which the
layer is amplifying showing that perfection
is not at all inconsistent with amplification.\\

\noindent
{\bf 2.}
  For the Maxwell equations and
$\sigma_1(x_1)\in W^{2,\infty}(\RR)$
well posedness is proved in
the remark following
Theorem \ref{th:vacus}
and we deduce perfection.
\end{example}

%%%%%%%%%%%%%%%%%%%%%%%%%%%
%  New section
%%%%%%%%%%%%%%%%%%%%%%%%%%
\section{Plane waves, geometric optics, and amplifying layers}\label{sec:amplification}

This section includes a series of ideas all
related to plane waves and short
wavelength asymptotic solutions of
WKB type.   We first recall the derivation
of such solutions from exact plane wave
solutions by Fourier synthesis.
Then we review the
construction of short wavelength asymptotic
expansions.  These are then applied to
examine the proposed absorption by the
$\sigma_j$.   In
many common cases the supposedly
absorbing layers lead to asymptotic
solutions which grow in the layer.
Related phenomena are
studied by Hu, and Becache, Fauqueux, Joly
\cite{HU:1996:ABC},  \cite{Becache:2003:SPM}.  For the Maxwell equations
for which the PML were designed, the
layers are not
amplifying.  At the end of
Section \ref{subsec:amplifying}, situations
where the amplification does not occur
are identified.

\subsection{Geometric optics by Fourier synthesis}
\label{subsec:synthesis}

When the coefficient $\sigma$ vanishes identically, both $L$ and $\til$ are homogeneous constant coefficient systems. When $(\utau,\uxi)$ is a smooth point of the characteristic variety, denote by $\tau=\tau(\xi) $ the smooth parameterization, and $\Pi_{L}\txi$ and $\Pi_{\til}\txi$ the associated spectral projections for $\xi\approx \uxi$, see \eqref{eq:proj}. The function $\tau(\xi)$ is homogeneous of degree 1, while the projectors are  homogeneous of degree 0.  The next argument works equally well for $L$ and $\til$.

For $G(\xi)\in C^\infty_0(\RR^d)$  construct exact solutions
for
$0<\eps << 1$,
$$
U^\eps(t,x)
\ :=\
\int
e^{i(\xi\cdot x+\tau(\xi) t)}
\
\Pi_{L}(\tau,\xi) \
G(\xi-\uxi/\eps)
\
d\xi\,.
%\eqno(2.1)
$$
Make the change of variable
$$
\xi-\uxi/\eps \ :=\ \
\zeta
\,,
\qquad
\xi=(\uxi  +\eps \zeta)/\eps
\,,
$$
and extract the rapidly oscillating term $e^{i(\uxi\cdot x +\utau t)/\eps}$ to find,
\begin{align}\label{eq:fa}
U^\eps(t,x)
&\ :=\
e^{i(\uxi\cdot x +\utau t)/\eps}
\
\int
e^{i\bigl((\tau(\uxi+\eps\zeta)-\tau(\uxi))t/\eps+\zeta x \bigr)}
\
\Pi_{L}(\tau(\uxi +\eps \zeta),\uxi +\eps \zeta)\
G(\zeta)
\
d\zeta\cr
&\ : =\
e^{i(\uxi\cdot x +\utau t)/\eps}
\  a(\eps,t,x)
\,.
\end{align}
Expanding in $\eps$
and keeping just the leading term yields the principal term in the geometric
optics approximation
$$
U^\eps \ \approx\
e^{i(\uxi\cdot x +\utau t)/\eps}
\
\int
e^{i(x\cdot\zeta - \bfv(\uxi)\cdot\zeta t)}
\
\Pi_{L}(\utau,\uxi)\
G(\zeta)
\
d\zeta
\,,
\qquad
\bfv(\uxi) \ := -
\partial_{\xi}\tau(\uxi)
\,.
%\eqno(2.3)
$$
One has
$$
U^\eps \ \approx\
e^{i(\uxi x +\utau t)/\eps}
\
a_0(x\,-\,\bfv(\uxi) t)
\,,
\qquad
a_0(x):=\int e^{ix\cdot\zeta}\,\Pi_{L}(\utau,\uxi)\,G(\zeta)\,d\zeta
\,.
%\eqno(2.4)
$$

A complete Taylor expansion yields the corrected
approximations which satisfy the equation with
an error $O(\eps^N)$ for all $N$.  We write
$O(\eps^\infty)$ for short.
This yields  infinitely accurate solutions,
\begin{equation}\label{eq:WKB}
U^\eps(t,x)
\ :=\
e^{i(\uxi x +\utau t)/\eps}
a(t,x,\eps)\,,
\qquad
a(t,x,\eps)\ \sim\
a_0(x-\bfv t) + \eps a_1(t,x) + \cdots \ \,.
%(2.5)
\end{equation}

If 0 is a semi-simple eigenvalue of $L(\tau,\xi)$, and $\Phi_0\in \ker\,L(\tau,\xi)\setminus \{0\}$ of dimension 1,  then the leading amplitude $a_0$ in the case of \eqref{hyp} (resp. \eqref{PMLgen}) is of the form
$$
\alpha(t,x)\,\Phi_0\,,
\qquad
\biggl(
{\rm resp.} \quad
\alpha(t,x)
\bigl(\frac{\xi_1}{\tau}A_1\Phi_0,\ldots,\frac{\xi_d}{\tau}A_d\Phi_0\bigr)
\biggr)
%\eqno(2.6)
$$
with scalar valued amplitude $\alpha$ satisfying,
$$
\big(
\partial_t +\bfv.\partial_x\big)\alpha\ =\ 0\,.
%\eqno(2.7)
$$
This shows that $\alpha$ is constant on the {\bf rays}
which are lines with velocity equal to the
group velocity $\bfv(\xi):=-\partial_{\xi}\tau(\xi)$.

For $g\in C^\infty_0\setminus 0$ the solutions
do not
have compact spatial support.  This weakness is easily
overcome.  Choose $\chi\in C^\infty_0(\RR^d)$
with $\chi=1$ on a neighborhood of the origin.
For  $g\in \caS(\RR^d)$,
define exact solutions by cutting off the integrand
outside the domains of definition of $\tau(\xi)$
and $\Pi_{L}(\tau(\xi),\xi)$,
\begin{equation}\label{eq:fourier}
u^\eps(t,x)\ :=\
\int e^{i(\xi x+\tau(\xi) t)}
\
\Pi_{L}(\tau(\xi),\xi)\
g(\xi-\uxi/\eps)
\ \chi(\sqrt \eps(\xi-\uxi/\eps))\
d\xi\,.
\end{equation}
The analysis above applies with the only change being the
initial values.  In the preceding case these values were equal to the
transform of $\Pi_{L}(\tau(\xi),\xi)\,g(\xi-\uxi/\eps)$ and in the present case they
are infinitely close to that
quantity,
\begin{equation*}
u^\eps(0,x) \ =\
\int
e^{ix.\xi}\ \Pi_{L}(\tau(\xi),\xi)\ g(\xi-\uxi/\eps)\ d\xi
\ +\
O(\eps^\infty)\,.
\end{equation*}
This yields infinitely accurate approximate
solutions (2.5) which have support in the tube of rays
with feet in the support of
$\int
e^{ix.\xi}\ \Pi_{L}(\tau(\xi),\xi)\ g(\xi)\ d\xi$.

\subsection{Geometric optics with variable coefficients}
\label{subsec:WKB}
The Fourier transform method of the preceding sections
is limited to problems with constant coefficients.  In this
section the WKB method which works for
variable coefficients is introduced.  It will also serve
for the analysis of reflected waves.

Let $\cl$ be the general operator in \eqref{gen}. Fix $\txi\in   {\rm Char} \,\cl$ and seek asymptotic solutions
\begin{equation}\label{eq:expansion}
    U^{\eps}\ \sim\  e^{iS/\eps} \sum_{j=0}^{+\infty}
    \eps ^j a_j(t,x),
    \qquad \mbox{ with the phase }
    S(t, x,\xi)\ =\
    t\tau  + x\xi.
  \end{equation}
More precisely we construct smooth functions $a_j(t,x)$ with ${\rm supp}\, a_j\ \cap ([0,T]\times\RR^d)$ compact so that if
  $$
  a(t,x,\eps)\ \sim\
  \sum_{j=0}^\infty\ \eps^j\ a_j(t,x)\,,
  $$
  in the sense of Taylor series at $\eps=0$,
  and
  $ {\rm supp}\, a\ \cap ([0,T]\times\RR^d\times ]0,\eps])$
  is  {compact},
  then
  $$
  U^\eps\ :=\
  e^{iS/\eps} \
  a(t,x,\eps)
 $$
 satisfies for all $s,N$,
 $$
\big\|\cl\,
U^\eps
\big\|_{H^s([0,T]\times\RR^d)}
\ =\
O(\eps^N)\,.
$$
In this case we say that \eqref{eq:expansion}
is a {\sl infinitely accurate approximate solution}.
The next result recalls some
facts about such solutions.
\begin{theorem} \label{th:stabilitegen}
  Suppose Problem \eqref{gen} is hyperbolic,  $(\tau,\xi)\in {\rm Char}(\cl)$ satisfies the smooth variety hypothesis, and $0$ is a semi-simple eigenvalue of $\cl(\tau,\xi)$.
\begin{enumerate}[label={\it(\roman{*})}, ref={\it(\roman{*})},leftmargin=0.7 cm]
\item   If  the coefficients $a_j$
  satisfy the recursion relation
  \begin{subequations}\label{eq:rec}
    \begin{align}
      &a_0(t,x) \in \ker\,\cl_1(\tau,\xi),\label{eq:rec1}\\[1mm]
      &\forall j \geq 0, \quad
      i\cl_1(\tau,\xi)\, a_{j+1}(t,x)+ \cl(\partial
      _t,\partial
      _x)\,a_j(t,x)=0, \label{eq:rec2}
    \end{align}
  \end{subequations}
  then \eqref{eq:expansion} is an
  infinitely accurate approximate solution of \eqref{gen}.
 \item   If $g_j= \Pi_\cl(\tau,\xi)g_j\in C^\infty_0(\RR^d_x)$
 are supported in a fixed
 compact $K$,
 then
  there is one and only one family of $a_j$ satisfying \eqref{eq:rec}
 together with the initial conditions, $\Pi_\cl(\tau,\xi)a_j(0,\cdot)=g_j$
 and the polarization $\Pi_\cl(\tau,\xi)a_0=a_0$.
 They have support in the
 tube of rays with feet in $K$ and speed of
 propagation $\bfv(\xi)=-\partial_\xi\tau(\xi)$.
\item   The principal term $a_0$ is a solution of the transport equation
  \begin{equation}\label{eq:transport}
    \begin{array}{c}
      \ds\partial _ta_0 +\bfv(\xi) \cdot\partial_x a_0+\Pi_\cl(\tau,\xi){\cal B}(x)\Pi_\cl(\tau,\xi) a_0=0.
    \end{array}
  \end{equation}
\end{enumerate}
\end{theorem}

\begin{proof}
For simplicity note $\Pi_\cl:=\Pi_\cl\txi$ when no ambiguity is to be feared.

The equations \eqref{eq:rec} are obtained by injecting $U^{\eps}$ in (\ref{gen}), to find an expression $\sim e^{iS/\eps}\sum_{j\ge 0} \eps^j w_j(t,x)$. In order that the $w_j$ vanish it is necessary and sufficient that the equations \eqref{eq:rec}  are satisfied.

Next examine the leading order terms to find the relations determining $a_0$. Projecting the case $j=0$ of \eqref{eq:rec}   onto $\ker\cl_1$ yields,
$$
\Pi_{\cl}\, (\partial _t+\sum_{l=1}^d {\cal A}_l\partial _{l} +{\cal B}(x))\,a_0=0.
$$
This yields a first order system satisfied by  $a_0= \Pi_{\cl}\,a_0$,
\begin{equation} \label{eq:transp1}
  \partial _t a_0 +\sum_{l=1}^d \Pi_\cl\,  {\caA}_l\,
  \Pi_\cl\, \partial _{l}
  a_0+\Pi_\cl\, {\cal B}(x)\,\Pi_\cl\, a_0=0.
\end{equation}
The leading order part of this equation is
a  scalar transport operator.  To see this
differentiate
$\cl_1(\tau(\xi),\xi)\Pi_{\cl}(\tau(\xi),\xi) =0$
with respect to
$\xi_l$ to find
$$
\left( {\cal A}_l+\frac{\partial \tau(\xi)}{\partial
    \xi_l} Id \right)\Pi_{\cl}(\tau(\xi),\xi)
    +\cl_1(\tau(\xi),\xi)    \frac{\partial\ }{\partial\xi_l}
  \left(\Pi_{\cl}(\tau(\xi),\xi)\right) =0.
  $$
Multiplying on the left by $\Pi_{\cl}(\tau(\xi),\xi)$
eliminates  the second term yielding,
$$
  \Pi_{\cl}{\cal A}_l\Pi_{\cl}+
  \frac{\partial \tau(\xi)}{\partial
    \xi_l} \Pi_{\cl}=0.
$$
Injecting this in  \eqref{eq:transp1} yields \eqref{eq:transport}.

In order to compute the coefficients recursively,  multiply \eqref{eq:rec2} on the left by the partial inverse $Q_{\cl}\txi$, using the identity in \eqref{eq:inverse},
to obtain for $j \ge 1$,
\begin{equation}\label{eq:recur1}
(I-\Pi_{\cl})a_j
\ =\
i\, Q_\cl\, \cl(\partial
      _t,\partial
      _x)\,a_{j-1}.
\end{equation}
Projecting \eqref{eq:rec2} on the kernel yields,
\[
\Pi_{\cl}\,\cl(\partial
      _t,\partial
      _x)\,a_j=0\,.
\]
Writing $a_j$ as
\[
a_j=\Pi_{\cl} \,a_j + (I-\Pi_{\cl})\,a_j,
\quad
{\rm yields}
\quad
\Pi_{\cl}\,\cl(\partial_t,\partial_x)\,\Pi_{\cl}\,  a_j=
      -\Pi_{\cl} \,\cl(\partial
      _t,\partial
      _x)\,(I-\Pi_{\cl})\,a_j.
\]
This is again a transport equation, but with a righthand side,
\begin{equation}\label{eq:recur2}
      \ds\partial _t\Pi_{\cl}a_j +\bfv \cdot\partial_x \Pi_{\cl}a_j+\Pi_\cl\, B\,\Pi_\cl\, a_j=-\Pi_{\cl}\,\cl(\partial
      _t,\partial
      _x)\,(I-\Pi_{\cl})\,a_j .
\end{equation}
\eqref{eq:recur1} and \eqref{eq:recur2} permit to calculate the coefficients recursively, knowing the initial values.
\end{proof}

Next apply the above algorithm
 to the PML operator
 $\widetilde L$. Fix $\txi\in  {\rm Char} \,L$ and seek asymptotic solutions
\begin{equation}\label{eq:expansionpml}
    \tiu^{\eps}\ \sim\  e^{iS/\eps} \sum_{j=0}^{+\infty}
    \eps ^j \tilde{a}_j(t,x),
    \qquad \mbox{ with the phase }\quad
    S(t, x)\ =\
    t\tau+ x\cdot\xi.
  \end{equation}
\begin{corollary} \label{th:stabilitepml}
  Suppose Problem \eqref{hyp} is strongly well posed,  $\txi \in {\rm Char} L$ satisfies the smooth variety hypothesis, and $0$ is a semi-simple eigenvalue of $L\txi$.
\begin{enumerate}[label={\it(\roman{*})}, ref={\it(\roman{*})},leftmargin=0.7 cm]
\item   If  the coefficients $\tilde{a}_j$
  satisfy the recursion relation
  \begin{subequations}\label{eq:recpml}
    \begin{align}
      &\tilde{a}_0(t,x) \in \ker\,\til_1\txi,\label{eq:recpml1}\\[1mm]
      &\forall j \geq 0, \quad
      (I-\Pi_{\til})\,\tilde{a}_j(t,x)=i Q_{\til} \,\til(\partial_t,\partial_x)\,\tilde{a}_{j-1}(t,x), \label{eq:recpml2}\\
      &\ds\partial _t\Pi_{\til}\,\tilde{a}_j +\bfv \cdot\partial_x \Pi_{\til}\,\tilde{a}_j
      +\beta(x)\,\Pi_{\til} \, \tilde{a}_j=-\Pi_{\til}\,\til(\partial
      _t,\partial
      _x)(I-\Pi_{\til})\,\tilde{a}_j.\label{eq:recpml3}
    \end{align}
  \end{subequations}
  then \eqref{eq:expansionpml} is an
  infinitely accurate approximate solution of \eqref{PMLgen}.

 \item   If $\tilde{g}_j(x)= \Pi_{\til}\,\tilde{g}_j\in C^\infty_0(\RR^d)$
 are supported in a fixed
 compact $K$,
 then
  there is one and only one family of $\tilde{a}_j$ satisfying \eqref{eq:recpml}
 together with the initial conditions, $\Pi_{L}\txi\tilde{a}_j(0,x)=\tilde{g}_j$
 and the polarization $\Pi_{L}\txi\tilde{a}_0=\tilde{a}_0$.
 They have support in the
 tube of rays with feet in $K$ and speed of
 propagation $\bfv=-\partial_\xi\tau(\xi)$.

\item   The principal term $\tilde{a}_0$ is a solution of the transport equation
  \begin{equation}\label{eq:transportpml}
    \begin{array}{c}
      \ds\partial _t\tilde{a}_0 +\bfv \cdot\partial_x \tilde{a}_0+\beta(x) \tilde{a}_0
      \ =\ 0,
      \quad\mbox{ with }\quad
      \ds\beta(x)
      \ =\
      \sum_{l=1}^d
      \frac{\sigma_l(x_l) \xi_l}{ \tau(\xi)}\
       \frac{\partial
        \tau(\xi)}{\partial \xi_l}.
    \end{array}
  \end{equation}
\end{enumerate}
\end{corollary}
\begin{proof}
  We need only identify the constant term in \eqref{eq:transport}.
  Use the form of the projector given in Proposition  \ref{prop:projection}, to obtain
  \begin{equation*}
  %\label{eq:idinterm}
     \Pi_{\til} \, B(x)\,\Pi_{\til}\,   =\beta(x) \Pi_{\til}\,.
  \end{equation*}
\end{proof}
\subsection{Amplifying layers}\label{subsec:amplifying}

The coefficient $\sigma_1(x_1)\ge 0$ is introduced
with
 the idea that waves will be damped in the
 layer.
In this section, we show that sometimes the anticipated
decay is not achieved, and waves may  be
amplified.  This was observed in \cite{Becache:2003:SPM}.
The authors
analysed the phenomenon
for $\sigma$ constant in the layer.
They showed
that in an infinite layer solutions
can in certain cases grow
infinitely large.
We present a related analysis
 using
WKB solutions  which has three advantages,

\vskip.1cm

\noindent
{\bf 1.}  The analysis is valid for variable coefficients
$\sigma_1(x_1)$ which corresponds to common practice.

\noindent
{\bf 2.}  The growth is seen immediately and not
expressed in terms of large time asymptotics.

\noindent
{\bf 3.}  The  analysis in  \cite{Becache:2003:SPM} was in part restricted to
$d=2$ and eigenvectors of multiplicity one.  We
remove these restrictions.
\vskip.1cm

\noindent
It is because of {\bf 2} that we  choose not to follow the authors
of \cite{Becache:2003:SPM} in calling this phenomenon
instability.

\begin{theorem}
Suppose $(\tau,\xi) \in {\rm Char}(L)$ satisfies
 the smooth variety hypothesis and $\beta(x)$ is
 as in
  \eqref{eq:transportpml}. Suppose in addition there is an interval
  on a ray
$$
\Gamma:=\{ (0,\ux)+t(1,-\partial_\xi\tau(\xi)),\quad 0 \le t \le t_0\},
\qquad
{\rm so\ that},
\qquad
 \int_0^{t_0}
\beta(\ux-t\partial_\xi\tau(\xi))) \ dt\  <\  0.
$$
Then the corresponding WKB solution grows in the layer.
\end{theorem}
\begin{proof}
  The solution of the transport equation
  \eqref{eq:transportpml} evaluated on $\Gamma$ is
  \begin{equation*}\label{transportexpl}
    \tilde{a}_0(t_0,x) = \exp \left( -\int_0^{t_0}
      \beta(\ux-s\partial_\xi\tau(\xi))))ds \right)
    \tilde{a}_0(0,\ux+t_0\partial_\xi\tau(\xi)) .
  \end{equation*}
  The exponential is strictly greater than 1, so
  \begin{equation*}\label{ampl}
    |\tilde{a}_0(t_0,\ux+t_0\partial_\xi\tau(\xi))| > |a_0(0,\ux )| .
  \end{equation*}
\end{proof}

\begin{example}[No amplification for Maxwell/D'Alembert]
If the dispersion relation is $\tau^2 =|\xi|^2$
and  $\sigma\ge 0$ then there is no amplification since,
$$
\beta \ =\ \sum_{j=1}^d \sigma_{j}(x_j)\frac{\xi_j^2 }{\|\xi\|^2}\ \ge \ 0\,.
$$
\end{example}

\begin{example} [Amplification is common]
For the dispersion relation
$\tau^2=q(\xi)$ where $q$ is a positive definite
quadratic form so that the $\xi$ axes are not
major and minor axes of the ellipse $q=1$,
 there are always $\tau>0,\xi$ so that
 $x_1$  layers with $\sigma_1>0$ are
amplifying (\cite{Becache:2003:SPM}).
There are two lines on $\{\tau=q(\xi)^{1/2}\}$  where
$\partial q/\partial \xi_1=0$.  The half cone on which
$\partial q/\partial \xi_1 <0$ corresponds to rays
on which $x_1$ is increasing so they enter a layer
$x_1>0$. The half cone
$\{\partial q/\partial \xi_1 <0\}$ is divided into two
sectors by the plane $\xi_1=0$.  The sector on which
$\xi_1>0$ (resp. $\xi_1 < 0$) corresponds to growing
(resp. decaying) solutions (see Figure \ref{fig:ampli} on the left).  This example
shows that amplification is very common.
\begin{figure}[ph]
\begin{center}
\psfrag{x1}{$\xi_1=0$}
\psfrag{outgoing}{\small outgoing}
\psfrag{v}{outgoing, $\frac{\partial q}{\partial \xi_1} < 0$}
\includegraphics[width=0.50\textwidth]{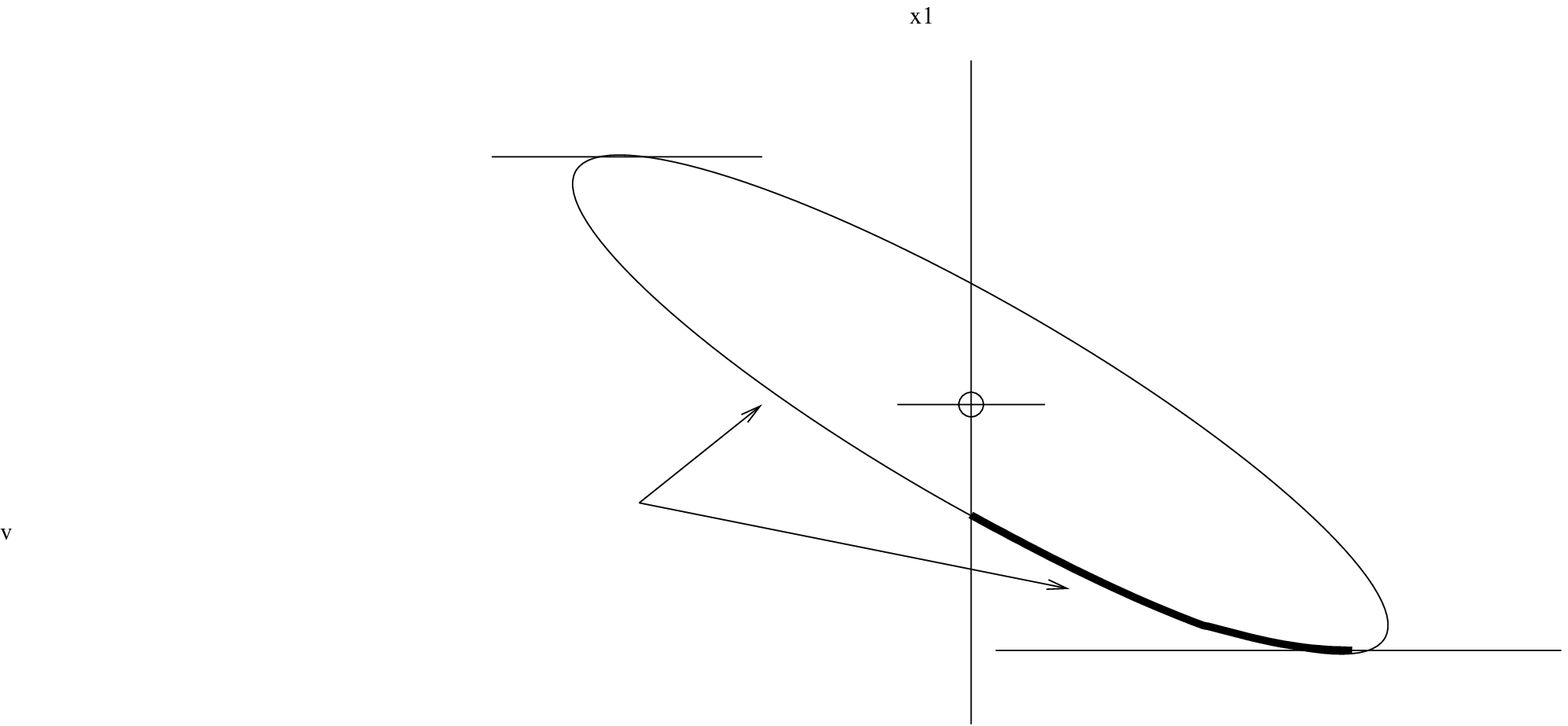}
\qquad\qquad\qquad
\includegraphics[width=0.18\textwidth]{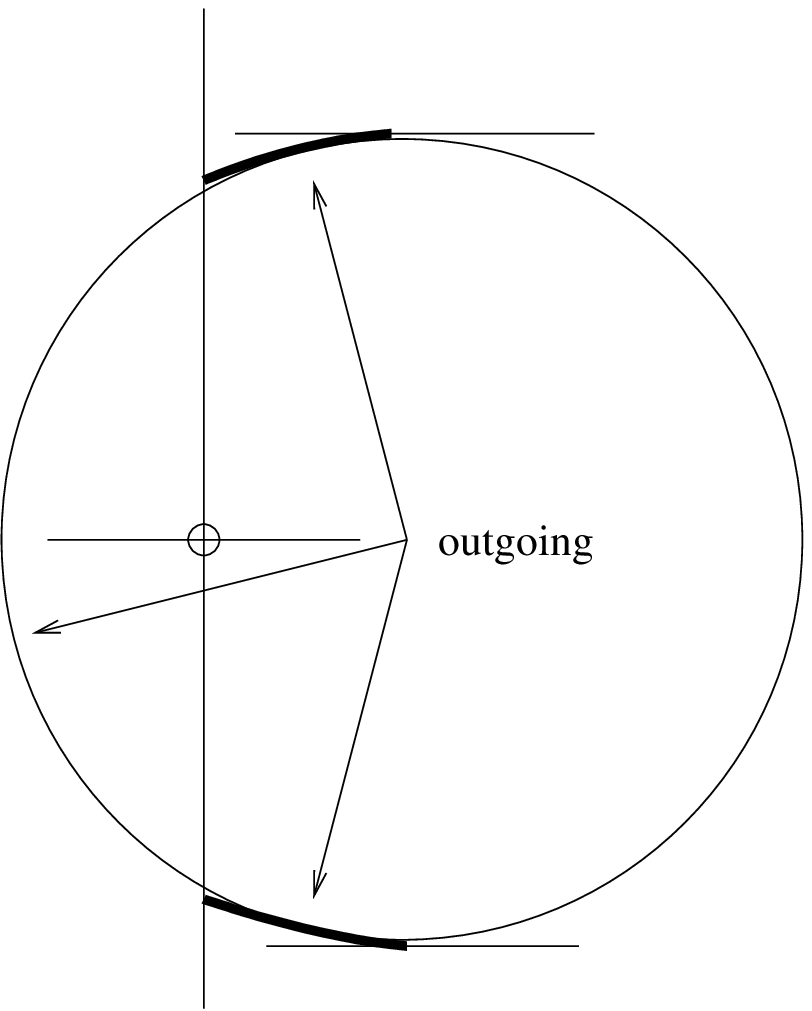}
\caption{Amplified outgoing wave numbers in bold \label{fig:ampli}}
\end{center}
\end{figure}
Consequently for the dispersion relation $\tau^2=q(\xi)$
it is wise to align coordinates along the major and minor
axes of the ellipse to avoid amplification.  However,
if  $(\tau^2-q_1)(\tau^2-q_2)$ divides the characteristic
polynomial and the axes of $q_1$ and $q_2$ are distinct
from each other then no linear change of coordinates
can avoid amplification in the layer.

A second example from \cite{Diaz:2006:TDA} is the linearized
compressible Euler equation with nonzero background
velocity $(c,0), c>0$ for which amplified wave numbers
at a right hand boundary are indicated in bold
in the right hand figure \ref{fig:ampli}.
\end{example}

\noindent
{\bf Summary.}  There is no amplification when the
characteristic polynomial is a product of factors
$\tau$ and $\tau^2-q$ where $q$ is a positive definite
quadratic form with axes of inertia parallel to the coordinate
axes.  This includes the cases of Maxwell's equations
in vacuum,
for which the method was developed by B\'erenger,
 the
linearized Euler equations about the stationary state,
and the linear isotrope elasticity equations.  For  these
the quadratic forms $q$ are multiples of $|\xi|^2$.

\begin{example}[Methods related to B\'erenger, continued]
For the model developed in Section \ref{subsub:PMRB}
for the Maxwell equations,  one can compute
%We only need to compute $\Pi_R B\Pi_R$. Since $\ker L\txi$ is monodimensional, so is $\ker R_1\txi$. Therefore $\Pi_R=\Phi\lmj{R}^0{\Phi\lmj{R}^0}^*/{\Phi\lmj{R}^0}^*\Phi\lmj{R}^0$, with
%\[
%\Phi\lmj{R}^0\txi=
%\begin{pmatrix}
%\ds \Phi\txi\\
%0\end{pmatrix}
%.\]
%Then
%\[
%\Pi_R B\Pi_R=\ds\frac{{\Phi\lmj{R}^0}^* B {\Phi\lmj{R}^0}}{{\Phi\lmj{R}^0}^*{\Phi\lmj{R}^0}} \Pi_R
%:=\gamma \Pi_R\,.
%\]
%Compute now $\beta$:
%\[
%\beta={\Phi}^*\Sigma^{(1)}{\Phi}
%\]
%Inserting the value of $\Phi$ into the previous formula gives
\[
\beta
\ =\
 2\
\frac{\xi_2^2+\xi_3^2}{\xi_1^2}
\
\ds \sum_{j=1}^3\ \sig_j(x_j)\xi_j^2\  >\ 0\,.
\]
Thus,
 this model has exactly the same good properties as B\'erenger's, and is strongly well posed. For Maxwell equations, it is therefore an attractive alternative. The advantage is two fold.   The system with the auxiliary variables is very  compact.
And it is strongly well posed, even for discontinuous
$\sigma$.

On the surface this result
 sounds almost too good to be true.
However the B\'erenger system in the case of Maxwell's
equations has almost exactly the same
structure.  The energy method proof when
$\sigma_j^{\prime\prime}\in L^\infty$ shows that
there is a large vector $\VV$ consisting of the components
of $U$ together with differential operators $P_\alpha(D)$
applied to $U$ and a strongly well posed equation
for $\VV$.   This means that if one were to introduce the
additional variables in $\VV$ one obtains a system
with some of the desirable properties of SPML (strong PML).  However,
the SPML reduction is much more compact, and, has a good energy
estimate even when $\sigma$ is discontinuous.
The extension of this
 strategy to other equations is not straightforward.
 For elastodynamic models, see \cite{Rahmouni:2004:WPP}.

\end{example}

%%%%%%%%%%%%%%%%%%%%%%%%%%%
%  New section
%%%%%%%%%%%%%%%%%%%%%%%%%%
\section{Harmoniously matched layers}
\label{sec:hml}

This section introduces a new absorbing layer method.
It is based on the
following strategy.
Start with an operator $L=L_1(\partial)$ on the left
and consider a smart layer on the right
\begin{equation}
\label{eq:newC}
R(t,x,\partial)= L_1(\partial) + C(t,x),
\qquad
C\ =\
\sigma(x_1)\Big(\pi_+(A_1) + \nu \pi_-(A_1)\Big),
\qquad
{\rm supp}\,\sigma \subset[0,\infty[,
\end{equation}
generalizing
 \eqref{eq:smartB}.
This method is  embedded in a family of absorbing layers parameterized by
$\mu\ge 0$,
\begin{equation}
\label{eq:Lmu}
R^\mu
\ :=\
L_1 \ +\ \mu\, C\,.
\end{equation}
The method is nonreflective when $\mu=0$ and
is both reflective and dissipative for $\mu> 0$.
When $\sigma$ is discontinuous the leading order reflection
coefficient for wave packets of amplitude $1$
oscillating as  $e^{i(\tau t + x\xi)/\eps}$
is of the form
$\eps\, \mu \,r(\tau,\xi)$.
The leading order reflections can be removed by
an extrapolation method using two values of $\mu$.
This simultaneously removes the leading reflections at all
angles of incidence.  We call the resulting method the
{\it harmoniously matched layer}.

\subsection{Reflection  is linear in $\mu$ by scaling}
\label{subsec:relinear}

In this subsection the linearity in $\mu$
of leading order reflections by the layer with $R^\mu$ from
\eqref{eq:Lmu}
is demonstrated by
a scaling argument when
 $\sigma(x_1)={\bf 1}_{x_1>0}$.
In the next subsection
the reflection is computed exactly for Maxwell's equations
yielding additional information.

If $(L_1 + C)U=0$ then
$$
\underline U(t,x) := U(\mu t,\mu x),
\qquad
{\rm satisfies}
\qquad
R^{\mu}\,\underline U \ =\ 0\,.
$$

Suppose that  $U$ has an incoming wave of wavelength
$\eps$ and reflected waves $U_\ell$ with
amplitudes $\rho_\ell \eps$.
Then $\underline U$ has an incoming wave with
wavelength $\underline \eps := \eps/\mu$.
The reflected waves have amplitudes
$$
 \rho_\ell \, \eps \ = \
 \rho_\ell\,\mu \ \frac{\eps}{\mu}
 \ = \
  \rho_\ell\,\mu \
  \underline \eps
 \,.
$$
Denote by
$\underline \rho_\ell$
the reflection coefficient of $R^{\mu}$.
The leading
amplitude of the reflected $\ell$ wave
is then  $\underline \rho_\ell \underline\eps$.
The preceding identity shows that
$\underline \rho_\ell=  \rho_\ell\mu$
showing that the reflection coefficients are  linear
in $\mu$.
\footnote{This argument
can be made rigorous under the following
conditions.  The incoming wave is a wave packet
with oscillatory part $e^{i(\tau t + x\xi)/\eps}$
with $(\tau,\xi)\in {\rm Char} L$.  Denote
$(\tau,\xi^\prime)$ the part determining the oscillations
in $x_1=0$. Consider  the roots
$\xi_1$ of each of the equation,
$
{\rm det}\, L_1(\tau,\xi_1, \xi^\prime)= 0$.
The nonreal roots are called {\it elliptic}.  They lead to
waves which have the structure of a boundary
layer of thickness $\sim \eps$.  The real roots are called
{\it hyperbolic}.  The favorable situation is when all the hyperbolic
roots are at smooth points of the characteristic variety
and the group  velocities are transverse to
the boundary.   In that case one
can construct infinitely accurate asymptotic solutions
of the transmission problem consisiting of incoming, reflected,
and transmitted wave packets corresponding to the hyperbolic
roots,  and, a finite
number of boundary layers corresponding to  elliptic roots.
 As this is a long story,  we content ourselves
with the Maxwell computation of the next subsection.}

\subsection{Reflection for Maxwell with smart  layers}\label{subsec:pmlgeometricoptics}

In this section $\cl$ may denote one of two
distinct operators.  One option is  the Maxwell operator $L_1$
from
\eqref{matricemaxwell2D} for the  $\CC^3$ valued field
$E+iB$.
The  lower order term is $\caB:=\mu C$ from the smart layer
\eqref{eq:newC}.
Alternatively $\cl$ may denote  the B\'erenger operator
operator
$\til$ with lower order term $\caB=\mu\, C$
with
$$
C =\begin{pmatrix}
  \sigma(x_1)\,I_N&0 &\dots & 0 \\
  \vdots &0&  \ddots & \vdots \\
  0 & \dots &0& 0
  \end{pmatrix}\,,
  \qquad
  {\rm supp}\, \sigma\subset [0,\infty[\,.
  $$
  In both cases the absorption term is linear in
  $\mu$.  We compute the dependence of the
  reflection coefficient on $\mu$.

Lemma \ref{lem:trace} shows that
  the Cauchy problem is equivalent to homogeneous problems
   in each half-space with a transmission condition on $\Gamma:=\{x_1=0\}$,
 \begin{equation}
 \label{eq:transcond}
\big[\caA_1  \,\cau \big]_{\Gamma}\ =\ 0\,.
 \end{equation}
 {\it In order to cover both cases
 the operator, coefficients, and unknown are indicated with
 round letters.}

We study the reflection of
high frequency waves in $x_1\le 0$ which
approach the boundary $x_1=0$.
The input is
an incident wave with phase
$
S^I(t,x) := \tau t + \xi ^\cdot x$,
where  $\tau \ne 0$ and  $\tau(\xi)=\pm |\xi|$. The phase is chosen so that the group velocity $\bfv=-\xi/\tau$ satisfies $v_1>0$.
Denote $x':=(x_2,x_3)$, $\xi'=(\xi_2,\xi_3)$.
Theorem \ref{th:stabilitegen} applies to the incident wave with ${\cal B}\equiv 0$.  In $x_1\le 0$,
the incident  wave is
\begin{equation}
\label{eq:incidente}
\cau^\eps :=e^{iS^I(t,x)/\eps}\,a \imj{I}(t,x,\eps) \,,\qquad
a \imj{I}(t,x,\eps)\sim\sum_{j=0}^\infty
\ \eps^j\,a \imj{I}_j(t,x)
,\qquad
\cl_1(\partial_t,\partial_x)\,{\cal U}^\eps = O(\eps^\infty)
\,.
\end{equation}
%The first term of the expansion satisfies,
%\[
%\big(
%\partial_t + \bfv.\partial_x
%% + \sigma h(x_1)
%\big)a \imj{I}_0
%\ =\ 0,\qquad a \imj{I}_0\in \ker L_1(\tau,\xi).
%\]
%
%
Suppose that the
amplitudes $a\imj{I}_j$ are supported in a tube, $\cal T$,  of rays
with compact temporal crossections
${\cal T}\cap \{t=0\}\subset\subset \{x_1<0\}$.

We construct  a transmitted wave with the same
phase, and a reflected wave with phase
$S\imj{R}(t,x):=
\tau t +\xi\imj{R}x$, with $\xi\imj{R}:=(-\xi_1,\xi')$.
We first show that there are uniquely determined
reflected and transmitted waves.  Then we compute
exactly the leading terms in their asymptotic expansions.

The  reflected wave $\cav^\eps$ is also supported
in $x_1\le 0$.  The group velocity for the
reflected wave  is equal to $\bfv\imj{R}:=(-v_1,v')$, and
in $x_1\le 0$,
\begin{equation}
\label{eq:reflechie}
{\cal V}^\eps=e^{iS\imj{R}(t,x)/\eps}\,a \imj{R}(t,x,\eps),
\qquad
a \imj{R} (t,x,\eps)\ \sim\
\ds\sum_{j=0}^\infty
\ \eps^j\, a \imj{R}_j(t,x),
\qquad
\cl_1(\partial_t,\partial_x){\cal V}^\eps = O(\eps^\infty)
\,.
\end{equation}
%The first term  of the expansion satisfies,
%\[
%\big(
%\partial_t + \bfv^R\cdot\partial_x
%% + \sigma h(x_1)
%\big)a \imj{R}_0
%\ =\ 0,\quad a\imj{R}_0\in \ker L_1(\tau,-\xi_1,\xi').
%\]
%%
The  transmitted wave
is supported  in $x_1\ge 0$,
\begin{equation}
\label{eq:transmise}
{\cal W}^\eps:=e^{iS^I(t,x)/\eps}\,a\imj{T}(t,x,\eps),
\qquad
a \imj{T}(t,x,\eps) \ \sim\
\ds \sum_{j=0}^\infty
\ \eps^j\, a \imj{T}_j(t,x)
,\qquad
\cl(\partial_t,\partial_x){\cal W}^\eps = O(\eps^\infty)
\,.
\end{equation}
%with,
%\[\big(
%\partial_t + \bfv\cdot\partial_x
% + \,\sigma\, \Pi_{L}C\Pi_{L}\,
%\big)a \imj{T}_0
%\ =\ 0,\quad a \imj{T}_0\in \ker L_1(\tau,\xi).
%\]

\begin{theorem}\label{th:transmission}
\leavevmode\par
\begin{enumerate}
[label={\it(\roman{*})}, ref={\it(\roman{*})}]

        \item   Given the  incoming amplitudes $a\imj{I}_j$
there are uniquely determined amplitudes $a\imj{T}_j$
and $a \imj{R}_j$ so that for any choice
of the $a^{I,R,T}(t,x,\eps)\sim \sum \eps^j\, a^{I,R,T}(t,x)$,
the $\cau^\eps$, ${\cal V}^\eps$ and ${\cal W}^\eps$ are infinitely accurate
solutions of the differential equations
and the transmission condition is also satisfied to infinite order,
\begin{equation}\label{eq:transmissioncondition2}
\forall (t,x') \in \RR\times \RR^2,
\qquad
\caA_1(\cau^\eps  + \cav^\eps )(t,0\lmj{-},x')
\ =\
\caA_1(\caw^\eps)(t,0\lmj{+},x') + O(\eps^\infty)\,.
\end{equation}

        \item  In the case of B\'erenger's  PML, the coefficients $\tilde{a}_j^R$  vanish identically for $j \ge 0$.

        \item For the smart layer (\ref{eq:newC}, \ref{eq:Lmu}) with $\sig={\bf 1}_{x_1>0}$,
    the coefficient $a_0^R$ vanishes identically. The reflection coefficient of the layer is equal to
$$
R(\tau,\xi)
\ =\
 i \mu\, \ds(1+\nu)\
\frac{\xi_1^2-\tau^2}{8\tau\xi_1^2}
\ =\
i \frac{\mu(1+\nu)}{8\tau}\ds
\frac{v_1^2-1}{v_1^2} \,.
$$
 That is, if $a\imj{I}_0(t,0\lmj{-},x')=\alpha(t,x')\Phi\txi\in  \kl$,
 then
 $$
 a\imj{R}(t,0\lmj{-},x')\ =\
\eps R\txi \,\alpha(t,x')\,
\Phi(\tau,\xi\imj{R})
\ +\
 O(\eps^2)\,.
   $$
Furthermore, the amplitudes $a\imj{T,R}$ are such that on the interface $\Gamma$, we have for all $(i,j)\in \NN^2$, $i=0$ and $j \le 1$, or $i \ge 1$ and $j\ge 0$,\footnote{$\CC_j[\mu]$ denotes the space of
polynomials of degree less than or equal to $j$ with complex
coefficients.  $\CC_j[\mu]\otimes \CC^3$ is the corresponding
 space of polynomials with coefficients in $\CC^3$.}
\begin{equation}\label{eq:linear}
\partial_1^j a\imj{T}_i-\partial_1^ja\imj{I}_i
\ \in\
\mu\,\big(\CC_{i+j-1}[\mu]\otimes\CC^3\big),\qquad
\partial_1^j a\imj{R}_i
\ \in\
\mu\,\big(\CC_{i+j-1}[\mu]\otimes\CC^3\big)\,.
\end{equation}
\item The smart layer with $\sigma(x_1)$ satisfying
      $\sig(0)=\cdots=\sig^{(k-1)}(0)=0$, $\sig^{(k)}(0)\ne 0$ is nonreflecting at order $k$ for any angle of incidence, \textit{i.e.}
if $a\imj{I}_0(t,0\lmj{-},x')=\alpha(t,x')\Phi\txi$,
 there exists $R_k\txi$ such that
 $$
 a\imj{R}(t,0\lmj{-},x')\ =\
\eps^k \sig^{(k)}(0) \, R_k \txi \,\alpha(t,x')\,
\Phi(\tau,\xi\imj{R}) \ +\
 O(\eps^{k+1})\,.
   $$
Furthermore the amplitudes $a\imj{T,R}$ are linear functions of $\mu$ on the interface $\Gamma$.  That is denoting
$c_i\imj{R}(\mu)= {a\imj{R}_i}\trace{\Gamma}$
and
$c_i\imj{T}(\mu)= {a\imj{T}_i}\trace{\Gamma}-{a\imj{I}_i}\trace{\Gamma}$,
we have for all $i \ge k$ in $\NN$,
\[
c_i\imj{R,T}(\mu)
\ \in\
 \mu\,\big(\CC_{i-1}[\mu]\otimes\CC^3\big).
\]
\end{enumerate}
\end{theorem}

\begin{remark}
\noindent{\bf 1.}
 There exist choices of $a^{I,R,T}$ so that
 $\cau^\eps$, ${\cal V}^\eps$, and ${\cal W}^\eps$  is an exact solution.
Since the transmission problem is well posed, there is
a uniquely  determined corrector $c^\eps$ smooth and infinitely
small on  both sides so that adding $c^\eps$ yields an exact solution.
Adding $c^\eps$ to the left corresponds to adding the infinitely
small term $c^\eps\,e^{iS^I/\eps}$ to $a^I$ with a similar remark
on the right.

\noindent
{\bf 2.}  Part $(iv)$ of the theorem with $k=0$ generalizes
part $(iii)$ to discontinuous and variable $\sigma(x_1)$.

\noindent
{\bf 3.} The basis elements, $\Phi\imj{R}$ for $a_1\imj{R}$ and $\Phi\imj{I}$ for $a_0\imj{I}$
 are  homogeneous of degree 2 in $\tau,\xi$.
 Doubling $\tau,\xi$ and also $\eps$ leaves the incoming and
 reflected waves unchanged.   Therefore
 $\eps R(\tau,\xi)$ must be equal to $2\eps R(2\tau,2\xi)$.
 This explains why
 $R$ is  homogeneous of degree -1.

\noindent
{\bf 4.} The reflection coefficient vanishes when $\xi'=0$.
Since it is an even function of $\xi$, $\nabla_\xi R=0$ too.
\end{remark}

\begin{proof}
The incoming solution is given.

\noindent $(i)$
Seek the leading amplitudes $a \imj{T}_0$ and $a\imj{R}_0$. We will show that $a \imj{R}_0=0$ so it is actually $a \imj{R}_1$ that is the leading amplitude of the reflected wave. A jump discontinuity in a lower order coefficient does not lead to reflection at leading order. Denote
\[
\cl_T\ :=\ \partial_t+\caA_2\partial_2+\caA_3\partial_3\ ;
\qquad
\cl_1\ :=\ \cl_T+\caA_1\partial_1\ ;
\qquad
\cl\ :=\ \cl_1+\mu C\,.
\]
${\cal T}:=\partial_t+v_2\partial_2+v_3\partial_3$ is the tangential transport operator.
By Theorem \ref{th:stabilitegen}, the amplitudes are polarized, \textit{i.e.} $a\imj{I,R,T}_0=\Pi_\cl a\imj{I,R,T}_0$, and  $a\imj{T}_0$ (resp. $a \imj{R}_0$) satisfies a forward transport equation in $x_1\ge 0$ (resp.  backward in $x_1\le 0)$ with zero initial values in time,
\begin{equation}
\label{eq:transptoutes}
\begin{array}{ll}
(v_1\partial_1+ {\cal T})\,a\imj{I}_0=0,& \quad x_1 \in \RR,\\
(-v_1\partial_1+ {\cal T})\,a\imj{R}_0=0,& \quad  x_1 \in \RR_-,\\
(v_1\partial_1+ {\cal T}+\mu \Pi_{\cl}\,C\,\Pi_{\cl})\,a\imj{T}_0=0,& \quad x_1 \in \RR_+.\\
\end{array}
\end{equation}

Therefore, to determine $a\imj{T}_0$ and $a\imj{R}_0$ everywhere, it suffices to know $a\imj{T}_0(t,0\lmj{+},x')$ and $a \imj{R}_0(t,0\lmj{-},x')$. These values are determined from the transmission condition \eqref{eq:transcond}:
\begin{equation}
\label{eq:transcondordrezero}
\caA_1
(a\imj{I}_0(t,0\lmj{-},x')
\ +\
a\imj{R}_0(t,0\lmj{-},x'))
\ =\
\caA_1 a\imj{T}_0(t,0\lmj{+},x').
\end{equation}
The
matrix $\caA_1$ is singular. It is easy to see that $(\ker\cl\txi \,\oplus\, \ker\cl\txim )\,\cap \,\ker \caA_1={0}$. Therefore,
$\caA_1\ker\cl\txi$ and $\caA_1 \ker\cl\txim$ are complementary subspaces and generate $\mathrm{Range} \, \caA_1$. This proves that
\begin{equation}\label{eq:transcondordrezeromax}
a\imj{R}_0(t,0\lmj{-},x')
\ =\
0,
\qquad
a\imj{T}_0(t,0\lmj{+},x')
\ -\
a\imj{I}_0(t,0\lmj{-},x')
\ =\
0.
\end{equation}
By the transport equation, we conclude that
\begin{equation}\label{eq:reflexpartout}
a\imj{R}_0
\ \equiv\
 0,\qquad \mbox{for }x_1 < 0.
\end{equation}
The reflected zeroth order term vanishes identically when
$x\in\RR^3_-$.
We also deduce from the transport equation \eqref{eq:transptoutes} that
\begin{equation}\label{eq:saut1}
v_1(\partial_1 a\imj{T}_0-\partial_1 a\imj{I}_0)
\ +\
\mu\, \Pi_{L}\,C\,\Pi_{L}a\imj{T}_0
\ =\
0 \quad
\mbox{  on }
\quad
\Gamma.
\end{equation}%

Next
 determine inductively the correctors. For simplicity, throughout the proof we note $\Pi_\cl:=\Pi_\cl\txi$ and $\Pi_\cl\imj{R}:=\Pi_{L}\txim$.
Write the recursion relation \eqref{eq:rec} for $j \ge 1$ for the incident, reflected and transmitted waves. Split the amplitudes as
$$
\begin{array}{rcl}
a\imj{I,T}_j(t,x)
\ &=&\
\Pi_{\cl} a \imj{I,T}_j(t,x)
\ +\
 \ipc \, a \imj{I,T}_{j}(t,x) ,\\
a \imj{R}_j(t,x)\ &=&\
\Pi_{\cl}\imj{R} a \imj{R}_j(t,x)
\ +\
 \ipcr\, a \imj{R}_{j}(t,x).
\end{array}
$$
$\ipc a \imj{T}_{j}(t,x)$ and $\ipcr a \imj{R}_{j}(t,x)$ are determined directly by \eqref{eq:recur1}.
To determine the projection on the kernel, split the transmission condition \eqref{eq:transcond} and insert \eqref{eq:recur1} on the interface to get,
\begin{equation}\label{eq:condtrans}
\begin{split}
\ds
\caA_1(\, &\Pi_{\cl}  a \imj{I}_j(t,0 ,x')
       -\Pi_{\cl} a \imj{T}_j(t,0\lmj{+},x')
       +\Pi_{\cl}\imj{R} a\imj{R}_j(t,0\lmj{-},x'))
\ =\\
&-\,\caA_1(
\ipc a \imj{I}_j(t,0 ,x')
- \ipc a \imj{T}_j(t,0\lmj{+},x')%\\
%&\hspace{3cm}
+\ipcr a \imj{R}_j(t,0\lmj{-},x')
).
\end{split}
\end{equation}
 As for the terms of order 0, this determines  $ \Pi_{\cl} a \imj{T}_j(t,0\lmj{+},x')$ and
$ \Pi_{\cl}\imj{R} a \imj{R}_j(t,0\lmj{-},x')$.
By \eqref{eq:recur2}, the projections are solution of a transport equation, therefore uniquely determined by initial data and the values on the boundary . Borel's theorem allows one to construct
$$
a \imj{I}(t,x,\eps),
\qquad
a \imj{T}(t,x,\eps),
\qquad
{\rm and,}
\qquad
a \imj{R}(t,x,\eps),
$$
so that the transmission condition is exactly satisfied.  With this choice the
approximate solution satisfies the transmission problem with infinitely small
residual.\\

\noindent $(ii)$
Theorem \ref{thm:berengerperfect} implies that the exact
solution in $x_1\le 0$ is equal to $\cau^\eps+O(\eps^\infty)$.
The error of the approximation is $O(\eps^\infty)$ so
the exact solution is equal to $\cau^\eps + \cav^\eps +O(\eps^\infty)$.
Therefore
$\cav^\eps = (\cau^\eps +\cav^\eps)-\cau^\eps=O(\eps^\infty)$ which is the desired conclusion.\\

\noindent $(iii)$  For the smart layer (\ref{eq:newC}, \ref{eq:Lmu}) with $\sig={\bf 1}_{x_1>0}$, compute the first order term by \eqref{eq:recpml} with $j=1$.
First deduce from \eqref{eq:recur1} that
\begin{equation}\label{eq:nopolar1}
\begin{array}{lll}
\ip a_1\imj{I}(t,x)&=\
i Q_{L} L_1(\partial
      _t,\partial
      _x)\,a_{0}\imj{I}(t,x),&\quad x_1\in \RR,\\
\ipr a_1\imj{R}(t,x)&=\ 0,&\quad x_1\in \RR_-,\\
\ip a_1\imj{T}(t,x)&=\
i Q_{L}( L_1(\partial
      _t,\partial
      _x)+\mu \,C)\,a_{0}\imj{T}(t,x),&\quad  x_1\in \RR_+.
\end{array}
\end{equation}
Replace in \eqref{eq:nopolar1} the $x_1$ derivatives using \eqref{eq:transptoutes},
\begin{equation*}
\begin{split}
&\ip a_1\imj{I}(t,x)
=i Q_{L} (L_T+A_1\partial_1)\,
a_{0}\imj{I}(t,x)%,\\
=i Q_{L} (L_T-\frac{1}{v_1}A_1{\cal T})\,
a_{0}\imj{I}(t,x),\\
&(I-\Pi_{L})a_1\imj{T}(t,x)
=i Q_{L} (L_T+A_1\partial_1+\mu\, C)\,
a_{0}\imj{T}(t,x)%,\\
=i Q_{L} \bigl(
L_T-\frac{1}{v_1}A_1{\cal T}
+\mu\,(- \frac{1}{v_1} A_1\Pi_{L}\,C\,\Pi_{L}+  C\,)
\bigr)\,
a_{0}\imj{T}(t,x)\,,
\end{split}
\end{equation*}
to obtain,  with $\ds C_1:=  \ds C -\frac{1}{v_1} A_1\Pi_{L}\,C\,\Pi_{L}
$,
\begin{equation}\label{eq:nopolar2}
(I-\Pi_{L})a_1\imj{T}(t,0\lmj{+},x')
\ -\
(I-\Pi_{L})a_1\imj{I}(t,0,x')
\ =\
i \mu Q_{L} C_1
a_{0}\imj{I}(t,0,x').
\end{equation}
Using  $(I-\Pi_{L})a_1\imj{R}=0$
in the transmission condition yields,
\begin{equation}\label{eq:trans1}
A_1\big(
\Pi_{L}\imj{R}a_1\imj{R}\ +\ \Pi_{L}a_1\imj{I}-\Pi_{L}a_1\imj{T}
\big)
\ =\
i \mu
A_1 Q_{L} C_1 a_0\imj{I}\,.
\end{equation}

The eigenvalues of $A_1$ are $0$ and $\pm 1$, with associated
 orthonormal set of eigenvectors
$\Phi_0=e_1$ and $\Phi_\pm=\big(0,1,\pm i\big)/\sqrt 2$.
The projection operators on the positive and negative eigenspaces are  $\pi_\pm(A_1)=\Phi_\pm \Phi_\pm^*$, and $C=\Phi_+\Phi_+^*+\nu \Phi_-\Phi_-^*$. The kernel of $\ltx$  is
one-dimensional, it is spanned by
\begin{equation}\label{eq:phi}
%\begin{array}{lcl}
\Phi\txi
\ =\
\xi -\ds\frac{\tau^2}{\xi_1}e_1
+i\ds\frac{\tau}{\xi_1}\xi\wedge e_1
\ =\
\left(\xi_1\,-\ds\frac{\tau^2}{\xi_1},
\ds i\frac{\tau}{\xi_1}\xi_3\,+\,\xi_2,
\ds -i\frac{\tau}{\xi_1}\xi_2\,+\,\xi_3
\right), \end{equation}
and the projection on $\kl$ is
$\Pi_{L}=\frac{\Phi \Phi^*}{\Phi ^*\Phi}$.  Compute
\[
\begin{array}{ll}
\Pi_{L} C \Pi_{L}&=\ds \frac{\Phi \Phi^*}{\Phi ^*\Phi}(\Phi_+\Phi_+^*+\nu \Phi_-\Phi_-^*)\frac{\Phi \Phi^*}{\Phi ^*\Phi}
\\
&=\ds \frac{1}{(\Phi ^*\Phi)^2}
\Phi (\Phi^*\Phi_+)(\Phi_+^*\Phi )\Phi^*+\nu
\Phi (\Phi^*\Phi_+)(\Phi_+^*\Phi) \Phi^*\\
&=\ds \frac{|\Phi^* \Phi_+|^2+\nu |\Phi^* \Phi_-|^2}{\Phi ^*\Phi}\  \Pi_{L}.
\end{array}
\]
Define
\begin{equation}\label{eq:pcp}
\gamma \ :=\
 \ds \frac{|\Phi^* \Phi_+|^2+\nu |\Phi^* \Phi_-|^2}{\Phi ^*\Phi},
 \qquad
 {\rm so},
 \qquad
\Pi_{L} C \Pi_{L} = \gamma \Pi_{L}\,.
\end{equation}
Since $a_0\imj{I}$ is polarized,
\[
\ds C_1a_0\imj{I}=  \ds (C -\frac{\gamma}{v_1} A_1\Pi_{L})\,a_0\imj{I}
=
(C -\frac{\gamma}{v_1} A_1)\,a_0\imj{I}
:=
\widetilde{C}_1\,a_0\imj{I}.
\]
To compute the righthand side of \eqref{eq:trans1}, use
\[
C=
\begin{pmatrix}
0 \ &  0  & 0 \\[2mm]
0 & \ds \frac{\nu+1}{2}& \ds i\,\frac{\nu-1}{2}\\[3mm]
0 & \ds -i\, \frac{\nu-1}{2}& \ds \frac{\nu+1}{2}\\
\end{pmatrix},
\gamma \ =\
\ds\frac{(\tau-\xi_1)^2+\nu(\tau+\xi_1)^2}{4\tau^2}
\ =\ \ds\frac{1}{4}((1+v_1)^2+\nu(1-v_1)^2)
\,,
\]
%
%\[
%C=
%\begin{pmatrix}
%0 \ &  0  & 0 \\[2mm]
%0 & \ds \frac{\nu+1}{2}& \ds i\,\frac{\nu-1}{2}\\[3mm]
%0 & \ds -i\, \frac{\nu-1}{2}& \ds \frac{\nu+1}{2}\\
%\end{pmatrix},\qquad
%\ds \widetilde{C}_1\ =\
%\begin{pmatrix}
%0 \ &  0  & 0 \\[2mm]
%0 & \ds \frac{\nu+1}{2}& \ds i\,\frac{\nu-1}{2}\\
%0 & \ds -i\, \frac{\nu-1}{2}& \ds \frac{\nu+1}{2}\\
%\end{pmatrix}
%\ -\
%\ds\frac{\gamma}{v_1}
%\begin{pmatrix}
%0 \ &  0  & 0 \\
%0 & 0& \ds -i\\
%0 & \ds i\, & 0\\
%\end{pmatrix}\,.
%\]
%Simplify
to find
\begin{equation*}
 \widetilde{C}_1 =
\frac{\nu+1}{2}
\begin{pmatrix}
0 \quad &  0  & 0 \\
0 \quad & 1&\ds i\,\frac{v_1^2+1}{2v_1}\\
0 \quad &-\ds i\,\frac{v_1^2+1}{2v_1} & 1
\end{pmatrix}\,.
\end{equation*}
Write $a_{0}\imj{I}=\alpha_0\imj{I}\Phi$, and compute
\[
\begin{split}
\widetilde{C}_1\,\Phi &=
\ds(1+\nu)\
\frac{\xi_1^2-\tau^2}{4\tau\xi_1^2}\ \Psi,
\qquad
\Psi=\begin{pmatrix}
0  \\
\ds \xi_2\tau-i\xi_1\xi_3,
\\
 \xi_3\tau+i\xi_2\xi_3,
 \end{pmatrix}
\end{split}
\, ,
\]
to find a new version of \eqref{eq:trans1},
\begin{equation}\label{eq:trans3}
A_1(\Pi_{L}\imj{R}a_1\imj{R}+ \Pi_{L}a_1\imj{I}-\Pi_{L}a_1\imj{T})=
i \mu \ds(1+\nu)\
\frac{\xi_1^2-\tau^2}{4\tau\xi_1^2}\ \alpha_0\imj{I}\
A_1
Q_L\, \Psi\,.
\end{equation}

Next compute $Q_L\, \Psi$. First compute a basis of eigenvectors for $L\txi$.
$\Phi_2$ is such that $\ltx\Phi_2=\tau \Phi_2$, $\Phi_3$ is such that $\ltx\Phi_3=2\tau \Phi_3$.   Choose
\[
\Phi_2=\xi,\qquad
\Phi_3=\Phi(-\tau,\xi) .
\]

Note that
\[
\Psi=\tau \xi -\xi_1(\tau e_1+i \xi \wedge e_1),
\]
and
\[
\Phi=\xi -\ds \frac{\tau}{\xi_1}(\tau e_1-i \xi \wedge e_1),\quad
\Phi_3=\xi -\ds \frac{\tau}{\xi_1}(\tau e_1+i \xi \wedge e_1),
\]
which gives
\[
(\tau e_1+i \xi \wedge e_1)
\ =\
\ds \frac{\xi_1}{\tau}(\xi -\Phi_3),
\]
and
\[
\Psi
\ =\
\tau \xi -\frac{\xi_1^2}{\tau}(\xi -\Phi_2)
\ =\
\frac{\tau ^2-\xi_1^2}{\tau}\xi
 + \frac{\xi_1^2}{\tau}\Phi_3.
\]

Since $Q_L$ is the left inverse of $L$, we have $Q_L\xi=\frac{1}{\tau}\xi$, and $Q_L\Phi_3=\frac{1}{2\tau}\Phi_3$, which gives

\[
Q_L \Psi
\ =\
\frac{\tau ^2-\xi_1^2}{\tau^2}\xi
 \ +\
  \frac{\xi_1^2}{2\tau^2}\Phi_3
  \,.
\]

Write the coefficients on $\Gamma$ as
\[
\Pi_{L} a\imj{I,T}_1
\ =\
\alpha\imj{I,T}_1\Phi,\qquad
\Pi_{L}\imj{R} a\imj{R}_1
\ =\
\alpha\imj{R}_1\Phi\imj{R} \,,
\]
and inject  into the transmission condition to obtain
\[
\alpha\imj{R}_1 A_1\Phi +
(\alpha\imj{I}_1-\alpha\imj{T}_1)A_1\Phi\imj{R}=
i \mu \ds(1+\nu)\
\frac{\xi_1^2-\tau^2}{4\tau\xi_1^2}\ \alpha_0\imj{I}\
A_1
(\frac{\tau ^2-\xi_1^2}{\tau^2}\xi
 + \frac{\xi_1^2}{2\tau^2}\Phi_3).
\]
Since the kernel of $A_1$ is $e_1$,  $A_1$ in the preceding identity
may be replaced by the projection on $(e_2,e_3)$. The projection of $\Phi$ is
$\phi = \xi'-i\frac{\tau}{\xi_1} \xi\wedge e_1$, and
note that $\Phi_3$ and $\Phi\txim$ have the same projection, which is $\phi_3 = \xi'+i\frac{\tau}{\xi_1} \xi\wedge e_1$. Write

\[
\begin{array}{lcl}
\alpha\imj{R}_1\phi_3 +
(\alpha\imj{I}_1-\alpha\imj{T}_1)\phi
&=&
i \mu \ds(1+\nu)\
\frac{\xi_1^2-\tau^2}{4\tau\xi_1^2}\ \alpha_0\imj{I}\
(\frac{\tau ^2-\xi_1^2}{\tau^2}\xi'
 + \frac{\xi_1^2}{2\tau^2}\phi_3)\\
&=&
i \mu \ds(1+\nu)\
\frac{\xi_1^2-\tau^2}{4\tau\xi_1^2}\ \alpha_0\imj{I}\
(\frac{\tau ^2-\xi_1^2}{2\tau^2}(\phi+\phi_3)
 + \frac{\xi_1^2}{2\tau^2}\phi_3)\\
 &=&
i \mu \ds(1+\nu)\
\frac{\xi_1^2-\tau^2}{4\tau\xi_1^2}\ \alpha_0\imj{I}\
( \frac{\tau ^2-\xi_1^2}{2\tau^2}\phi+ \frac{1}{2}\phi_3).
\end{array}
\]
The solutions are parameterized by
$\alpha_0\imj{I}$,
\[
\alpha\imj{R}_1
\ =\
i \mu \ds(1+\nu)\
\frac{\xi_1^2-\tau^2}{8\tau\xi_1^2}\ \alpha_0\imj{I},\qquad
\alpha\imj{I}_1-\alpha\imj{T}_1
\ =\
-i \mu \ds(1+\nu)\
\frac{(\xi_1^2-\tau^2)^2}{8\tau^3\xi_1^2}\ \alpha_0\imj{I}\,.
\]

Now use the results in  $(i)$ and prove \eqref{eq:linear} by
induction on $i$. Equation
\eqref{eq:reflexpartout} asserts  that $a\imj{R}_0=0$ for $x_1 < 0$.
Equations \eqref{eq:transptoutes},  and \eqref{eq:transcondordrezero},
imply that at $x_1=0$, $a\imj{I,T}_0=\Pi_{L}a\imj{I,T}_0$, and
\[
v_1(\partial_1a\imj{T}_0-\partial_1 a\imj{I}_0)+
{\cal T}(a\imj{T}_0- a\imj{I}_0)
+\mu\gamma(a\imj{T}_0- a\imj{I}_0)
+\mu\gamma a\imj{I}_0=0,
\qquad
 x_1 \ge 0.
\]
Differentiation in $x_1$ several times yields
$ \partial_1^ja\imj{T}_0-\partial_1^j a\imj{I}_0 \in
\mu\,\CC_{j-1}[\mu]\otimes\CC^3$, giving the results for $i=0$ and $j \ge 1$.

Assuming the inductive hypothesis is true for $i$ we prove it for
$i+1$.
Write  \eqref{eq:recpml} in the form
\[
\begin{array}{lcl}
(I-\Pi_{L}\imj{R})\,a_{i+1}\imj{R}(t,x)&=&
i Q_{L}\imj{R}
(L_T(\partial_t,\partial_{x'})+A_1\partial_1)
\,a_{i}\imj{I}(t,x),\\
(I-\Pi_{L})\,a_{i+1}\imj{I}(t,x)&=&
i Q_{L}(L_T(\partial_t,\partial_{x'})+A_1\partial_1)
\,a_{i}\imj{I}(t,x),\\
(I-\Pi_{L})\,a_{i+1}\imj{T}(t,x)&=&
i Q_{L}(L_T(\partial_t,\partial_{x'})+A_1\partial_1
+\mu C)\,a_{i}\imj{T}(t,x).
\end{array}
\]
By induction, $\partial_1^j(I-\Pi_{L}\imj{R})\,a_{i}\imj{R}(t,x) \in
\mu\,\CC_{i+j-1}[\sig]\otimes\CC^3
$ on the interface.  Write for $x_1 \ge 0$,
\begin{equation}\label{eq:nopolardec}
%\begin{array}{lcl}
 (I-\Pi_{L})\,(a_{i+1}\imj{T}-a_{i+1}\imj{I})(t,x)
% &=&
=
i Q_{L}\,(L_T+A_1\partial_1+\mu C)
\,(a_{i}\imj{T}-a_{i}\imj{I})(t,x)+
i \mu \,Q_{L}\,C\,a_{i}\imj{I}(t,x).
%\end{array}
\end{equation}
The inductive hypothesis,
shows that  on $\Gamma$,
$\partial_1^j(L_T+A_1\partial_1+\mu C)
\,(a_{i}\imj{T}-a_{i}\imj{I}) \in
\mu\,\CC_{i+j}[\mu]\otimes\CC^3
$
and $\partial_1^j(\mu C\,a_{i}\imj{I}) \in
\mu\,\CC_{0}[\sig]\otimes\CC^3$.  The result follows for $(I-\Pi_{L})\,(a_{i+1}\imj{T}-a_{i+1}\imj{I})$, and for
$(I-\Pi_{L}\imj{R})\,a_{i+1}\imj{R}$ in the same way. The transmission condition extends the assertion to the other parts $\Pi_{L}\,(a_{i+1}\imj{T}-a_{i+1}\imj{I})$ and $\Pi_{L}\imj{R}\,a_{i+1}\imj{R}$.

\textit{\textrm{(iv)}}
Here $\sigma$ vanishes to order $k$ at $x_1=0$.
\eqref{eq:transcondordrezeromax}
and \eqref{eq:reflexpartout} are still valid, and  the transport equations \eqref{eq:transptoutes} implies  on the interface $\Gamma$  that
\begin{equation}\label{eq:smooth0}
\begin{array}{l}
%a_0\imj{I,R,T} \equiv \Pi_{L} a_0\imj{I,R,T},  \qquad
%a_0\imj{R}\equiv 0,\qquad \mbox{ everywhere},  \\
\partial_1^j \Pi_{L}a_0\imj{T} - \partial_1^j \Pi_{L}a_0\imj{I}
\ =\ 0, \qquad j=0,\cdots,k,\quad\\
\partial_1^{k+1} \Pi_{L}a_0\imj{T} - \partial_1^{k+1} \Pi_{L}a_0\imj{I}\ =\
-\mu\sig^{(k)}(0) \ds \frac{\gamma}{v_1}\Pi_{L}a_0\imj{I} .
\end{array}
\end{equation}

From \eqref{eq:nopolardec} for $i=0$, \eqref{eq:nopolar1}  and \eqref{eq:smooth0}, derive
\begin{equation}\label{eq:smooth1ip}
\begin{array}{l}
a_1\imj{R} \equiv \Pi_{L} a_1\imj{R} \qquad \mbox{ everywhere,} \\
\partial_1^j \ip a_1\imj{T} - \partial_1^j \ip a_1\imj{I}=0, \qquad j=0,\cdots,k-1,\quad\mbox{ on } \Gamma, \\
\partial_1^{k} \ip a_1\imj{T} - \partial_1^k \ip a_1\imj{I}
=i\mu\sig^{(k)}(0) \ds Q_LC_1\Pi_{L} a_0\imj{I}\qquad\mbox{ on } \Gamma .
\end{array}
\end{equation}
Using the transmission conditions to obtain on the interface $\Gamma$,
\[
\Pi_{L}\imj{R} a_1\imj{R} = 0, \qquad \Pi_{L} a_1\imj{I}=\Pi_{L} a_1\imj{T} .
\]

Insert into the transport equations \eqref{eq:recur2}
to find $ a_1\imj{R} = 0 $ in $\RR_-$, and
\begin{equation}\label{eq:smooth1}
\begin{array}{l}
 \partial_1^j \Pi_{L} a_1\imj{T} - \partial_1^j \Pi_{L} a_1\imj{I}=0, \ j=0,\cdots,k-1,\quad\mbox{ on } \Gamma,\\
\partial_1^{k} \Pi_{L} a_1\imj{T} - \partial_1^k \Pi_{L} a_1\imj{I}
=-\frac{1}{v_1}\Pi_{L} A_1
(\partial_1^{k} \ip a_1\imj{T} - \partial_1^k \ip a_1\imj{I})
=-i\frac{\mu}{v_1}\sig^{(k)}(0) \ds \Pi_{L}A_1Q_LC_1\Pi_{L} a_0\imj{I}\quad\mbox{ on } \Gamma .
\end{array}
\end{equation}
Recover
\begin{equation}\label{eq:smooth1bis}
\begin{array}{l}
 \partial_1^j a_1\imj{T} - \partial_1^j a_1\imj{I}=0, \ j=0,\cdots,k-1,\qquad\mbox{ on } \Gamma,\\
\partial_1^{k}  a_1\imj{T} - \partial_1^k a_1\imj{I}
=i\mu\sig^{(k)}(0)(I-\frac{1}{v_1}\Pi_{L}A_1) \Pi_{L} A_1
 \ds Q_LC_1\Pi_{L} a_0\imj{I}
 \qquad
 \mbox{ on } \Gamma .
\end{array}
\end{equation}
Now proceed iteratively, to see that, for $i < k+1$,
\begin{equation}\label{eq:smoothi}
\begin{array}{l}
 a_i\imj{R} \equiv \Pi_{L} a_i\imj{R} \equiv 0\qquad \mbox{ everywhere}, \\
\partial_1^j a_i\imj{T} - \partial_1^j a_1\imj{I}=0, \qquad j=0,\cdots,k-i,\quad\mbox{ on } \Gamma,\\
\partial_1^{k-i+1}  \ip a_i\imj{T} - \partial_1^{k-i+1}\ip a_i\imj{I}
=iQ_LA_1
(\partial_1^{k-i+2}  a_{i-1}\imj{T} - \partial_1^{k-i+2} a_{i-1}\imj{I}) \quad\mbox{ on } \Gamma\, ,\\
\partial_1^{k-i+1} \Pi_{L} a_i\imj{T} - \partial_1^{k-i+1}\Pi_{L} a_i
\imj{I}
=
-\ds\frac{1}{v_1}\Pi_{L} A_1
(\partial_1^{k-i+1}  \ip a_i\imj{T} - \partial_1^{k-i+1}\ip a_i\imj{I}
)\,.
\end{array}
\end{equation}

Denote by  $s_i$ the value of $ \partial_1^{k-i+1} a_i\imj{T} - \partial_1^{k-i+1} a_i\imj{I}$ on $\Gamma$.
Equation  \eqref{eq:smoothi}
yields the recursion relation
\[
%\begin{array}{l}
s_j=i(I-\ds\frac{1}{v_1}\Pi_{L} A_1)Q_LA_1s_{j-1}\,,
\quad %\\
s_1=i\mu\sig^{(k)}(0)(I-\ds\frac{1}{v_1}\Pi_{L} A_1)Q_LC_1\Pi_{L} a_0\imj{I}\,,
% \end{array}
\]
which can be solved as
\[
 s_{k+1}=i\mu\,\sig^{(k)}(0)M\Pi_{L} a_0\imj{I}, \qquad
 {\rm with}\qquad
 M:=(i(I-\ds\frac{1}{v_1}\Pi_{L} A_1)Q_LA_1)^{k} Q_LC_1.
\]
The first nonzero reflected term is therefore $a_{k+1}\imj{R}=\Pi_{L}\imj{R} a_{k+1}\imj{R}$, and using the transmission condition yields
\[
A_1(\Pi_{L}a_{k+1}\imj{R}
+\Pi_{L}a_{k+1}\imj{I}-\Pi_{L}a_{k+1}\imj{T}
+s_{k+1})\ =\ 0\,.
\]
The incoming amplitude on $\Gamma$ is $ a_0\imj{I} = \alpha(t,x') \Phi$, the leading reflection is $a_{k+1}\imj{R} =\alpha_{k+1}\imj{R}\Phi\imj{R}$ and the leading transmission is $a_{k+1}\imj{T} =\alpha_{k+1}\imj{R}\Phi\imj{T}$.
 Using again the notation $\Phi'$ to denote the projection of a vector $\Phi$ on  $Vec(e_2,e_3)$, this linear system is solved as
 \[
 \alpha_{k+1}\imj{R}= -\ds\frac{\Phi'\wedge s_{k+1}'}{\Phi'\wedge (\Phi\imj{R})'}
,\quad \alpha_{k+1}\imj{T}-\alpha_{k+1}\imj{I}= \ds\frac{(\Phi\imj{R})'\wedge s_{k+1}'}{\Phi'\wedge (\Phi\imj{R})'}
.
 \]
 \[
 \alpha_{k+1}\imj{R}= -i\mu\sig^{(k)}(0)\alpha\ds\frac{\Phi'\wedge (M\Phi)'}{\Phi'\wedge (\Phi\imj{R})'}
,\quad \alpha_{k+1}\imj{T}-\alpha_{k+1}\imj{I}= i\mu\sig^{(k)}(0)\alpha \ds\frac{(\Phi\imj{R})'\wedge (M\Phi)'}{\Phi'\wedge (\Phi\imj{R})'}
.
 \]
The proof of the linearity follows the same path as in $(iii)$.
\end{proof}

%
%
%We compute
%a geometric form of the reflection coefficient.
% The {\it angle of incidence} $0\le \theta<\pi/2$ is
%defined by
%$
%\cos\theta = v_1
%$.
%The reflection coefficient is equal to
%$$
%i\
%\sig\
% \ds\frac{1+\nu}{16\tau}
% \
%\tan^2\,\theta
%\,.
%$$
%

\subsection{Harmoniously matched layers}\label{subsec:hml}

Based on Theorem \ref{th:transmission} we construct an extrapolation method for symmetric hyperbolic operators with smart layers which eliminates the leading order reflection. The resulting method has desirable stability properties and is nearly as good as B\'erenger's algorithm for the Maxwell equations where his method is at its best.   We think that the new method provides a good alternative in situations where B\'erenger's  method is not so effective.

Consider the computational domain
$x_1\le b_1$.  The domain of interest is
the interval
$x_1\le a_1<b_1$.
The absorbing layer is located in $a_1\le x_1\le b_1$.
 The differential operator
in the computational domain is  symmetric hyperbolic $L$
with smart layer
\[
L\,U\, +\sigma_1(x_1)
\big(\pi_{+}(A_1)+ \nu\, \pi_{-}(A_1)\big)
\,U
\ =\
0\,,
\qquad \sigma_1\ge 0\,,
\qquad
{\rm supp}\,\sigma_1\subset \{x_1\ge a_1\}\,.
\]
At the outer boundary $x_1=b_1$ of the absorbing layer impose
the simplest weakly reflecting boundary condition
\[
 \pi_{-}(A_1)\,U\ =\ 0
 \qquad {\rm when}\quad x_1=b_1\,.
\]
This is a well posed problem provided that $A_1$ has constant
rank on $x_1=b_1$.
When $L=L_1(\partial)$ has constant coefficients
it
generates a contraction group in $L^2(\{x_1\le b_1\})$.

The {\bf hamoniously matched layer algorithms}
compute a smart layer with coefficient $\sig_1$ and also
 with coefficient $2\sig_1$.
In view of Theorem \ref{th:transmission},
 subtracting the second from twice the first,
 $\,2\,U(\sigma_1)-U(2\sigma_1)$,
yields a field with  one more vanishing term in the reflected wave
at the interface $x_1=a_1$.
This extrapolation removes the leading reflection.

The harmonious matched layer algorithms in a rectangular
domain $\cal R$ perform the same extrapolation with absorptions
in  all directions.  With
\[
LU\,
\  +\
\sum_{j=1}^d
\sigma_j(x_j) (\pi_{+}(A_j)+ \nu\, \pi_{-}(A_j))U=0,
\qquad \sigma_j\ge 0,
\qquad
{\rm supp}\,\sigma_j\subset \{|x_j|\ge a_j\}\,.
\]
with
\[
 \pi_{\mp}(A_j)\,U\ =\ 0
 \qquad {\rm when} \quad x_j=\pm b_j\,.
\]
This initial boundary value problem
on a rectangle has weak solutions. \footnote{This can be proved by penalisation.  Denote
by $\Omega$ the rectangular computational domain.
Add
$\Lambda \,{\bf 1}_{\RR^d\setminus \Omega}$
 to $L$
and solve on $\RR^{1+d}_{t,x}$.  The limit
as $\Lambda\to \infty$ provides a solution
in
$L^\infty\big([0,T]\,:\, L^2({\cal R})\big)$
\cite{Bardos:1982:MPB}.
}
  When $L=L_1(\partial)$,
the
$L^2( {\cal R} )$ norm is nonincreasing in time.
The extrapolation is
$
2\,U(\sig_1,\cdots,\sig_d)
-
U(2\sig_1,\cdots,2\sig_d)
$.
\vskip.2cm
\noindent
{\bf Open Problem.}  For discontinuous $\sigma_j$, the uniqueness of solutions
to the initial boundary value problem on the rectangular
computational domain
is not known because of
 the discontinuity of the boundary  space
${\rm ker}\,A_j$ at the corner.  Solutions
are typically discontinuous.
 Uniqueness of strong solutions
  and existence of weak solutions is
proved by the energy method.
{\sl  We do not know how to prove uniqueness of solutions with regularity
not exceeding that of  solutions known to  exist.
}  Similar problems plague virtually all methods on
rectangular domains with absorbing boundary conditions
imposed on the computational domain with corners.
The present problem is one of  the simplest
of its kind.  The fact that algorithms designed
to compute solutions encounter no difficulties
is reason for optimism.

\vskip.2cm

\subsection{Numerical experiments}\label{subsec:numerics}

Simulations are performed
 for  the 2-D transverse electric Maxwell system in the $(x,y)$ coordinates,
\begin{equation} \label{maxwell2D}
\begin{array}{l}
\pd t E_x - \pd y H_{z} =0\,, \\
\pd t E_y +\pd x H_{z}=0\,,\\
\pd t H_{z}+\pd x E_y -\pd y E_x=0\,,
\end{array}
\end{equation}
in a rectangle, with boundary conditions $n\wedge E=0$ on the west, north and south boundaries. The layer will be imposed on the east boundary.
Maxwell B\'erenger is given by

\begin{equation} \label{PMLmaxwellber}
\begin{array}{l}
\pd t E_x - \pd y H_{z} =0\,, \\
\pd t E_y +\pd x H_{z}+\sigma(x) E_y=0\,,\\
\pd t H_{zx}+\pd x E_y +\sigma(x) H_{zx}=0\,, \\
\pd t H_{zy}-\pd y E_x=0\,,\\
H_{z}=H_{zx}+H_{zy}\,.
\end{array}
\end{equation}
For the computation, these equations are used in the whole rectangle (see the discussion in the introduction), with $\sig=0$ outside the layer. The boundary conditions are
\begin{equation}\label{eq:abc}
E_y=H_{z} \mbox{ and }
 E_x=0\,\mbox{ on the east }, \quad
n\wedge E=0 \mbox{ on the other boundaries}.
\end{equation}
Since
$\Pi_+(A_1)= \frac{E_y+H_z}{2} (0,\ 1,\ 1)$ and $\Pi_-(A_1)= \frac{E_y-H_z}{2} (0,\ 1,\ -1)$, the smart layers are:
\begin{equation} \label{PMLmaxwellhpr}
\begin{array}{l}
\pd t E_x - \pd y H_z =0\,, \\
\pd t E_y +\pd x H_z+\frac{\sigma(x)}{2}(E_y+H_z+\nu(E_y-H_z))=0\,,\\
\pd t H_z-\pd y E_x+\pd x E_y +\frac{\sigma(x)}{2}(E_y+H_z-\nu(E_y-H_z))=0\,.\\
\end{array}
\end{equation}
The boundary conditions \eqref{eq:abc} are imposed.

The Yee scheme for Maxwell is
\begin{subequations}\label{eq:yeemax}
\renewcommand{\theequation}{\theparentequation-\alph{equation}}
\begin{align}
&\ds\frac{(E_x)_{i+\fr,j}^{n}-(E_x)_{i+\fr,j}^{n-1}}{\Delta t}
-\frac{(H_z)_{i+\fr,j+\fr}^{n-\fr}-(H_z)_{i+\fr,j-\fr}^{n-\fr}}{\Delta y}=0
% \label{eq:yeemax_a}
\nonumber
\\
&\ds\frac{(E_y)_{i,j+\fr}^{n}-(E_y)_{i,j+\fr}^{n-1}}{\Delta t}+
\frac{(H_z)_{i+\fr,j+\fr}^{n-\fr}-(H_z)_{i-\fr,j+\fr}^{n-\fr}}{\Delta x}=0
% \label{eq:yeemax_b}
\nonumber
\\
&\ds\frac{(H_{z})_{i+\fr,j+\fr}^{n+\fr}-(H_{z})_{i+\fr,j+\fr}^{n-\fr}}{\Delta t}
+\frac{(E_y)_{i+1,j+\fr}^{n}-(E_y)_{i,j+\fr}^{n}}{\Delta x}
-\frac{(E_x)_{i+\fr,j+1}^{n}-(E_x)_{i+\fr,j}^{n}}{\Delta y}=0.
% \label{eq:yeemax_c}
\nonumber
\end{align}
\end{subequations}
The Yee scheme for Maxwell B\'erenger using  the notations $\sig_i=\sigma(x_i)$ and
$\sig_{i+\fr}=\sigma(x_{i+\fr})$ is,
\begin{subequations}\label{eq:yeeber}
\begin{align}
&\ds\frac{(E_x)_{i+\fr,j}^{n}-(E_x)_{i+\fr,j}^{n-1}}{\Delta t}
\ -\
\frac{(H_z)_{i+\fr,j+\fr}^{n-\fr}-(H_z)_{i+\fr,j-\fr}^{n-\fr}}{\Delta y}\ =\
0\,, %\label{eq:yeeber_c}
\nonumber
\\
&\ds\frac{(E_y)_{i,j+\fr}^{n}-(E_y)_{i,j+\fr}^{n-1}}{\Delta t}
\ +\
\frac{(H_z)_{i+\fr,j+\fr}^{n-\fr}-(H_z)_{i-\fr,j+\fr}^{n-\fr}}{\Delta x}
\ +\ \sig_{i}\ \frac{(E_y)_{i,j+\fr}^{n}+(E_y)_{i,j+\fr}^{n-1}}{2}\ =\
0\,, % \label{eq:yeeber_d}
\nonumber
\\
&\ds\frac{(H_{zx})_{i+\fr,j+\fr}^{n+\fr}-(H_{zx})_{i+\fr,j+\fr}^{n-\fr}}{\Delta t}
\ +\
\frac{(E_y)_{i+1,j+\fr}^{n}-(E_y)_{i,j+\fr}^{n}}{\Delta x}
\ +\
\sig_{i+\fr} \frac{(H_{zx})_{i+\fr,j+\fr}^{n+\fr}+(H_{zx})_{i+\fr,j+\fr}^{n-\fr}}{2}
\ =\ 0\,,
%\label{eq:yeeber_a}
\nonumber
\\
&\ds\frac{(H_{zy})_{i+\fr,j+\fr}^{n+\fr}-(H_{zy})_{i+\fr,j+\fr}^{n-\fr}}{\Delta t}
\ -\
\frac{(E_x)_{i+\fr,j+1}^{n}-(E_x)_{i+\fr,j}^{n}}{\Delta y}
\ =\ 0\,,% \label{eq:yeeber_b}
\nonumber\\
&(H_z)_{i+\fr,j+\fr}^{n+\fr}
\ =\
(H_{zx})_{i+\fr,j+\fr}^{n+\fr}
                            \ +\
                            (H_{zy})_{i+\fr,j+\fr}^{n+\fr}\,.
\nonumber
\end{align}
\end{subequations}

%The computation of $H_{zx}$ and $H_{zy}$ can be explicited as
%\begin{subequations}\label{eq:yeeberbis}
%\begin{align}
%&\ds (1+\frac{\Delta\,t}{2}\sigma(x_{i+\fr}))
%(H_{zx})_{i+\fr,j+\fr}^{n+\fr}=
%(1-\frac{\Delta\,t}{2}\sigma(x_{i+\fr}))(H_{zx})_{i+\fr,j+\fr}^{n-\fr}
%-\frac{\Delta\,t}{\Delta x}((E_y)_{i+1,j+\fr}^{n}-(E_y)_{i,j+\fr}^{n})
%\nonumber
%\\
%&
%\ds
%(H_{zy})_{i+\fr,j+\fr}^{n+\fr}=
%(H_{zy})_{i+\fr,j+\fr}^{n-\fr}
%+\frac{\Delta\,t}{\Delta y}((E_x)_{i+\fr,j+1}^{n}-(E_x)_{i+\fr,j}^{n})\nonumber
%\\
%&\ds (E_x)_{i+\fr,j}^{n+1}=
%(E_x)_{i+\fr,j}^{n}
%+\frac{\Delta\,t}{\Delta y}((H_z)_{i+\fr,j+\fr}^{n+\fr}-(H_z)_{i+\fr,j-\fr}^{n+\fr})\nonumber
%\\
%&\ds (1+\frac{\Delta\,t}{2}\sigma(x_{i}))
%(E_y)_{i,j+\fr}^{n+1}=
%(1-\frac{\Delta\,t}{2}\sigma(x_{i}))(E_y)_{i,j+\fr}^{n}
%-\frac{\Delta\,t}{\Delta x}
%((H_z)_{i+\fr,j+\fr}^{n+\fr}-(H_z)_{i-\fr,j+\fr}^{n+\fr})
%\nonumber
%\\
%\end{align}
%\end{subequations}
The Yee scheme for the smart layers is
\begin{subequations}\label{eq:yeesmart}
\begin{align}
&\ds\frac{(E_x)_{i+\fr,j}^{n}-(E_x)_{i+\fr,j}^{n-1}}{\Delta t}
\ -\
\frac{(H_z)_{i+\fr,j+\fr}^{n-\fr}-(H_z)_{i+\fr,j-\fr}^{n-\fr}}{\Delta y}\ =\ 0\,,
\label{eq:yeesmart_a}
\\
&
\ds\frac{(E_y)_{i,j+\fr}^{n}-(E_y)_{i,j+\fr}^{n-1}}{\Delta t}\ +\
\frac{(H_z)_{i+\fr,j+\fr}^{n-\fr}-(H_z)_{i-\fr,j+\fr}^{n-\fr}}{\Delta x}
\nonumber\\
&
%\hspace{10mm}
\ds +\  \frac{(1+\nu)\sigma_{i}}{2}
%\Biggl[
\frac{(E_y)_{i,j+\fr}^{n}+(E_y)_{i,j+\fr}^{n-1}}{2}
%\Biggr]
%\nonumber\\
%&
%\hspace{15mm}
\ds \ +\
\frac{(1-\nu)}{2}
\frac{\sigma_{i+\fr}(H_{z})_{i+\fr,j+\fr}^{n-\fr}
      +\sigma_{i-\fr}(H_{z})_{i-\fr,j+\fr}^{n-\fr}}{2}
\ =\ 0\,,
\label{eq:yeesmart_b}
\\
&\ds\frac{(H_{z})_{i+\fr,j+\fr}^{n+\fr}-(H_{z})_{i+\fr,j+\fr}^{n-\fr}}{\Delta t}
\ +\
\frac{(E_y)_{i+1,j+\fr}^{n}-(E_y)_{i,j+\fr}^{n}}{\Delta x}
\ -\
\frac{(E_x)_{i+\fr,j+1}^{n}-(E_x)_{i+\fr,j}^{n}}{\Delta y}\nonumber\\
&
%\hspace{10mm}
\ds \ +\
\frac{(1-\nu)}{2}\frac{\sigma_{i+1}(E_y)_{i+1,j+\fr}^{n}+\sigma_i(E_y)_{i,j+\fr}^{n}}{2}
%\nonumber\\
%&\hspace{15mm}
\ds \ + \ \frac{(1+\nu)\sigma_{i+\fr}}{2}\
%\Biggl[
\frac{(H_{z})_{i+\fr,j+\fr}^{n+\fr}+(H_{z})_{i+\fr,j+\fr}^{n-\fr}}{2}
%\Biggr]
\ =\ 0\,.
\label{eq:yeesmart_c}
\end{align}
\end{subequations}
%
%The computation  can be explicited as
%\[
%\begin{cases}
%\ds (E_x)_{i+\fr,j}^{n+1}=(E_x)_{i+\fr,j}^{n}
%+ \frac{\Delta\,t}{\Delta\,y}((H_z)_{i+\fr,j+\fr}^{n+\fr}-(H_z)_{i+\fr,j-\fr}^{n+\fr}) \\
%%
%\begin{split}
%&\ds (1+\frac{\Delta\,t}{4}(1+\nu)\sigma_i)(E_y)_{i,j+\fr}^{n+1}
%=
%(1-\frac{\Delta\,t}{4}(1+\nu)\sigma_i)(E_y)_{i,j+\fr}^{n}
%-\frac{\Delta\,t}{\Delta\,x}((H_z)_{i+\fr,j+\fr}^{n+\fr}-(H_z)_{i-\fr,j+\fr}^{n+\fr})\\
%&\hspace{70mm}\ds - \frac{\Delta\,t}{4}(1-\nu)
%( \sigma_{i+\fr}(H_{z})_{i+\fr,j+\fr}^{n+\fr}+\sigma_{i-\fr}(H_{z})_{i-\fr,j+\fr}^{n+\fr})
%\end{split}
%\\[5mm]
%%
%\begin{split}
%&\ds(1+\frac{\Delta\,t}{4}(1+\nu)\sigma_{i+\fr})(H_{z})_{i+\fr,j+\fr}^{n+\fr}
%=(1-\frac{\Delta\,t}{4}(1+\nu)\sigma_{i+\fr})(H_{z})_{i+\fr,j+\fr}^{n-\fr}\\
%&\hspace{55mm}-\frac{\Delta\,t}{\Delta\,x}((E_y)_{i+1,j+\fr}^{n}-(E_y)_{i,j+\fr}^{n})\\
%&\hspace{60mm}+\frac{\Delta\,t}{\Delta\,y}((E_x)_{i+\fr,j+1}^{n}-(E_x)_{i+\fr,j}^{n})\\
%&\hspace{60mm}
%\ds -\frac{\Delta\,t}{4}(1-\nu)
% (\sigma_{i+1}(E_y)_{i+1,j+\fr}^{n}+\sigma_{i}(E_y)_{i,j+\fr}^{n})
%\end{split}
%%
%\end{cases}
%\]
The schemes are implemented using time windows to save memory.

The harmoniously matched layers can be implemented in several ways that we compare.
The function $\sig(x)$ is  as above.
\begin{description}
\item[HML Version 1.  Global extrapolation.]
Compute the solution of \eqref{eq:yeesmart} with an absorption of $\sig$, $(E^1,H^1)$ and $2\sig$, $(E^2,H^2)$ over the whole time window. Then $E_{x,y}=2*E_{x,y}^1-E_{x,y}^2$ and $H_{z}=2*H_{z}^1-H_{z}^2$.\\

\item[HML Version 2.  Local extrapolation.]
%HMLV2
Compute \textit{at each time step} the solution of \eqref{eq:yeesmart} with an absorption of $\sig$, $(E^1,H^1)$ and $2\sig$, $(E^2,H^2)$ over the whole time interval. Then $E_{x,y}=2*E_{x,y}^1-E_{x,y}^2$ and $H_{z}=2*H_{z}^1-H_{z}^2$.
Save computation by taking advantage of the fact that the computation of $E_x$ does not involve the absorption parameter. At each time step,
\begin{enumerate}
\item $E_x$ is computed by \eqref{eq:yeesmart_a},
\item two values of $E_y$ are computed by \eqref{eq:yeesmart_b}: $E_y^1$ with an absorption parameter equal to $\sig$, $E_y^2$ with an absorption parameter equal to $2\sig$.
\item two values of $H_z$ are computed by \eqref{eq:yeesmart_c}: $H_z^1$ with an absorption parameter equal to $\sig$, $H_z^2$ with an absorption parameter equal to $2\sig$.
\end{enumerate}
    Then $E_y=2*E_y^1-E_y^2$ and $H_z=2*H_z^1-H_z^2$.\\

\item[HML Version 3.  Split field local extrapolation.]
% HMLV3
At each time step,
\begin{enumerate}
\item $E_x$ is computed by \eqref{eq:yeesmart_a},
\item two values of $E_y$ are computed by \eqref{eq:yeesmart_b}: $E_y^1$ with an absorption parameter equal to $\sig$, $E_y^2$ with an absorption parameter equal to $2\sig$. Then $E_y=2*E_y^1-E_y^2$.
\item two values of $H_z$ are computed by \eqref{eq:yeesmart_c}: $H_z^1$ with an absorption parameter equal to $\sig$, $H_z^2$ with an absorption parameter equal to $2\sig$. Then $H_z=2*H_z^1-H_z^2$.
\end{enumerate}

\end{description}

We perform a series of experiments to illustrate the transmission properties of the layers.
The coefficient $\nu$ is meant to achieve  backward absorption and is taken equal to zero.
The domain of interest is $(0,6)\times (0,10)$, the
coefficient $\sigma(x)$ is supported in  $6\le x\le 10$.
The time of computation is 4, the initial electric field is zero.
 The initial transverse magnetic field,
 \[
H_z^0\ =\  \cos^2(\pi\frac{\left|\vec{x}-\vec{x}_c\right|}{r})
\
\cos(k\pi\,\vec{v\cdot}\frac{\vec{x}-\vec{x}_c}{r})\
\chi_{\left|\vec{x}-\vec{x}_c\right| \le r}
\]
is compactly supported in the ball $B(\vec{x}_c,r)$, with $\vec{x}_c=(5,5)$ and $r=0.8$.

The time of computation is fixed such that there is no reflection on the exterior walls.
The initial mesh is taken to be $\Delta x = \Delta y =0.1$, $ \Delta t=0.0702$, and then divided by $2$ twice.

In the first set of experiments, the absorption coefficient is constant in the layer, equal to $2$. The initial magnetic field hits the layer   at incidence $0^\circ$( $\vec{v}=(1,0)$) or $45^\circ$  ( $\vec{v}=(1,1)$).

In Table \ref{table:hfdisc} we compare the performances on a high frequency wave ($k=10$), while in Table \ref{table:bfdisc} we consider a low frequency wave ($k=1$).

In Table \ref{table:hfdisc} we compare the performances on a high frequency wave ($k=10$), while in Table \ref{table:bfdisc} we consider a low frequency wave ($k=1$).

\begin{table}[H]
\tbl{Comparison of the $L^\infty$ errors for high frequency, discontinuous absorption.   \label{table:hfdisc}
}
{\centering\tabcolsep0.3em
\begin {tabular}{|c||c|c|c||c|c|c|}\hline
& \multicolumn{3}{c|}{normal incidence} & \multicolumn{3}{c|}{$45^\circ$ incidence} \\ \hline
refinement& 0 & 1 & 2 & 0 & 1 & 2 \\\hline\hline
  \textbf{B\'erenger}  & \textbf{9.4e-02} & \textbf{3.9e-02}  &  \textbf{7.9e-03}
            & \textbf{1.3e-01 } &  \textbf{2.9e-02}  & \textbf{ 5.6e-03}\\
  Smart       & 5.2e-02   & 1.3e-02  &  5.1e-04
            & 6.2e-02  &  1.1e-02  &  5.3e-03\\
  HMLV1       & 3.4e-02  &  3.1e-03  &  2.1e-05
            & 4.5e-02  &  1.2e-03  &  5.5e-04\\
  HMLV2       & 2.5e-02   & 6.0e-03  &  1.2e-03
            & 7.4e-02  &  1.1e-02  &  1.7e-03\\
 \textbf{ HMLV3}       & \textbf{2.1e-02 } &  \textbf{4.2e-03}  &  \textbf{ 5.1e-04}
            &\textbf{4.5e-02 } &  \textbf{5.3e-03}  &  \textbf{5.7e-04}\\
  \hline
\end{tabular}}
\end{table}

\begin{table}[H]
\tbl{Comparison of the $L^\infty$ errors for low frequency, discontinuous absorption.  \label{table:bfdisc}
}
{\centering\tabcolsep0.3em
\begin {tabular}{|c||c|c|c||c|c|c|}\hline
& \multicolumn{3}{c|}{normal incidence} & \multicolumn{3}{c|}{$45^\circ$ incidence} \\ \hline
refinement& 0 & 1 & 2 & 0 & 1 & 2 \\\hline\hline
  \textbf{B\'erenger}  & \textbf{1.5e-02}  &  \textbf{7.1e-03}  &  \textbf{3.5e-03}
            & \textbf{1.3e-02}  &  \textbf{6.1e-03}  &  \textbf{3.0e-03}\\
  Smart       & 2.0e-02  &  2.0e-02  &  2.01e-02
            & 4.3e-02  &  4.2e-02  &  4.2e-02\\
  HMLV1       & 1.7e-02  &  1.60e-02  &  1.6e-02
            & 3.4e-02  &  3.3e-02  &  3.2e-02\\
  HMLV2       & 1.8e-02  &  1.1e-02  &  6.7e-03
            & 3.1e-02  &  1.9e-02  &  1.1e-02\\
  \textbf{HMLV3}       &\textbf{4.3e-03}  &  \textbf{2.6e-03}  &  \textbf{1.4e-03}
            &\textbf{8.2e-03 } &  \textbf{4.8e-03 } &  \textbf{2.6e-03}\\
  \hline
\end{tabular}}
\end{table}

In Tables \ref{table:hfcont} and \ref{table:bfcont}, we perform the same set of experiments, but the absorption coefficient is now a third degree polynomial in the layer, equal to $(x-6)^3/8$.

\begin{table}[H]
\tbl{Comparison of the $L^\infty$ errors for high frequency, continuous absorption.  \label{table:hfcont}
}
{\centering%\tabcolsep0.3em
\begin {tabular}{|c||c|c|c||c|c|c|}\hline
& \multicolumn{3}{c|}{normal incidence} & \multicolumn{3}{c|}{$45^\circ$ incidence} \\ \hline
refinement& 0 & 1 & 2 & 0 & 1 & 2 \\\hline\hline
  \textbf{B\'erenger}  & \textbf{3.8e-05}  &  \textbf{1.9e-07}  &  \textbf{2.1e-09}
            & \textbf{2.0e-04}  &  \textbf{9.1e-07}  &  \textbf{1.6e-09}\\
  Smart       & 2.7e-05  &  2.2e-07  &  1.7e-07
            & 1.7e-04  &  9.0e-07  &  3.1e-08\\
  HMLV1       & 5.5e-07  &  6.0e-08  &  5.6e-08
            &5.6e-06  &  1.2e-08  &  4.7e-09\\
  HMLV2       &6.8e-07  &  6.5e-08  &  3.1e-08
            & 2.6e-06  &  8.1e-09  &  2.8e-09\\
  \textbf{HMLV3}       & \textbf{5.8e-08}  &  \textbf{2.4e-09}  &  \textbf{1.1e-09}
            &\textbf{1.5e-06}  &  \textbf{9.5e-10}  &  \textbf{9.0e-11}\\
  \hline
\end{tabular}}
\end{table}

\begin{table}[H]
\tbl{Comparison of the $L^\infty$ errors for low frequency, continuous absorption.   \label{table:bfcont}
}
{\centering\tabcolsep0.3em
\begin {tabular}{|c||c|c|c||c|c|c|}\hline
& \multicolumn{3}{c|}{normal incidence} & \multicolumn{3}{c|}{$45^\circ$ incidence} \\ \hline
refinement& 0 & 1 & 2 & 0 & 1 & 2 \\\hline\hline
  \textbf{B\'erenger}  & \textbf{6.2e-07}  &  \textbf{3.2e-08}  &  \textbf{7.8e-010}
            &\textbf{5.2e-07}  &  \textbf{2.9e-08}  &  \textbf{6.5e-010}\\
  Smart       & 5.3e-04  &  5.3e-04  &  5.2e-04
            & 3.9e-04  &  3.8e-04  &  3.7e-04\\
  HMLV1       & 1.6e-04  &  1.6e-04  &  1.5e-04
            &8.6e-05  &  8.3e-05  &  8.2e-05\\
  HMLV2       &4.1e-04  &  2.0e-04  &  9.6e-05
            & 2.0e-04  &  9.8e-05  &  4.8e-05\\
  \textbf{HMLV3}       & \textbf{1.1e-05}  &  \textbf{5.4e-06}  &  \textbf{2.7e-06}
            & \textbf{5.9e-06}  &  \textbf{2.9e-06}  &  \textbf{1.4e-06}\\
  \hline
\end{tabular}}
\end{table}

The B\'erenger layer performs well on every frequency and every angle of incidence. Among the 3 versions for the HML, the third version is the best, which should be analyzed thoroughly.

Next compare the method on a gaussian initial value, supported in $(0,6)\times (0,10)$. Table \ref{table:gaussdiscont} uses a constant absorption in the layer, while Table \ref{table:gausscont} uses the same smooth absorption as before.
\begin{table}[H]
\tbl{Comparison of the $L^\infty$ errors for a gaussian initial magnetic field, constant absorption.\label{table:gaussdiscont}}
{\centering
\begin{tabular}{|c||c|c|c|}
\hline
refinement   &0 &1&2\\
\hline\hline
\textbf{B\'erenger} & \textbf{1.5e-02}  &  \textbf{6.7e-03}  &  \textbf{ 3.3 e-03}
\\
\hline
Smart      &3.4 e-02 & 3.4e-02 & 3.3e-02\\
\hline
HML V1      &3.0e-02 & 2.9e-02 & 2.8e-02
\\
\hline
HML V2     &3.6e-02 & 2.5e-02 & 1.6e-02
\\
\hline
\textbf{HML V3 }    &\textbf{1.0e-02 }&   \textbf{6.6e-03 } &  \textbf{3.9e-03}
\\
\hline
\end{tabular}}
\end{table}

\begin{table}[H]
\tbl{Comparison of the $L^\infty$ errors for a gaussian initial magnetic field, continuous absorption.\label{table:gausscont}}
{\centering
\begin{tabular}{|c||c|c|c|}
\hline
refinement   &0 &1&2\\
\hline\hline
\textbf{B\'erenger} & \textbf{7.5e-07}  &  \textbf{2.0e-08}  & \textbf{8.3e-10}
\\
\hline
Smart      & 4.3e-04 & 4.2e-04 & 4.1e-04\\
\hline
HMLV1
    &1.3e-04 & 1.2e-04 & 1.2e-04
\\
\hline
HMLV2     &3.0e-04 & 1.5e-04 & 7.3e-05
\\
\hline
\textbf{HMLV3 }   &\textbf{8.8e-06}  &  \textbf{4.3e-06  } &  \textbf{2.1e-06}
\\
\hline
\end{tabular}}
\end{table}

Finally, take unstructured random initial value,
supported in the ball centered at $(5,5)$ and of radius $1$. In Table \ref{table:randdiscont}, the absorption coefficient is constant in the layer, equal to $3$.

\begin{table}[H]
\tbl{Comparison of the $L^\infty$ errors for a random initial magnetic field, constant absorption.\label{table:randdiscont}}
{\centering
\begin{tabular}{|c||c|c|c|}
\hline
refinement   &0 &1&2\\
\hline\hline
\textbf{B\'erenger} & \textbf{5.7e-02}  &  \textbf{4.9e-02}  &  \textbf{ 4.4e-02}
\\
\hline
Smart      &6.7 e-02 &  6.3e-02  &  5.4e-02\\
\hline
HML V1      &5.1 e-02 &  4.5e-02 &   4.0e-02
\\
\hline
HML V2     &6.4e-02  &  3.0e-02 &  1.9e-02
\\
\hline
\textbf{HML V3 }    &\textbf{3.2e-02 }&   \textbf{1.5e-02 } &  \textbf{6.7e-03}
\\
\hline
\end{tabular}}
\end{table}

In Table \ref{table:randcont}, the absorption coefficient is a function of $x$ in the layer, equal to $(x-6)^3/8$.

\begin{table}[H]
\tbl{Comparison of the $L^\infty$ errors for a random initial magnetic field, continuous absorption.\label{table:randcont}}
{\centering
\begin{tabular}{|c||c|c|c|}
\hline
refinement   &0 &1&2\\
\hline\hline
\textbf{B\'erenger} & \textbf{1.1e-04}  &  \textbf{5.0e-05}  & \textbf{ 4.4e-06}
\\
\hline
Smart      &7.2e-04    &6.9e-04 &   6.4e-04\\
\hline
HMLV1
    &2.1e-04  &  2.2e-04 &   2.0e-04
\\
\hline
HMLV2     &5.0e-04  &  2.7e-04  &  1.2e-04
\\
\hline
\textbf{HMLV3 }   &\textbf{1.5e-05}  &  \textbf{7.9e-06 } &  \textbf{3.7e-06}
\\
\hline
\end{tabular}}
\end{table}

{\bf Summary.}  When comparing the reflection properties, the harmoniously matched layer, version 3, is
competitive with the B\'erenger layer.  For very regular data, the B\'erenger layers outperform everything. The performance of the HMLV3 gives hope the method with its   stronger  well posedness,
more robust absorption,  and small reflection at all
angles will be a good method where B\'erenger has proven less
good.  For example,
  for non constant coefficients and nonlinear problems.
  We have taken pains to make the comparison where
  B\'erenger is at its best.
  In 2D with a layer in a single direction the HML has an extra cost.  Since there are 5 quantities to compute at each time step instead of 4 for B\'erenger. This  is no longer  the case in
  three dimensions, since both strategies have to split 6 unknowns.

  \vskip.2cm

\noindent
{\bf Open problems.}  {\bf 1.}  Our  analysis does not explain the much better behavior with continuous absorption, nor the
advantages of HMLV3.
  {\bf 2.}  A comparison with other methods where only supplementary ordinary differential equations are added should be made.

\section*{Acknowledgments}This research project has spanned many years.  It ows a great deal
to the support of the University Paris 13 where J. Rauch was often invited for one month visits.

%%%%%%%%%%%%%%%%%%%%%%%
%% New section
%%%%%%%%%%%%%%%%%%%%%%%%%%
%\section*{References}
%\input{biblio.tex}
\bibliographystyle{plain}
\bibliography{pml}

 \end{document}